\documentclass[12pt,a4paper,intlimits,sumlimits]{amsart}

\usepackage{amsmath, mathtools}
\usepackage{amsfonts}
\usepackage[alphabetic]{amsrefs}
\usepackage{graphicx}
\usepackage{framed}
\usepackage{amssymb}
\usepackage{esvect}
\usepackage{xcolor}
\usepackage{eucal}

\usepackage{hyperref}

\usepackage{etex}
\usepackage{latexsym, mathtools}
\usepackage[english]{babel}
\usepackage[utf8x]{inputenc}

\def \N {\mathbb{N}}
\def \R {\mathbb{R}}
\def \Z {\mathbb{Z}}

\usepackage{amsthm}
\theoremstyle{definition}
\newtheorem{definition}{Definition}[section]

\theoremstyle{plain}
\newtheorem{theorem}[definition]{Theorem}

\newtheorem{conjecture}[definition]{Conjecture}

\numberwithin{equation}{section}
\usepackage{geometry}
\renewcommand{\epsilon}{\varepsilon}
\newcommand{\e}{\varepsilon}

\renewcommand{\leq}{\leqslant}
\renewcommand{\le}{\leqslant}
\renewcommand{\geq}{\geqslant}
\renewcommand{\ge}{\geqslant}

\newcommand{\Per}{\operatorname{Per}}
\renewcommand{\div}{\operatorname{div}}

\allowdisplaybreaks

\begin{document}

\begin{abstract} We present here some classical and modern results about
phase transitions and minimal surfaces, which are quite intertwined topics.

We start from scratch, revisiting the theory of phase transitions as put forth by Lev Landau.
Then, we relate the short-range phase transitions to the classical minimal surfaces,
whose basic regularity theory is presented, also in connection with a celebrated conjecture by Ennio De Giorgi.

With this, we explore the recently developed subject of long-range phase transitions
and relate its genuinely nonlocal regime to the analysis of fractional minimal surfaces.
 \end{abstract}
 
 \title[(Non)local phase transitions and minimal surfaces]{Some perspectives \\
 on (non)local phase transitions and minimal surfaces}
 
 \author{Serena Dipierro}
 \author{Enrico Valdinoci}
 
 \address{Department of Mathematics and Statistics
 \newline\indent University of Western Australia \newline\indent
 35 Stirling Highway, WA 6009 Crawley, Australia}
 \email{serena.dipierro@uwa.edu.au, enrico.valdinoci@uwa.edu.au}
 
 \date{}
 
 \maketitle\tableofcontents
 
\section{Introduction} 

In this paper we deal with phase coexistence models, in connection with
the geometric problems of minimal interfaces. We will analyze the classical case,
in which the limit interface is regulated by surface tension, as well as the long-range case, which,
at least in some circumstances, can be related to a variational problem of fractional perimeter type.
To fully understand the rather sophisticated theory for these problems and its concrete impact on applied sciences, and also to develop a sufficient intuition of the elaborated questions arising in this context,
we will also review the basic physical notions related with phase transition models.\medskip

A classical description of a system exhibiting two phases is given by the energy functional
\begin{equation}\label{WDW-CLASSi}
\int_\Omega\left(\frac{|\nabla u(x)|^2}2+W(u(x))\right)\,dx,\end{equation}
where~$W$ is a ``double-well potential'', whose minima correspond to the two phases under consideration.

Roughly speaking, a double-well potential is a (sufficiently smooth) function
with two absolute minima. For convenience, one can renormalize~$W$ to be nonnegative and
its minima to be located at~$-1$ and~$+1$.
For stronger geometric results, it is also customary to assume that these minima are nondegenerate
(though some geometric investigation on degenerate cases are also possible, see~\cite{MR3890785}).
In this way, we will suppose that
\begin{equation}\label{DOPP}\begin{split}
&{\mbox{$W\in C^2(\R)$,
$W(x)>0$ for all~$x\in\R\setminus\{-1,1\}$,}}\\
&{\mbox{$W(-1)=W(1)=0$
and~$W''(\pm1)>0$.}}\end{split}\end{equation}
A typical example of double-well potential is given by the symmetric quartic polynomial
\begin{equation}\label{WDW}
W(u)=\frac{(1-u^2)^2}4.\end{equation}

As discussed in quite detail in Sections~\ref{SEC:ORD:THE}
and~\ref{FIR:ORD:THE},
the functional in~\eqref{WDW-CLASSi}
is related to the so-called free energy of a material exhibiting the coexistence of two phases and the model is general enough to comprise both\footnote{The names of first- and second-order phase transitions are nowadays quite standard, though they may be a bit confusing, since they do not deal with the order of a partial differential equations, nor with the degree of a polynomial.
Also, usually, what is ``second-order'' is potentially more complicated than what is ``first-order'', but for phase transitions the opposite
feature occurs: first-order phase transitions are potentially more
complicated than second-order ones, since they show the additional
phenomenon of ``latent heat''. Since second-order phase transitions
are technically simpler, we discuss them in Section~\ref{SEC:ORD:THE},
namely before, not after, the first-order phase transitions,
which are presented in Section~\ref{FIR:ORD:THE}.}
``first-order phase transitions'' at the critical temperature (in which the production or absorption of a ``latent heat'' occurs during the phase transition
without changing the temperature of the medium) and ``second-order phase transitions'' at a temperature less than or equal to the critical one (in this case, the phase change occurs without ``latent heat'').

We will also consider a natural nonlocal counterpart of~\eqref{WDW-CLASSi}, given by
\begin{equation}\label{PRES} \frac14 \iint_{{\mathcal{Q}}(\Omega)} \frac{(u(x)-u(y))^2}{|x-y|^{n+\alpha}}\,dx\,dy+\int_\Omega W(u(x))\,dx,\end{equation}
with
$$ {\mathcal{Q}}(\Omega):=(\Omega\times\Omega)\cup(\Omega^c\times\Omega)\cup
(\Omega\times\Omega^c),$$
being~$\Omega^c:=\R^n\setminus\Omega$ (a more precise intuition of the set~$ {\mathcal{Q}}(\Omega)$ as the set collecting all the possible interactions that affect the states
in the container~$\Omega$ will be discussed in Section~\ref{ISS}).

In~\eqref{PRES}, the exponent~$\alpha$ lies in the range~$(0,2)$, but, as we will see more precisely
in Section~\ref{GAMMACO2}, the value~$\alpha=1$ provides a structural threshold:
specifically, we will see that at a large scale the model tends to reduce to the classical interface
problem when~$\alpha\in[1,2)$, while when~$\alpha\in(0,1)$ the nonlocal features persist at every
scale and the asymptotic behavior is related to minimal surfaces of nonlocal type.

To describe these nonlocal minimal surfaces, as introduced by~\cite{MR2675483},
given~$\alpha\in(0,1)$ and two (measurable) disjoints subsets~$E$ and~$F$ of~$\R^n$, one considers
the~$\alpha$-interaction between~$E$ and~$F$, defined by
\begin{equation}\label{per-I} {\mathcal{I}}_\alpha(E,F):=\iint_{E\times F}\frac{dx\,dy}{|x-y|^{n+\alpha}}.\end{equation}
This is the cornerstone to construct the $\alpha$-perimeter functional, as in~\cite{MR2675483}. Namely,
one sets
\begin{equation}\label{per-IO} {\rm Per}_\alpha(E):={\mathcal{I}}_\alpha(E,E^c).\end{equation}
More generally,
given a bounded reference domain~$\Omega$ with Lipschitz boundary, one defines the $\alpha$-perimeter of~$E$
in~$\Omega$ as the collection of all the interactions between
the set~$E$ and its complement in which at least one side of the interaction lies in~$\Omega$,
that is
\begin{equation}\label{per-a} 
{\rm Per}_\alpha(E,\Omega):=
{\mathcal{I}}_\alpha(E\cap\Omega,E^c\cap\Omega)+{\mathcal{I}}_\alpha(E\cap\Omega^c,E^c\cap\Omega)+
{\mathcal{I}}_\alpha(E\cap\Omega,E^c\cap\Omega^c).\end{equation}
\medskip

With this basic mathematical setting in mind, we can now dive into a detailed analysis of the models above,
in connection with classical and contemporary results aiming at understanding qualitative and
symmetry properties for phase transitions. To this end, we will start by recalling in Section~\ref{LANDATHE}
the classical theory of phase transitions introduced by Lev Landau: the exposition aims to be detailed
and accessible, and basically no prerequisite is assumed (up to some basic thermodynamics
combined with common sense physical intuition).

Then, 
Section~\ref{Skdfe0rthX234rt} showcases one of the most typical mathematical
formulations for phase coexistence problems, namely the
Allen-Cahn equation, which is analyzed in Section~\ref{GAMMACO}
in the light of~$\Gamma$-convergence.

This constitutes a strong link between phase transitions and
minimal surfaces, i.e. hypersurfaces locally minimizing a perimeter functional.
The regularity of these surfaces will be recalled in Section~\ref{1243546578679ow3etg245tPOkcv},
thus motivating a classical conjecture by Ennio De Giorgi in relation with
the one-dimensional symmetry of global and monotone phase transitions.

Sections~\ref{x-09i8uytf-0iuytfdfghjdP}, \ref{GAMMACO2} and~\ref{09i2w-e23rt5-5PKM78}
provide the long-range counterpart of the previous investigations
from Sections~\ref{Skdfe0rthX234rt}, \ref{GAMMACO} and~\ref{1243546578679ow3etg245tPOkcv}. Namely,
Section~\ref{x-09i8uytf-0iuytfdfghjdP}
clarifies the role
of long-range interactions in the theory of phase transitions
and presents the nonlocal Allen-Cahn equation, while Section~\ref{GAMMACO2}
presents a recently established~$\Gamma$-convergence theory for the long-range setting,
and then Section~\ref{09i2w-e23rt5-5PKM78}
deals with
nonlocal minimal surfaces and with the long-range version
of the above mentioned conjecture by De Giorgi.

The paper ends with some technical estimates collected in Appendices~\ref{appeadd10}
and~\ref{sec:jyhtgrfu7y6tg5rfet8569b7c5tnxm4utxc46}.

\section{Landau theory of phase transitions}\label{LANDATHE}

The description of coexisting phases in a given material is a very complex topic,
combining different perspectives from statistical physics, thermodynamics, material sciences
and mathematical analysis. 
Also, the setting may vary according to the specific case under consideration, given the
different underlying physical structures involved: for instance,
while in our everyday life we are mostly exposed to the changes of state
related to melting, freezing, vaporization and condensation,
other important phase transitions, such as the one from a conducting to
a superconducting state, are the outcome of brand new properties of the solid state, such as electron coupling.

Though it is virtually impossible to provide here an exhaustive account of phase transition
theories, we can recall some general facts which can be helpful to develop an intuition of the problem under consideration.

A common treat in the study of phase transitions is to describe a given system
in terms of relevant physical quantities, such as temperature~$T$,
pressure~$P$, volume~$V$, entropy~$S$, magnetic moment~$M$, etc.
The energy of a system is thus an outline of its
``ability to perform some tasks''.
In this performance, however,
the system 
typically ``wastes''
some energy in the form of heat.
The ``useful'' energy, that is the energy that is available to do work, thus consists in the difference between
the full internal energy of the system minus the energy that is ``unavailable to perform work''
since it gets lost through heat.
Being ``free to do the work'', such energy is often called ``free energy'', though\footnote{Our presentation of free energy is certainly quite inaccurate. \label{FOFREEN} More precisely, there are typically two versions of this quantity, depending on the physical variables that describe the state of the system.
On the one hand, when one models the system by using as main independent variables the temperature~$T$ and the volume~$V$, one obtains the Helmholtz free energy~$F=F(V,T)$, whose infinitesimal increment is defined by
$$ dF = -S\,dT - P\, dV.$$
On the other hand, if the system is described in terms of the independent variables pressure~$P$ and temperature~$T$, one obtains the Gibbs free energy~$G=G(P,T)$, whose infinitesimal increment is given by
$$ dG = -S\,dT + V \,dP.$$
So, some care is needed when dealing with free energy, and the precise dependence on the relevant physical values cannot be, in general, omitted.
Nonetheless, in our discussion here, we are considering the free energy in dependence of temperature only, therefore the two approaches, namely the one based on the Helmholtz free energy and that based on the Gibbs free energy, are essentially equivalent. In particular, from the infinitesimal increments, we see that when pressure and volumes are kept constant and only temperature varies, then
$$ S=-\frac{\partial F}{\partial T}=-\frac{\partial G}{\partial T},$$
that is, one reconstructs the entropy from the variation in temperature of the free energies (essentially, in an equivalent way with respect to the choice of free energy).

We should also stress that the interpretation of free energy as the part of internal energy which is keen to perform work is also rather simplistic. Indeed, looking at the increments defining the free energies, one sees that they comprise the available energy to do mechanical work (related to change of pressure and volumes) but also a term which is entropy-dependent and temperature-related of the form~$-S\,dT$, which is not directly related to mechanical work. The term free is likely to be related to the fact that the increments of the free energies do not present terms of the form~$T\,dS$, which, in reversible processes, would correspond to heat increments~$dQ$, via the Second Law of Thermodynamics (that is, very roughly speaking, the contributions of heat are removed from the free energy).}
the name is under an intense debate.

Hence, the minimizers (or, more generally, the stable critical points) of this energy correspond to observable states
of our system. Concretely, the system may present significantly different features, or ``phases'',
such as being in a solid or fluid state, or having a magnetic momentum, or presenting superfluid or superconductor properties,
and the appearance of these phases may be seen as an outcome of energy minimization.

The arising of different phases may be the outcome of a critical physical parameter involved in the free
energy, such as temperature: in practice,
a ``disordered phase'' typically corresponds to a high temperature, while an ``exceptionally ordered phase''
arises at low temperatures. This is the case, for instance,
for magnetization, since magnetic materials have no permanent magnetic moment above their\footnote{The Curie Temperature is named after Pierre Curie, who in 1895 related some magnetic properties to change in temperature.}
Curie Temperature (about $570$ degree Celsius for the usual magnetite)
but below this temperature 
the atoms tend to behave as tiny magnets which spontaneously align themselves,
so that the magnetic materials show a permanent magnetization oriented in a certain direction, see e.g.~\cite{zbMATH05046576}.

Specifically, when the temperature~$T$ is above the Curie Temperature~$T_c$,
in the absence of external sources a magnetic system lies in a zero field state,
which happens to be a minimal configuration for the free energy corresponding to its temperature.
When the temperature is decreased below such critical temperature~$T_c$,
the system will go through a state in which the magnetization is still zero,
but this corresponds only to a critical point,
not a minimum of the free energy, making this equilibrium configuration totally unstable.
For this reason, below the critical temperature~$T_c$,
small perturbations from the environment will inevitably
induce the system to reorganize its microscopic structure to possibly preserve a zero average but creating
regions with a nontrivial magnetic momentum.
The formation of these magnetic
domains will produce a supplementary interfacial energy, which, in some sense, favors domain segregation
with a phase separation which is ``as small as possible'' (in a sense which will be clarified below
in Section~\ref{ISS}).

It is also interesting to remark that such a phase separation also produces a symmetry breaking:
the free energy is symmetric (since it weighs equally, say, the magnetizations
oriented towards the North pole and the ones oriented towards the South pole), 
nonetheless the magnetization configuration reached by the system during the cooling
is somewhat accidental, as a result of small environmental perturbations, making
the final state reached by the system not necessarily symmetric.

\subsection{The second-order theory of phase transitions}\label{SEC:ORD:THE}
To account for the phenomena of phase transition and phase coexistence, one can consider an
``order parameter''~$\eta$ which describes, in some sense, how every point of the system is ``organized''
with respect to a given notion of phase. The name of this parameter is possibly inspired by its statistical mechanics
interpretation, relating different phases to the degree of ``order'', or ``disorder'', of a system.
In specific situations, the order parameter can be either a scalar or a vector:
for instance,
in a liquid-gas phase transition the order parameter corresponds to the difference of the densities
between the two phases, which is a scalar, while in
superfluidity and superconductivity it is a complex-valued wave function 
(or, equivalently, a two-dimensional vector), and for ferromagnetic
momenta it is in general a three-dimensional vector. Phase transitions also occur
in cosmology, since as the universe expanded and cooled,
a number of symmetry-breaking phase transitions occurred, and the description of these phenomena often
relies on an order parameter which is a tensor.
Here we only consider the case in which the order parameter~$\eta$ is a scalar function.

As in the above discussion about the Curie Temperature,
the phase transition theory introduced in~\cite{LANDAUPH} considers a critical situation in which
phase separation occurs. Assume that
the free energy
presents the Taylor expansion
\begin{equation} \label{CEFF2}
a_0 +a_1\eta+a_2\eta^2+a_3 \eta^3+a_4 \eta^4+\dots,\end{equation}
where the coefficients~$a_0$, $a_1$, $a_2$, $a_3$, $a_4$, $\dots$, depend on relevant physical quantities.
Here, for the sake of concreteness, we suppose that they depend on the temperature~$T$ of the system.
We also focus on the case in which the free energy is symmetric with respect to the order parameter
(say, assuming that the deviations from the neutral case equally affect the energy, as in the magnetic
case outlined above). In this situation, the odd coefficients of the expansion in~\eqref{CEFF2} must vanish,
thus reducing the free energy to
\begin{equation} \label{CEFF}
a_0 +a_2\eta^2+a_4 \eta^4,\end{equation}
where we have also neglected the higher order terms.
The coefficient~$a_0$ is irrelevant to determine the critical points of this energy, but the coefficients~$a_2$
and~$a_4$ play an essential role. So, we can think that~$a_0$ is just a constant, while~$a_2=a_2(T)$ and~$a_4=a_4(T)$
depend on the temperature~$T$.
In particular, the sign of these coefficients is determinant in the formation of new phases.
Indeed, if we expect the state parameter~$\eta$ to be confined in a bounded region
(which is typically the case, since we do not expect that a physical parameter tends to diverge),
it is convenient to assume that~$a_4(T)$ is positive for all~$T$ (in this way, the free energy~\eqref{CEFF} 
is bounded from below and possesses minima for all~$T$).
To model the spontaneous formations of new phases below the critical temperature, one may suppose that
\begin{equation}\label{LFL}
{\mbox{$a_2(T)>0$
for all~$T>T_c$ and~$a_2(T)<0$ for all~$T<T_c$.}}\end{equation}
Also, we assume that~$a_2$ varies continuously with respect to~$T$, hence
\begin{equation}\label{24BIS} a_2(T_c)=0.\end{equation}
In this way, one readily checks that the critical points of~\eqref{CEFF}
are
\begin{equation}\label{ANA-01} \begin{dcases}
\{0\} & {\mbox{ if }} T\ge T_c,\\
\\
\displaystyle\left\{-\sqrt{-\frac{a_2(T)}{2a_4(T)}},\,0,\,\sqrt{-\frac{a_2(T)}{2a_4(T)}}
\right\} & {\mbox{ if }} T< T_c.
\end{dcases} \end{equation}
Furthermore, 
\begin{equation}\label{ANA-02}\begin{split}&
{\mbox{the critical point~$0$ is a nondegenerate minimum when~$T>T_c$,}}\\
&{\mbox{a degenerate minimum when~$T=T_c$,}}\\
&{\mbox{and a local maximum when~$T<T_c$,}}\\
&{\mbox{while~$\displaystyle\pm\sqrt{-\frac{a_2(T)}{2a_4(T)}}$ are nondegenerate minima when~$T<T_c$.}}\end{split}\end{equation}

\begin{figure}[h]
\includegraphics[width=0.65\textwidth]{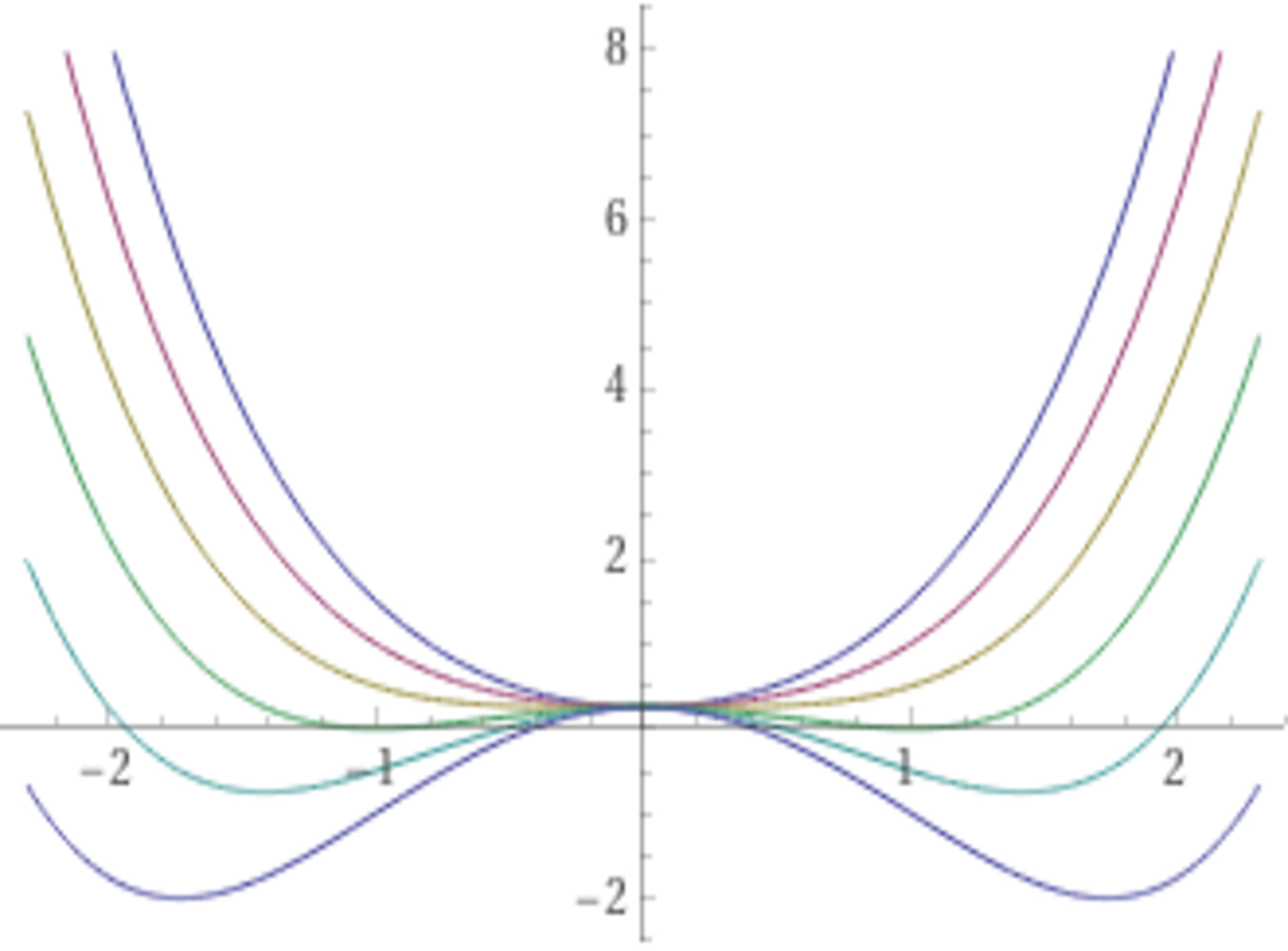}
\caption{Free energy corresponding to the parameters in~\eqref{AVPA}
for~$T\in\{4,3,2,1,0,-1\}$.}
        \label{9FIGURATBORDINP}
\end{figure}

\begin{figure}[h]
\includegraphics[width=0.85\textwidth]{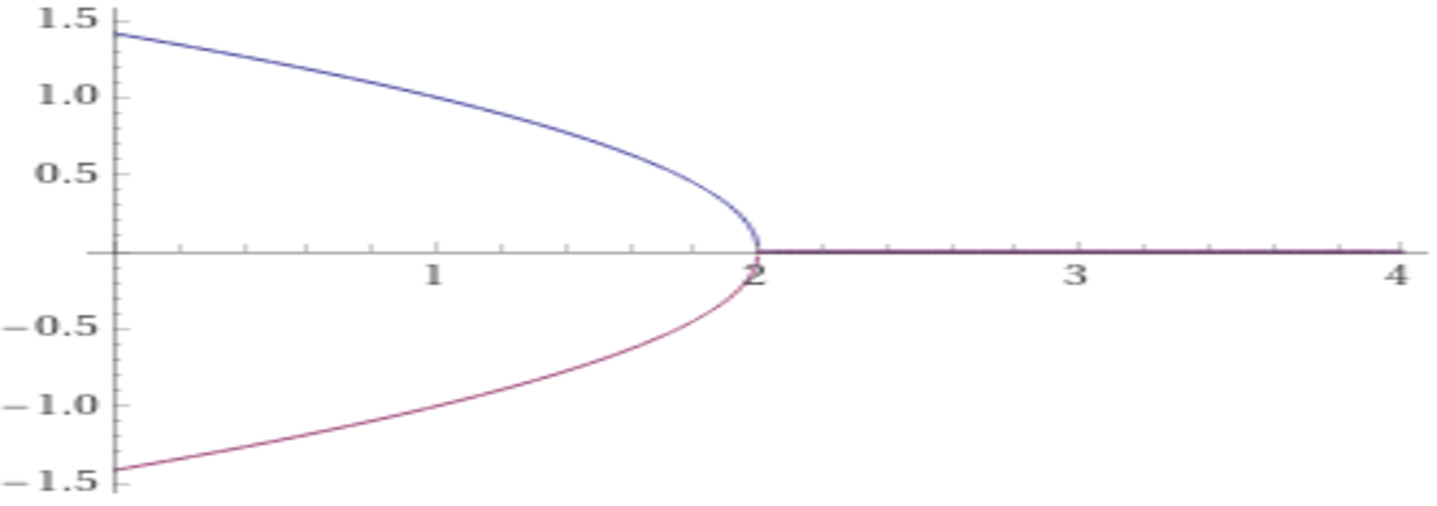}
\caption{Bifurcation diagram of the minimizers corresponding to the parameters in~\eqref{AVPA}
for~$T\in[0,4]$.}
        \label{9FIGURATBORDINP-1o}
\end{figure}

As an example, just to favor intuitive thinking, and without aiming at a realistic physical description of a specific material,
one can consider the case
\begin{equation}\label{AVPA}
T_c:=2,\qquad
a_0:=\frac14,\qquad a_2(T):=\frac{T-2}{2}
\qquad{\mbox{and}}\qquad a_4(T):=\frac14.\end{equation} Different plots of the free energy in~\eqref{CEFF}
for this case
related to various choices of the temperature~$T$ are given in Figure~\ref{9FIGURATBORDINP},
where the bifurcation diagram corresponding to the critical points\footnote{Obviously, deducing ``global'' properties of an energy functional, as in~\eqref{ANA-01} and~\eqref{ANA-02}, from its Taylor expansion in the vicinity of the origin, as in~\eqref{CEFF2}, is, to say the least, conceptually inadequate. Nevertheless, the local expansion
in~\eqref{CEFF2} here was mainly utilized to reduce the calculation to a four-degree polynomial, thus simplifying the notation and making the structure of the energy functional more transparent. Once one understands this simplified scenario for phase transitions, one also has a clue of the bifurcation occurring for the critical points of the energy in dependence of a variable parameter, thus recovering the situation in Figure~\ref{9FIGURATBORDINP}, and related ones, possibly in a more general, and more rigorous, way.}
is apparent.
Specifically, the minimizers of the free energy in the model case~\eqref{AVPA}
are
\begin{equation}\label{T-MINI} \begin{dcases}
\{0\} & {\mbox{ if }} T\ge 2,\\
\\
\displaystyle\left\{-\sqrt{2-T},\,\sqrt{2-T}
\right\} & {\mbox{ if }} T< 2,
\end{dcases} \end{equation}
see Figure~\ref{9FIGURATBORDINP-1o}.

Notice also that the parameters above produce the double-well potential in~\eqref{WDW}
when~$T:=1$.

The bifurcation diagram depicted in Figure~\ref{9FIGURATBORDINP}
represents a situation in which  the order parameter exhibited by the system
changes continuously with respect to the temperature: namely, in light of~\eqref{ANA-01},
if the free energy depends continuously on the temperature~$T$,
then the new minima when~$T<T_c$ can be seen as a continuous 
modification of the null state: for instance, in the setting of~\eqref{AVPA},
these minima are given by~$\pm\sqrt{2-T}$ for~$T<2$.
Interestingly, the dependence of these minima on the temperature~$T$ is not in general smooth,
due to the presence of the square root.

In jargon, the situations in which the observed order parameter depends continuously
on the temperature are called ``second-order phase transitions''.
The name possibly comes from this feature: if we evaluate the free energy
in~\eqref{CEFF} at the minima in~\eqref{ANA-01} we obtain
\begin{eqnarray*} {\mathcal{E}}(T)&:=&\begin{dcases}a_0 & \;{\mbox{ if }}T\ge T_c,\\
a_0-\frac{a_2^2(T)}{2a_4(T)}+\frac{a_4(T)\, a_2^2(T)}{4a_4^2(T)} &\; {\mbox{ if }}T<T_c,\end{dcases}\\
\\&=&
\begin{dcases}a_0 & \;{\mbox{ if }}T\ge T_c,\\
a_0-\frac{ a_2^2(T)}{4a_4(T)} &\; {\mbox{ if }}T<T_c.
\end{dcases}
\end{eqnarray*}
The function~${\mathcal{E}}$ is sometimes called ``the free energy as a function of temperature''.

Moreover,
$$ {\mathcal{E}}'(T)=\begin{dcases}0 & \;{\mbox{ if }}T> T_c,\\
-\frac{a_2(T)\,a_2'(T)}{2a_4(T)}
+\frac{ a_2^2(T)\,a_4'(T)}{4a_4^2(T)}
 &\; {\mbox{ if }}T<T_c,
\end{dcases}$$
whence, recalling~\eqref{24BIS},
$$ \lim_{T\nearrow T_c} {\mathcal{E}}'(T)=0=\lim_{T\searrow T_c} {\mathcal{E}}'(T).$$
Accordingly,
\begin{equation}\label{SECO01}
\begin{split}&
{\mbox{the first derivative with respect to temperature}}\\ &
{\mbox{of the ``free energy
as a function of temperature''}}\\&
{\mbox{vanishes continuously at the critical temperature.}}\end{split}
\end{equation}

Furthermore, 
$$ {\mathcal{E}}''(T)=\begin{dcases}0 & \;{\mbox{ if }}T> T_c,\\
\begin{matrix}\displaystyle
-\frac{ (a_2'(T))^2}{2a_4(T)}
-\frac{ a_2(T)\,a_2''(T)}{2a_4(T)}
+\frac{a_2(T)\,a_2'(T)\,a_4'(T)}{a_4^2(T)}\\ \displaystyle
+\frac{ a_2^2(T)\,a_4''(T)}{4a_4^2(T)}
-\frac{ a_2^2(T)\,(a_4'(T))^2}{2a_4^3(T)}\end{matrix}
 &\; {\mbox{ if }}T<T_c,
\end{dcases}$$
leading to
$$ \lim_{T\nearrow T_c} {\mathcal{E}}''(T)=
-\frac{ (a_2'(T_c))^2}{2a_4(T_c)}
\quad{\mbox{ while }}\quad
\lim_{T\searrow T_c} {\mathcal{E}}''(T)=0.$$
In particular, if~$T_c$ is a nondegenerate zero of~$a_2$ (as it happens for instance in
the model case~\eqref{AVPA}), we find that
$$ \lim_{T\nearrow T_c} {\mathcal{E}}''(T)<0=
\lim_{T\searrow T_c} {\mathcal{E}}''(T).$$
In this situation,
\begin{equation}\label{SECO02}
\begin{split}&
{\mbox{the second derivative with respect to temperature}}\\ &
{\mbox{of the ``free energy
as a function of temperature''}}\\&
{\mbox{presents a discontinuity at the critical temperature.}}\end{split}
\end{equation}

\begin{figure}[h]
\includegraphics[height=0.16\textheight]{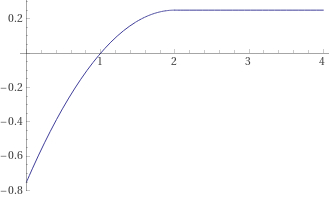}$\qquad$
\includegraphics[height=0.16\textheight]{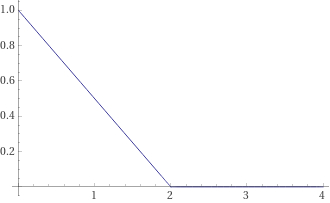}
\caption{The free energy
as a function of temperature, and its derivative,
corresponding to the parameters in~\eqref{AVPA}
for~$T\in[0,4]$.}
        \label{9FGHSJIGURATBORDINP}
\end{figure}

The phenomena in~\eqref{SECO01} and~\eqref{SECO02}
are likely to be the justification of the name of
``second-order phase transitions'' to describe these situations:
see Figure~\ref{9FGHSJIGURATBORDINP}
for a sketch of the functions~${\mathcal{E}}$ and~${\mathcal{E}}'$ in the model case~\eqref{AVPA}.

{F}rom a physical point of view, the discontinuities of~${\mathcal{E}}'$ at the critical temperature
are related\footnote{Without aiming at being thorough, a heuristic explanation of the link between
the ``latent heat'' \label{order-phase0}
and the possible discontinuities of~${\mathcal{E}}'$ can be obtained as follows. Consider a system undergoing a phase change at temperature~$T_c$, keeping the other physical parameters constant. Owing to the discussion in footnote~\ref{FOFREEN},
we know that the derivative with respect to temperature of the
free energy corresponds, up to a sign change, to entropy, hence we write~$S=-{\mathcal{E}}'$ and we use the Second Law of Thermodynamics~${dS}=\frac{dQ}T$ by supposing that a reversible process occurs
in the phase change due to the cooling of the system from temperature~$T_c+\e$ to~$T_c-\e$, being~$\e>0$ small.
In this way, we formally have that
$$ dQ=T\,dS=-T\,d{\mathcal{E}}',$$
though we need to interpret this equation with a pinch of salt, since~$d{\mathcal{E}}'={\mathcal{E}}''\,dT$ may be not classically defined at~$T_c$
(due to the possible discontinuities of~${\mathcal{E}}'$ at the critical temperature).
Therefore, it is convenient to integrate the above expression separately for~$T\in(T_c-\e,T_c)$ and~$T\in(T_c,T_c+\e)$. In this way, we find that, when~$\e\searrow0$,
\begin{eqnarray*} &&Q(T_c)-Q(T_c-\e)=\int_{T_c-\e}^{T_c} dQ=-\int_{T_c-\e}^{T_c}T\,d{\mathcal{E}}'=
-\int_{T_c-\e}^{T_c}(T_c+o(1))\,d{\mathcal{E}}'\\&&\qquad=-(T_c+o(1))
\big({\mathcal{E}}'(T_c)-{\mathcal{E}}'(T_c-\e)\big)
\end{eqnarray*}
and similarly
$$ Q(T_c+\e)-Q(T_c)=-(T_c+o(1))\big({\mathcal{E}}'(T_c+\e)-{\mathcal{E}}'(T_c)\big).$$
By summing up the latter two equations, one obtains that
$$ Q(T_c+\e)-Q(T_c-\e)=-(T_c+o(1))\big({\mathcal{E}}'(T_c+\e)-{\mathcal{E}}'(T_c-\e)\big).$$
Thus, assuming that the left- and right-limits as~$\e\searrow0$ exist,
$$ Q(T_c^+)-Q(T_c^-)=- T_c\big({\mathcal{E}}'(T_c^+)-{\mathcal{E}}'(T_c^-)\big).$$
This shows that a discontinuity of the derivative of the free energy at the critical temperature corresponds to
a ``latent heat'' proportional to minus such a discontinuity times the critical temperature.}
to the ``latent heat'' (roughly speaking,
the energy released or absorbed by the system in a phase change
without changing its temperature), hence~\eqref{SECO01}
corresponds to the absence of latent heat
in these types of phase transitions.

\subsection{The first-order theory of phase transitions}\label{FIR:ORD:THE}
In many physical situations, however, 
the change of the state of a substance at its critical temperature
is related to a latent heat which is supplied to or extracted from the medium
without changing its temperature. These types of phase transitions correspond to a discontinuity
of the first derivative of~${\mathcal{E}}$ and are called in jargon
``first-order phase transitions''.
In these situations, the observed
order parameters also jump discontinuously at
the transition temperature.

It is instructive to remark that a simple modification of the above
theory can also describe these phenomena.
For this, we retake the free energy expansion in~\eqref{CEFF2}, neglecting higher order terms,
and we aim at modeling a situation in which~$\eta=0$ is the observed value of the state parameter
above a critical temperature~$T_c$, but below~$T_c$ a new stable phase arises.
More precisely, we describe a model in which~$\eta=0$ is a nondegenerate
local minimum, hence a stable phase, for the free energy for
all values of the temperature~$T$
(and the only minimizer when~$T>T_c$), but a new stable phase arises when~$T\le T_c$, with the new phase
becoming a global minimizer when~$T<T_c$.
Here we are not assuming that the free energy is symmetric in~$\eta$.
{F}rom~\eqref{CEFF2}, the condition that~$\eta=0$ is a critical point gives that~$a_1$ must necessarily vanish for all~$T$,
therefore we can rewrite the free energy in this case as
\begin{equation}\label{GBSUJ-0} a_2(T)\,\eta^2+a_3(T) \eta^3+a_4(T) \eta^4,\end{equation}
where, for simplicity, we have dropped the term~$a_0$ since it does not modify the critical points of the system.

Accordingly, the condition that~$\eta=0$ is a nondegenerate
local minimum of the free energy yields that
\begin{equation}\label{BSNTT-1} a_2(T)>0.\end{equation}
The condition that the energy is bounded from below (thus producing minimizers) also gives that
\begin{equation}\label{BSNTT-2} a_4(T)>0.\end{equation}

The phase transition can be thus modeled on the specific properties of~$a_3(T)$.
Namely, the assumption that~$\eta=0$ is the only minimizer for~$T>T_c$ says that
\begin{equation}\label{BSNTT-3} a_2(T)+a_3(T) \eta+a_4(T) \eta^2>0\qquad{\mbox{for all~$\eta\in\R$ and~$T>T_c$.}}\end{equation}
Also, we assume that at~$T=T_c$ a new minimizer, say at~$\eta=\eta_\star>0$, occurs, whence
\begin{equation}\label{BSNTT-4} a_2(T_c)\,\eta_\star^2+a_3(T_c) \eta_\star^3+a_4(T_c) \eta^4_\star=0.\end{equation}
The existence of a global minimum different from~$\eta=0$ below the critical temperature translates into
\begin{equation}\label{BSNTT-5} \min_{\eta\in\R} a_2(T)\,\eta^2+a_3(T) \eta^3+a_4(T) \eta^4<0\qquad{\mbox{for all }}T<T_c.\end{equation}
An example of coefficients satisfying~\eqref{BSNTT-1}, \eqref{BSNTT-2}, \eqref{BSNTT-3}, \eqref{BSNTT-4}
and~\eqref{BSNTT-5} is, for instance,
\begin{equation}\label{BSNTT-6} \begin{split}&
T_c:=2,\qquad
a_2(T):=1,\\ &{\mbox{$a_3$ a smooth and increasing function such that~$a_3(T)<2$ for all~$T\in\R$}}\\
&{\mbox{with~$a_3(T)=T-4$ for all~$T\le5$,}}\\& a_4(T):=1
\qquad{\mbox{and}}\qquad \eta_\star:=1.
\end{split}\end{equation}
As a matter of fact, when~$T=T_c$ the free energy according to the parameters in~\eqref{BSNTT-6}
is~$\eta^2(\eta-1)^2$, which coincides, up to translations and scaling, to the double-well potential in~\eqref{WDW}.

\begin{figure}[h]
\includegraphics[width=0.8\textwidth]{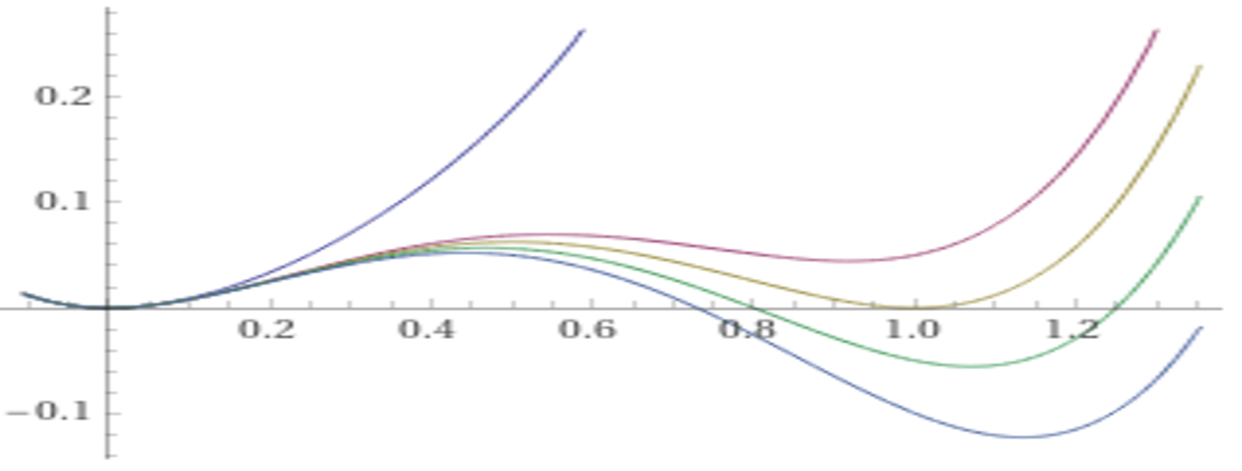}
\caption{Free energy corresponding to the parameters in~\eqref{BSNTT-6}
for~$T\in\{3,2.05,2,1.95,1.9\}$.}
        \label{9FIGURATBORDIbnNP}
\end{figure}

The free energy in~\eqref{GBSUJ-0} corresponding to the model case~\eqref{BSNTT-6}
for different values of the temperature~$T$ is sketched in Figure~\ref{9FIGURATBORDIbnNP}.
This situation has to be compared with that of second-order phase transitions in Figure~\ref{9FIGURATBORDINP}.

\begin{figure}[h]
\includegraphics[width=0.65\textwidth]{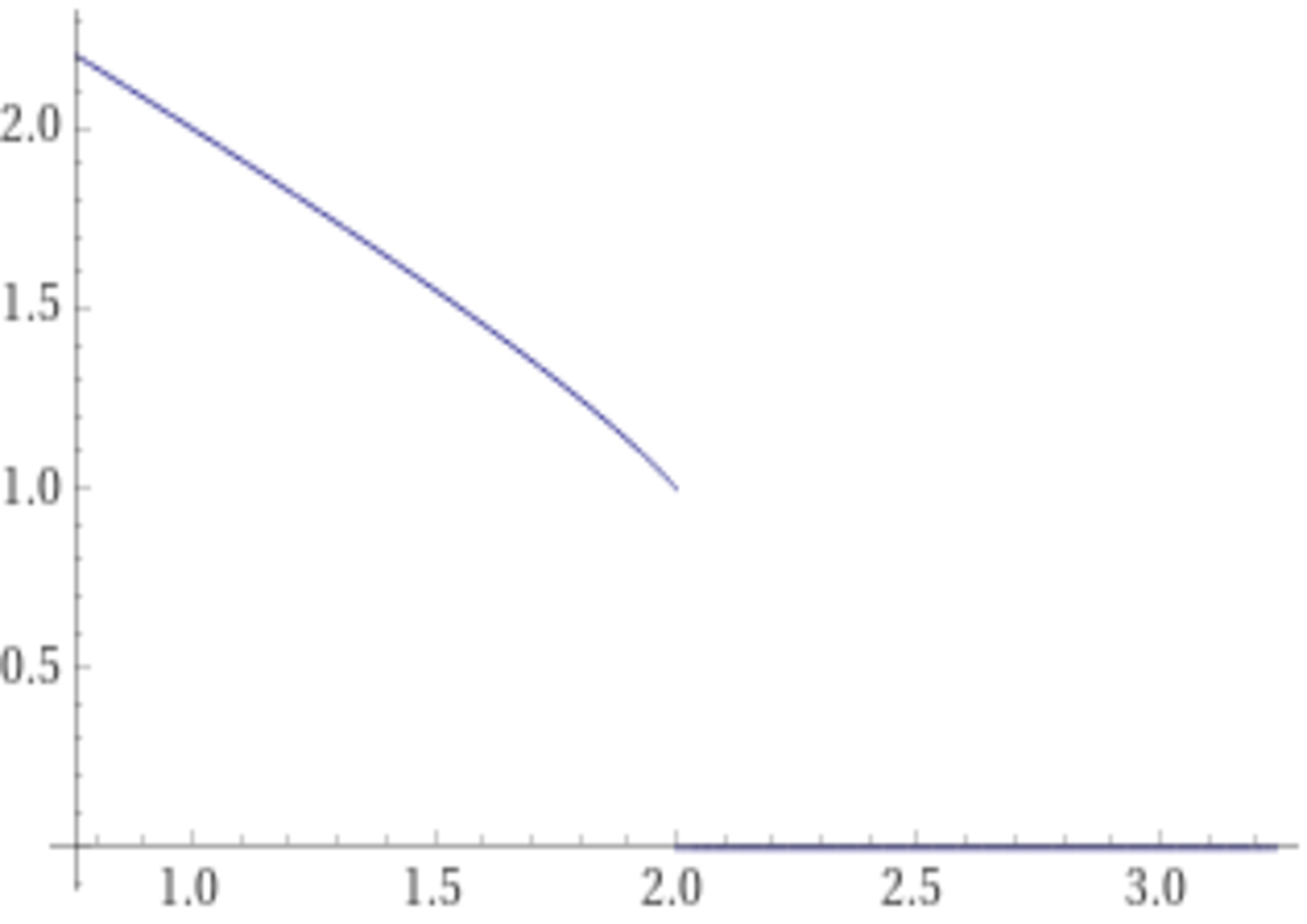}
\caption{Bifurcation diagram of the minimizers corresponding to the parameters in~\eqref{BSNTT-6}
for~$T\in[0,4]$.}
        \label{9FIGX2XcsURATBORDINP-1o}
\end{figure}

We stress that the nonzero critical points of the free energy in~\eqref{GBSUJ-0}
correspond to solutions~$\eta=\eta(T)$ of
\begin{equation*} 2a_2(T)+3a_3(T) \eta+4a_4(T) \eta^2=0\end{equation*}
and therefore when~$T<T_c$ the minimum 
takes the form
\begin{equation}\label{KAWPK-08ujn3rVFANOW} \eta(T):=\frac{-3 a_3(T)+\sqrt{9 a_3^2(T) - 32 \,a_2(T)\, a_4(T)}}{8 a_4(T)}
.\end{equation}
For instance, in the model case~\eqref{BSNTT-6} the global minima of the free energy are described by
\begin{equation}\label{T-MINI-2} \begin{dcases}
\{0\} & {\mbox{ if }} T> 2,\\
\{0,1\} & {\mbox{ if }} T= 2,
\\ \displaystyle\left\{
\frac{12 -3T+\sqrt{9 T^2 - 72 T + 112}}{8}
\right\} & {\mbox{ if }} T< 2,
\end{dcases} \end{equation}
see Figure~\ref{9FIGX2XcsURATBORDINP-1o}.

These minima should be compared with the situation in~\eqref{T-MINI}
(as well as Figure~\ref{9FIGX2XcsURATBORDINP-1o} should be compared
with the case of second-order phase transitions in Figure~\ref{9FIGURATBORDINP-1o}): in particular,
we stress that~\eqref{T-MINI-2} shows a discontinuous jump at the critical temperature
for the minima of the free energy, which 
corresponds to the abrupt formation of a new stable phase and which
constitutes a typical phenomenon for the so-called ``first-order
phase transitions''. 

As in Section~\ref{SEC:ORD:THE}, the name 
``first-order phase transition'' possibly arises from the regularity properties of
the free energy
in~\eqref{GBSUJ-0} evaluated at its global minima.
Namely, in this case the ``free energy
as a function of temperature'' takes the form
\begin{equation}\label{aeta2Tca3eta4Tc-XCVA} {\mathcal{E}}(T) :=
\begin{dcases}
0 & {\mbox{ if }}T\ge T_c,\\
a_2(T)\,\eta^2(T)+a_3(T) \,\eta^3(T)+a_4(T) \,\eta^4(T)
& {\mbox{ if }}T<T_c.
\end{dcases}
\end{equation}
As a result, by~\eqref{KAWPK-08ujn3rVFANOW},
\begin{equation}\label{aeta2Tca3eta4Tc}
\begin{split}&
\lim_{T\nearrow T_c} {\mathcal{E}}'(T) \\
\;=\,&
a_2'(T_c)\,\eta^2(T_c)+a_3'(T_c) \,\eta^3(T_c)+a_4'(T_c) \,\eta^4(T_c)\\&\qquad\qquad
+\big(2a_2(T_c)\,\eta(T_c)+3a_3(T_c) \,\eta^2(T_c)+4a_4(T_c) \,\eta^3(T_c)\big)\eta'(T_c)
\\
\;=\,&a_2'(T_c)\,\eta^2(T_c)+a_3'(T_c) \,\eta^3(T_c)+a_4'(T_c) \,\eta^4(T_c).
\end{split}\end{equation}
It is interesting to observe that this quantity is always nonnegative
(and strictly positive in ``nondegenerate'' cases): indeed, if
$$ \phi(T):= a_2(T)\eta^2(T_c)+a_3(T) \eta^3(T_c)+a_4(T) \eta^4(T_c),$$
it follows from~\eqref{BSNTT-3} and~\eqref{BSNTT-4} that~$\phi(T_c)=0\le\phi(T)$ for all~$T>T_c$,
whence
$$ 0\le\phi'(T_c)=a_2'(T_c)\eta^2(T_c)+a_3'(T_c) \eta^3(T_c)+a_4'(T_c) \eta^4(T_c).$$
{F}rom this and~\eqref{aeta2Tca3eta4Tc} we obtain that
\begin{equation}\label{aeta2Tca3eta4Tc2} \lim_{T\nearrow T_c} {\mathcal{E}}'(T)\ge0=\lim_{T\searrow T_c} {\mathcal{E}}'(T),\end{equation}
with strict inequality occurring whenever~$\phi'(T_c)\ne0$.

The strict inequality in~\eqref{aeta2Tca3eta4Tc2} corresponds to the typical situations
in the so-called first-order phase transitions, in which the derivative of the free energy with respect to temperature
is discontinuous at the critical temperature, which in turn corresponds to a physical situation in which a ``latent heat''
is emitted or absorbed by the system when the phase change occurs (recall footnote~\ref{order-phase0}).

As an example, one can consider the model situation in~\eqref{BSNTT-6}. In this case~\eqref{aeta2Tca3eta4Tc-XCVA} reduces to
\begin{eqnarray*} {\mathcal{E}}(T)& =&
\begin{dcases}
0 & \qquad{\mbox{ if }}T\ge 2,\\
\\
\begin{matrix}\displaystyle
\left( \frac{12 -3T+\sqrt{9 T^2 - 72 T + 112}}{8}\right)^2\\ \displaystyle\quad
+(T-4) \,\left(\frac{12 -3T+\sqrt{9 T^2 - 72 T + 112}}{8}\right)^3&\\ 
\displaystyle+\left(\frac{12 -3T+\sqrt{9 T^2 - 72 T + 112}}{8}\right)^4&
\end{matrix}
& \qquad{\mbox{ if }}T<2.
\end{dcases}\\
& =&
\begin{dcases}
0 & \qquad{\mbox{ if }}T\ge 2,\\
\\
\begin{matrix}\displaystyle
\frac{1}{512}\Big[
\sqrt{9 T^2 - 72 T + 112} \Big(9 T^3- 108 T^2 + 400 T -448\Big)\\
- 27 T^4 + 432 T^3- 2448 T^2  
+ 5760 T -4736 
\Big]
&
\end{matrix}
& \qquad{\mbox{ if }}T<2.
\end{dcases}
\end{eqnarray*}
where~\eqref{T-MINI-2} was used.

Accordingly, in this case,
$$ \lim_{T\nearrow 2} {\mathcal{E}}'(T)=1>0=\lim_{T\searrow 2} {\mathcal{E}}'(T),$$
showing the occurrence of the discontinuity
at the critical temperature of the first derivative of the free energy with respect to temperature.

For the sake of completeness, we also give an example in which,
in the framework of~\eqref{GBSUJ-0}, \eqref{BSNTT-1}, \eqref{BSNTT-2},
\eqref{BSNTT-3}, \eqref{BSNTT-4} and~\eqref{BSNTT-5}, a degenerate situation occurs, in which
the equality sign holds true in~\eqref{aeta2Tca3eta4Tc2}. This has some physical relevance because it entails
that the distinction between first- and second-order phase transitions has to be taken into account
with some caution: in particular, the next example shows that there are degenerate cases of
phase transitions whose free energy is modeled on~\eqref{GBSUJ-0}, \eqref{BSNTT-1}, \eqref{BSNTT-2},
\eqref{BSNTT-3}, \eqref{BSNTT-4} and~\eqref{BSNTT-5} but
whose first derivative of the free energy with respect to temperature happens to be continuous
(hence, not producing any ``latent heat'', these phase transitions should be effectively considered
second-order).

This degenerate example goes as follows.
We consider
\begin{equation}\label{BSNTT-6-DEGE} \begin{split}&
T_c:=2,\qquad
a_2(T):=1,\\
&{\mbox{$a_3$ a smooth and increasing function such that~$a_3(T)<2$ for all~$T\in\R$}}\\
&{\mbox{with~$a_3(T)=(T-2)^3-2$ for all~$T\le3$,}}\\& a_4(T):=1
\qquad{\mbox{and}}\qquad \eta_\star:=1.
\end{split}\end{equation}
This choice has to be compared with that in~\eqref{BSNTT-6}.
We stress that in this case~\eqref{BSNTT-1} and~\eqref{BSNTT-2} are obviously satisfied.
Moreover, if
$$g(T):=T^3 - 6 T^2 + 12 T - 12,$$ we have that~$g'(T)=3(T-2)^2\ge0$
and therefore,
when~$T\in(2,3]$,
$$ \frac{\big((T-2)^3-2\big)^2-4}{(T-2)^3}=
T^3 - 6 T^2 + 12 T - 12=g(T)<g(3)=-3<0
.$$
For this reason, we have that the discriminant of the polynomial
$$ x\mapsto P_T(x):= a_2(T) +a_3(T)\,x+a_4(T)\,x^2$$
is equal to~$a_3^2(T)-4a_2(T)a_4(T)=\big((T-2)^3-2\big)^2-4<0$
whenever~$T\in(2,3]$. Accordingly, since~$P_T(0)=a_2(T)=1>0$,
we have that~$P_T(x)>0$ for every~$T\in(2,3]$ and~$x\in\R$.

As a consequence,
if~$T\in(2,3]$ and~$\eta\in\R$,
\begin{eqnarray*}&&
a_2(T)+a_3(T) \eta+a_4(T) \eta^2=P_T(\eta)>0.
\end{eqnarray*}
Moreover, if~$T>3$ then~$a_3(T)\in(-1,2)$ and thus~$a_3^2-4<0$. As a result,
if~$T>3$ and~$\eta\in\R$,
$$ a_2(T)+a_3(T) \eta+a_4(T) \eta^2=1+a_3(T) \eta+\eta^2>0.$$
These considerations yield that~\eqref{BSNTT-3} is also satisfied.

Additionally, here
$$ a_2(T_c)\,\eta_\star^2+a_3(T_c) \eta_\star^3+a_4(T_c) \eta^4_\star=
1-2+1=
0,$$
which is~\eqref{BSNTT-4}, and, if~$T<2$,
\begin{eqnarray*}&&
\min_{\eta\in\R} a_2(T)\,\eta^2+a_3(T) \eta^3+a_4(T) \eta^4
\le a_2(T)+a_3(T)+a_4(T)\\&&\qquad=1+\big((T-2)^3-2\big)+1=
(T-2)^3 <0,\end{eqnarray*} 
which is~\eqref{BSNTT-5}, thus confirming that~\eqref{BSNTT-6-DEGE} fulfills the desired setting.

\begin{figure}[h]
\includegraphics[width=0.65\textwidth]{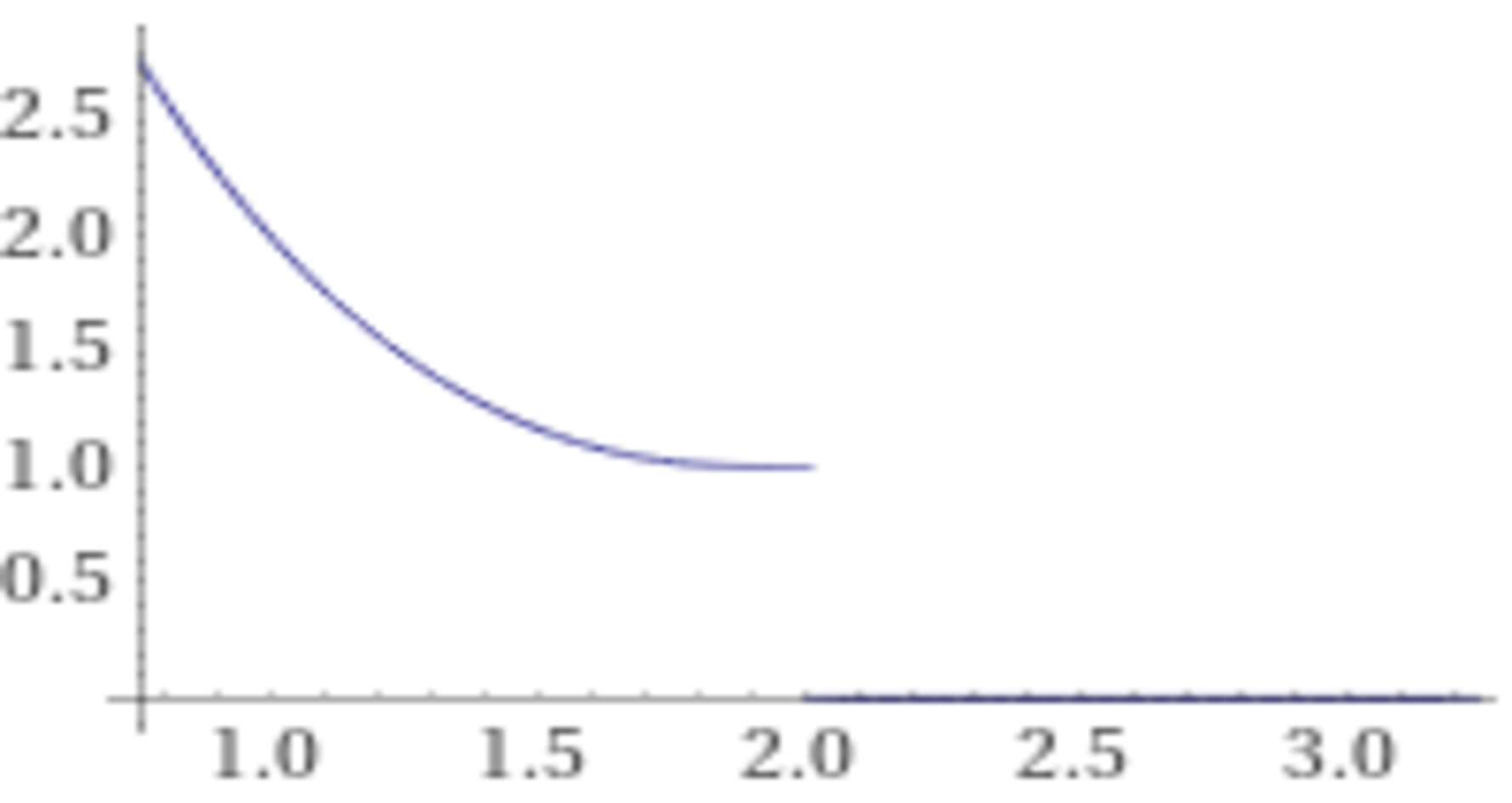}
\caption{Bifurcation diagram of the minimizers corresponding to the parameters in~\eqref{BSNTT-6-DEGE}
for~$T\in[0,4]$.}
        \label{9FIGX2XcsURATBORDINP-1o-DEGHE}
\end{figure}

In this case, by~\eqref{KAWPK-08ujn3rVFANOW}, we have that
the global minima of the free energy are described by
\begin{equation}\label{T-MINI-2-DEGE} \begin{dcases}
\{0\} & {\mbox{ if }} T> 2,\\
\{0,1\} & {\mbox{ if }} T= 2,
\\ \displaystyle\left\{
\begin{matrix}\displaystyle
\frac{1}{8}\Big(\sqrt{9 T^6- 108 T^5 + 540 T^4- 1476 T^3+ 2376 T^2- 2160 T+868}\\
- 3 T^3+ 18 T^2- 36 T+30\Big)\end{matrix}
\right\} & {\mbox{ if }} T< 2,
\end{dcases} \end{equation}
to be compared with~\eqref{T-MINI-2}.

The bifurcation diagram of the minima in~\eqref{T-MINI-2-DEGE} is sketched in Figure~\ref{9FIGX2XcsURATBORDINP-1o-DEGHE},
and it has to be compared with Figure~ \ref{9FIGX2XcsURATBORDINP-1o}.

\begin{figure}[h]
\includegraphics[width=0.8\textwidth]{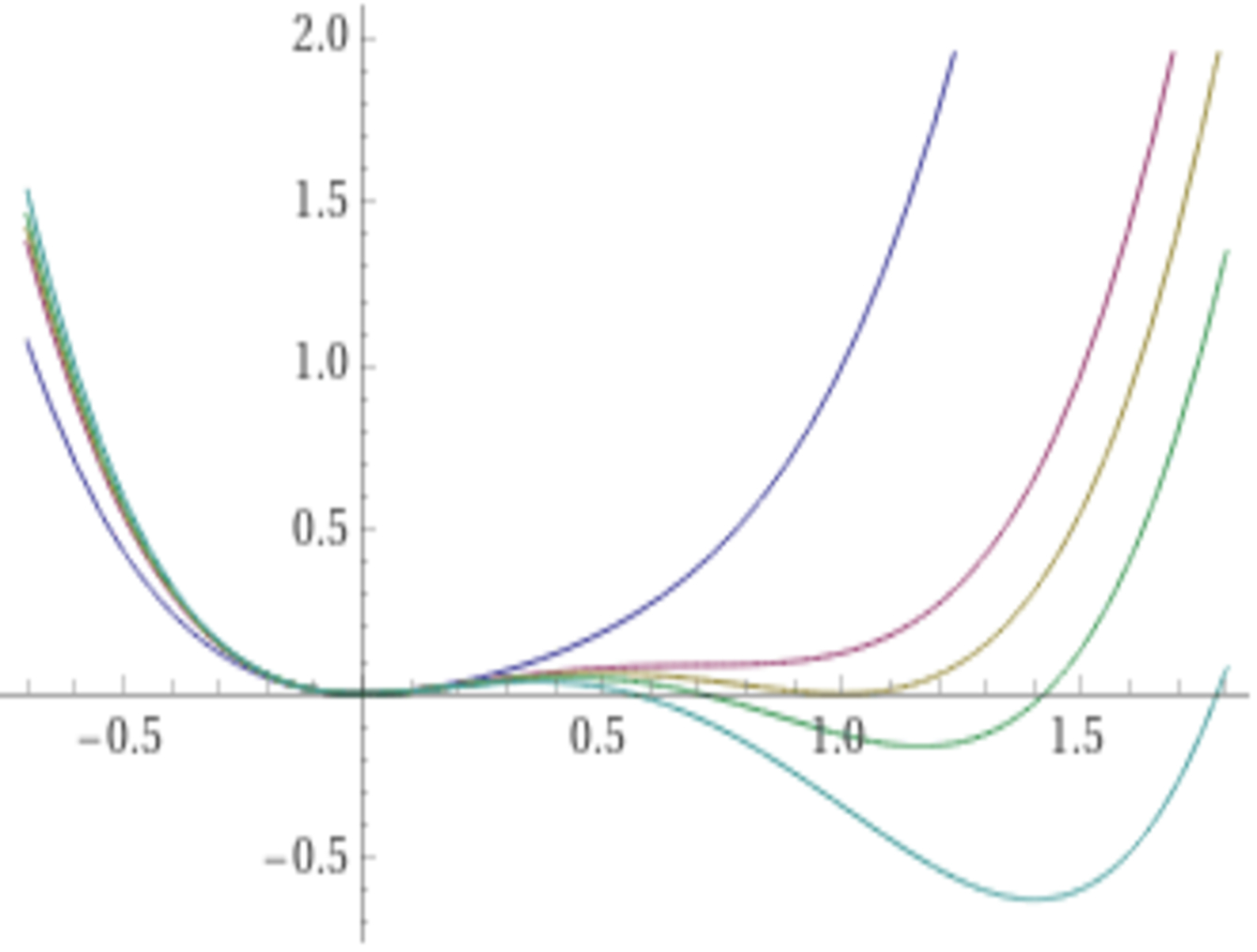}
\caption{Free energy corresponding to the parameters in~\eqref{BSNTT-6-DEGE}
for~$T\in\{3,2.5,2,1.5,1.3\}$.}
        \label{9FIGURATBORDIbnNP-090}
\end{figure}

The change of minimal levels for the free energy in this case is depicted in Figure~\ref{9FIGURATBORDIbnNP-090},
to be compared with Figure~\ref{9FIGURATBORDIbnNP}.

\begin{figure}[h]
\includegraphics[width=0.73\textwidth]{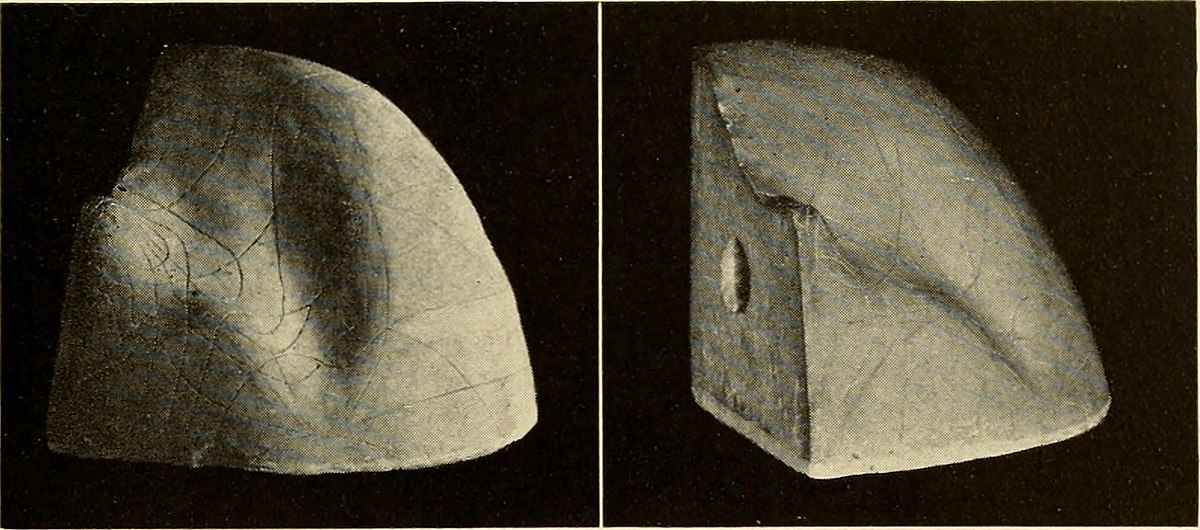}
\caption{Maxwell's thermodynamic surface (public domain image
from~\cite{COMME}).}
        \label{19FIGUR2COMME36t4y54eu65ATBORDIbnNP-090}
\end{figure}

\begin{figure}[h]
\includegraphics[height=0.31\textwidth]{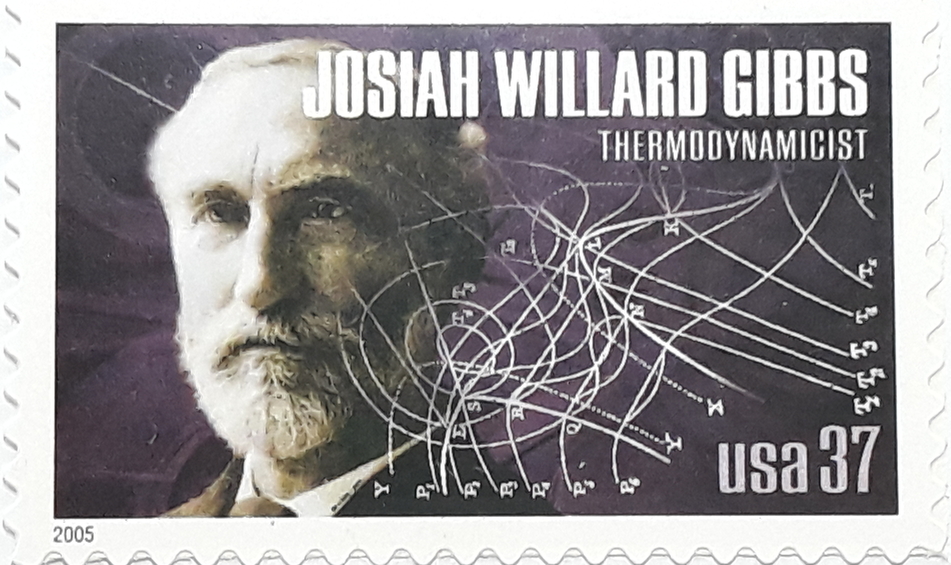} $\quad$
\includegraphics[height=0.31\textwidth]{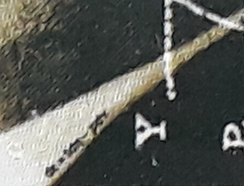}
\caption{Commemorative postage stamp
and zoom on the microprinting of the collar of Gibbs' portrait depicting his original mathematical equation for the change in internal energy (image
from~{\tt http://gigapan.com/gigapans/61023}).}
        \label{d19FIGUR267utgbjkCGIME36t4y54eu65ATBORDIbnNP-090}
\end{figure}

The relevant feature here is that, owing to~\eqref{aeta2Tca3eta4Tc},
\begin{eqnarray*}&&
\lim_{T\nearrow T_c} {\mathcal{E}}'(T)=
a_2'(T_c)\,\eta^2(T_c)+a_3'(T_c) \,\eta^3(T_c)+a_4'(T_c) \,\eta^4(T_c)\\&&\qquad=
0+3(2-2)^3+0
=0=\lim_{T\searrow T_c} {\mathcal{E}}'(T),
\end{eqnarray*}
hence the derivative of the free energy as a function of temperature is continuous at the critical temperature,
differently from what happened in case~\eqref{BSNTT-6}. Hence, some care is needed
to distinguish first- and second-order phase transitions and for this
it may not be sufficient\footnote{But figures are certainly essential to consolidate visual thinking.
For instance, inspired by Gibbs' work, in 1874 Maxwell spent an entire winter
to make a three-dimensional clay sculpture (also replicated in several plaster casts, one of which
was sent by Maxwell to Gibbs as a gift)
representing the energy of a fictitious substance with respect to volume and entropy,
see Figure~\ref{19FIGUR2COMME36t4y54eu65ATBORDIbnNP-090}.
In this sculpture, one can recognize the principal features of phase transitions
and latent heat formation via
simple geometrical operations,
such as placing a flat sheet of glass to mimic the tangent plane of the surface,
or placing the model in sunlight and tracing the curve when the rays graze the surface.

The specific setting of this sculpture is different from the one described
here in Figures~\ref{9FIGURATBORDINP}, \ref{9FIGURATBORDIbnNP}
and~\ref{9FIGURATBORDIbnNP-090}, since Maxwell was only describing
the internal energy of the system as a function of~$S$ and~$V$,
not the free energy corresponding to different order parameters
(roughly speaking, Figures~\ref{9FIGURATBORDINP}, \ref{9FIGURATBORDIbnNP}
and~\ref{9FIGURATBORDIbnNP-090} represent a free energy, at a given temperature, pressure and volume,
in dependence of the order parameter~$\eta$, and the minimization of this functional would produce
the phase observed in the system, which in turn would produce the corresponding observed free energy).

To appreciate the geometric arguments elucidated in Maxwell's thermodynamic surface,
we recall that the internal energy~$U$ has the form
$$ U=TS-PV,$$
therefore
$$ \frac{\partial U}{\partial S}=T\qquad{\mbox{and}}\qquad \frac{\partial U}{\partial V}=-P.$$
A normal vector to the surface~$U=U(S,V)$ is accordingly
$$ \nu:=\left( -\frac{\partial U}{\partial S},-\frac{\partial U}{\partial V}, 1
\right)=(-T,P,1).
$$
Consequently, if two points of the surface are touched by the same plane with normal~$\nu$, then
they necessarily
present the same temperature~$T$ and pressure~$P$
(but they possibly present different values of entropy
and internal energy). The physical system will then select the state with lower energy
and the difference of entropy (multiplied by temperature)
is related to latent heat (see footnote~\ref{order-phase0}).

Geometric considerations of this type allowed Maxwell to draw on the model the curves of
equal pressure and of equal temperature. 
In 2005, the United States Postal Service issued a $37$~cents commemorative postage stamp 
honoring Gibbs. Next to Gibbs's portrait,
the stamp features a diagram illustrating a thermodynamic surface.
Also, an almost invisible microprinting on the collar of Gibbs' portrait depicts 
the equation~$d\e=t\,d\eta-p\,dv$
(in Gibbs' notation, $\e$ stands for internal energy, $t$ for temperature, $\eta$ for entropy,
$p$ for pressure and~$v$ for volume, hence this formula is equivalent to the one that we used
to characterize the tangent vector to the thermodynamic surface).
See Figure~\ref{d19FIGUR267utgbjkCGIME36t4y54eu65ATBORDIbnNP-090}.

Here is another historical anecdote. In his scientific correspondence with other scientists,
Maxwell developed a jokey kind of code language. For instance, $\Sigma\phi\alpha\rho\xi$
stood
for spherical harmonics and~$\Theta\Delta$ics for thermodynamics.
See~\cite{zbMATH05047391} for a detailed
biography of James Clerk Maxwell.}
to only look at pictures
(for instance, it is difficult to distinguish from
Figures~\ref{9FIGURATBORDIbnNP}
and~\ref{9FIGURATBORDIbnNP-090} that
only the first corresponds to the creation of latent heat):
thus, as usual,
a detailed mathematical framework can come in handy.

\subsection{Interfacial energy of phase transitions}\label{ISS}
The description of phase models so far focused only on the favorable configurations of the free energy
which support one phase over the other, but we have not discussed how two different coexisting phases
are separated, that is what the geometry of the domains corresponding to each phase is.

For this, we start by observing that the coexistence of two phases occurs when they both attain the minimal value
of the free energy: this is precisely the case of first-order phase transitions at the critical temperature
and of second-order phase transitions at the critical temperature or below it,
see the discussions in Sections~\ref{SEC:ORD:THE} and~\ref{FIR:ORD:THE}, as well as Figures~\ref{9FIGURATBORDINP},
\ref{9FIGURATBORDIbnNP} and~\ref{9FIGURATBORDIbnNP-090}.

\begin{figure}[h]
\includegraphics[width=0.13\textwidth]{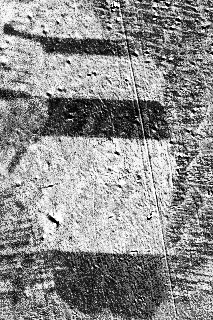}
\includegraphics[width=0.13\textwidth]{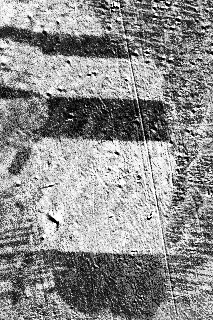}
\includegraphics[width=0.13\textwidth]{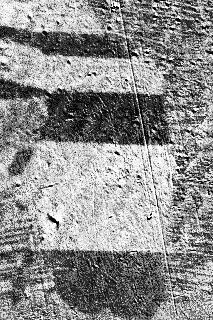}
\includegraphics[width=0.13\textwidth]{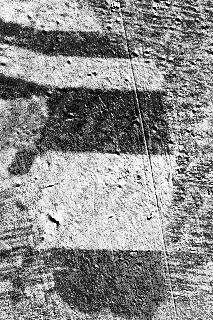}
\includegraphics[width=0.13\textwidth]{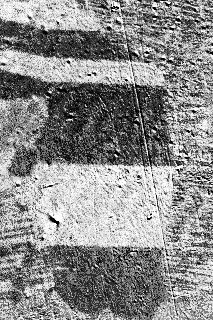}
\includegraphics[width=0.13\textwidth]{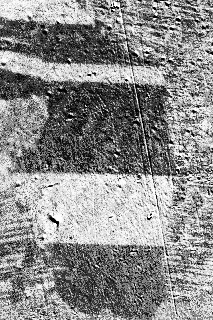}
\includegraphics[width=0.13\textwidth]{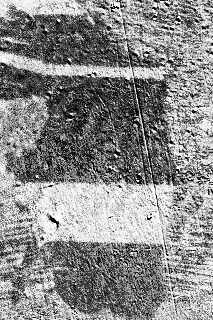} \\ \vskip0.1cm
\includegraphics[width=0.13\textwidth]{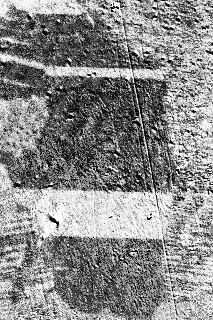}
\includegraphics[width=0.13\textwidth]{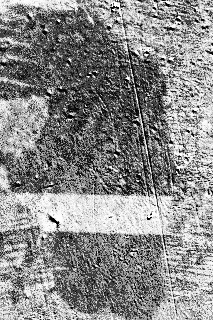}
\includegraphics[width=0.13\textwidth]{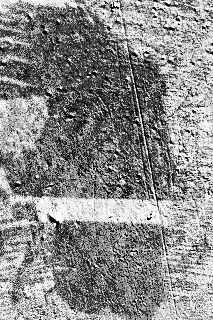}
\includegraphics[width=0.13\textwidth]{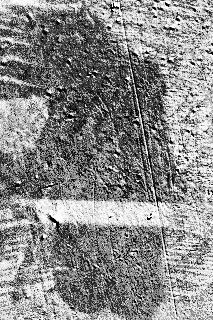}
\includegraphics[width=0.13\textwidth]{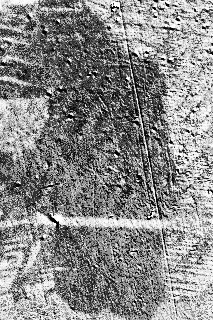}
\includegraphics[width=0.13\textwidth]{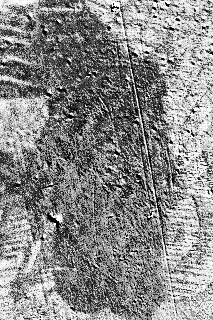}
\includegraphics[width=0.13\textwidth]{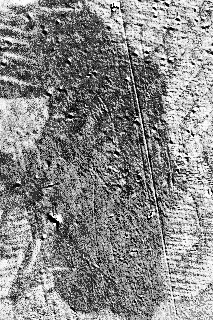}
\caption{Moving domain walls in a grain of silicon steel caused by an increasing external magnetic field. White areas are domains with magnetization directed up, dark areas are domains with magnetization directed down
(image from Wikipedia,
by Zureks, Chris Vardon,
CC BY-SA 3.0).}
        \label{19FIGUR236t4y54eu65ATBORDIbnNP-090}
\end{figure}

In principle, when two minima of the free energy occur at the same level,
the two phases are equally favorable from an energetic point of view, hence any configuration
in which any point of the state lies in any of the two phases is ``as good as any other''.
This however is in contradiction with common experience, since in many phenomena
the change of phase between different regions of the space occurs in very specific regions, see e.g.
Figure~\ref{19FIGUR236t4y54eu65ATBORDIbnNP-090}.

In practice, the model of ``pure phases'' is possibly ``too abrupt'', since small ``fluctuations'' produce values of the state
parameter which are not precisely equal to either of the two phases.
To model these fluctuations, one can consider the order parameter~$\eta$ as a function of the continuous spatial
coordinates and assume that the fluctuations are the byproduct of the mutual interaction between regions of spaces
corresponding to a different state parameter. Though a full understanding of these fluctuations is a very
delicate problem in statistical physics, a crude, but perhaps efficient, model is to assume that the
``effective'' energy of the system comes from the free energy, plus an interaction term between different phases.
The introduction of this interaction term dates back at least to Van der Waals, who considered
molecular interactions as averaged over long ranges. One of the additional benefits of such an interaction term
is, roughly speaking, to penalize the unnecessary changes of phase and somewhat favor, among all the configurations which minimize the free energy, the ones which present a ``minimal interface'' between regions with different phases.

The precise notion of minimal interface certainly depends on the additional interaction term that one takes
into account, thus we describe now some specific choices of interest. Considering the spatial domain the whole of~$\R^n$, one can take into account an interaction energy of the form
\begin{equation}\label{ITE56} \iint_{\R^n\times\R^n} (\eta(x)-\eta(y))^2\,{\mathcal{K}}(x,y)\,dx\,dy.\end{equation}
The quadratic power above is somewhat arbitrary\footnote{Ideally, one may want to determine a precise interaction kernel from general first principles.
However, due to the complexity of natural phenomena, in many concrete situations, the precise determination of an interaction kernel may be based on phenomenological considerations, interpretation of experimental data, or even, more frequently than one may think at a first, on the opportunity of finding useful computational simplifications.
As an example of ``convenient choice'' of an interaction kernel, we recall the fact that, in the
development of his new theory of gases based on statistical physics~\cite{zbMATH02723358},
James Clerk Maxwell introduced an interaction kernel
based on the inverse fifth power of the molecular distance.
Actually, it seems that the model was possibly taking into account
the general case of the $\kappa$th power of the distance: according to~\cite{BISTA},
``it has been argued that
Maxwell admitted that this choice of~$\kappa = 5$ [came] less from the
physical consequences of the choice than from the attractiveness of the possibility of explicit integration''.}
but it provides the advantage of corresponding to a linear equation (and, of course, being symmetric
if one exchanges the roles of~$x$ and~$y$ if~${\mathcal{K}}$ also presents this symmetry).
If we assume the kernel~${\mathcal{K}}$ to be symmetric under translations and rotations, we have that
\begin{equation}\label{ITE56a} {\mathcal{K}}(x,y)={\mathcal{K}}(x-y,y-y)={\mathcal{K}}(x-y,0)={\mathcal{K}}(|x-y|e_1,0)=:K(|x-y|).\end{equation}
If additionally the kernel is short-range, i.e. it vanishes whenever~$|x-y|\ge\varrho$, for some small~$\varrho>0$,
the quantity in~\eqref{ITE56} formally reduces to
\begin{equation}\label{JwOSef3LN4-5P579R0O}
\begin{split}
&
\int_{\R^n}\left[\int_{ B_\varrho(x)} (\eta(x)-\eta(y))^2\,K(|x-y|)\,dy\right]\,dx\\=\,&
\int_{\R^n}\left[\int_{ B_\varrho(x)} \Big(\nabla\eta(x)\cdot(x-y)+O(\sigma(x)\,|x-y|^2)\Big)^2\,K(|x-y|)\,dy\right]\,dx\\=\,&
\int_{\R^n}\left[\int_{ B_\varrho} \Big(\nabla\eta(x)\cdot z+O\big(\sigma(x)\,|z|^2\big)\Big)^2\,K(|z|)\,dz\right]\,dx\\
=\,&
\int_{\R^n}\left[\int_{ B_\varrho} \Big(\nabla\eta(x)\cdot z\Big)^2\,K(|z|)\,dz\right]\,dx
+
\int_{\R^n}\left[\int_{ B_\varrho} O(\sigma^2(x)\,|z|^3)\,K(|z|)\,dz\right]\,dx\\
=\,&
\int_{\R^n}\left[\int_{ B_\varrho} |\nabla\eta(x)|^2 \,z_1^2\,K(|z|)\,dz\right]\,dx
+
\int_{\R^n}\left[\int_{ B_\varrho} O(\sigma^2(x)\,|z|^3)\,K(|z|)\,dz\right]\,dx\\
=\,& C\int_{\R^n} |\nabla\eta(x)|^2 \,dx+O(C\varrho),
\end{split}\end{equation}
where
$$ \sigma(x):=\| \eta\|_{C^2(B_\varrho(x))}, \qquad
C:= \int_{ \R^n} z_1^2\,K(|z|)\,dz,$$
and a suitable decay on the kernel and on~$\sigma$ are assumed for integrability purposes.

Accordingly, for~$\varrho$ sufficiently small, the interaction term in~\eqref{ITE56}, up to normalizing constants,
can be approximated by $$ \int_{\R^n} |\nabla\eta(x)|^2 \,dx.$$
This and the discussions in  Sections~\ref{SEC:ORD:THE}
and~\ref{FIR:ORD:THE} give that, for short-range phase interactions,
we can approximately
describe the coexistence of two phases for first-order phase transitions at the critical temperature
and of second-order phase transitions at the critical temperature or below it via the energy functional
\begin{equation}\label{GUR} \frac12 \int_{\Omega} |\nabla\eta(x)|^2 \,dx+\int_\Omega W(\eta(x))\,dx,\end{equation}
where~$W$ is a double-well potential, as described in~\eqref{DOPP}
and~$\Omega\subseteq\R^n$ (we consider here the case in which the order parameter~$\eta$
is prescribed along~$\partial\Omega$; alternatively one can consider the case in which the average of~$\eta$, or of a function of~$\eta$,
in~$\Omega$ is prescribed).

The functional in~\eqref{GUR} is the prototypical example of classical phase coexistence energy studied
e.g. in~\cite{MR473971, MR866718, MR1097327}. The phase separation in this case
is dictated by the usual ``surface tension'' aiming at making the interface
a codimension~$1$ surface with the least possible $(n-1)$-dimensional area: to see this,
at least heuristically, one considers a rescaling of~\eqref{GUR} in which the gradient term
is explicitly a penalization of the double-well potential responsible of the phase separation.

That is, for a small parameter~$\e>0$, one takes into account the energy functional
\begin{equation}\label{GUR2} \frac{\e}2 \int_{\Omega} |\nabla\eta(x)|^2 \,dx+\frac1\e\int_\Omega W(\eta(x))\,dx.\end{equation}
By the Cauchy-Schwarz Inequality and the Coarea Formula, one can bound this quantity from below by
\begin{eqnarray*}
\int_{\Omega} |\nabla\eta(x)|\sqrt{2 W(\eta(x))}\,dx&=&\int_{-1}^1\left[
\int_{\Omega\cap\{\eta(x)=\tau\}}\sqrt{2 W(\tau)}\,d{\mathcal{H}}^{n-1}_x
\right]\,d\tau\\&=&
\int_{-1}^1\sqrt{ 2 W(\tau)}\,
{\mathcal{H}}^{n-1}\big( \Omega\cap\{\eta=\tau\}\big)\,d\tau,
\end{eqnarray*}
where~${\mathcal{H}}^{n-1}$ denotes the~$(n - 1)$-dimensional Hausdorff measure. We recall that $-1$ and~$1$ represent the pure phases of the system.

Also, for small~$\e$, we may think that the energy minimizers try to ``optimize'' the above lower bound
and to sit in the zeros (or close to the zeros) of the double-well potential as much as possible.
Therefore, for small~$\e$, the minimal energy in~\eqref{GUR2} is expected to be related
to
$$ c{\mathcal{H}}^{n-1}\big( \Omega\cap (\partial E)\big)
\qquad{\mbox{where}}\qquad
c:=\int_{-1}^1\sqrt{2 W(\tau)}\,d\tau,$$
being~$E$ a set in which the order parameter is ``essentially'' equal to~$+1$
and the complement of~$E$ a set in which the order parameter is ``essentially'' equal to~$-1$.
A precise formulation for this phenomenon will be given in Section~\ref{GAMMACO}.
See also~\cite{MR768066} for an alternative approach to interfacial energy for
phase transitions.

For long-range interactions, the gradient approximation in~\eqref{JwOSef3LN4-5P579R0O}
is not available anymore and, in light of~\eqref{ITE56}
and~\eqref{ITE56a}, it is more opportune to replace~\eqref{GUR} with a nonlocal energy term of the form 
\begin{equation}\label{GURNL} \frac14 \iint_{{\mathcal{Q}}(\Omega)} (\eta(x)-\eta(y))^2\,K(|x-y|)\,dx\,dy+\int_\Omega W(\eta(x))\,dx,\end{equation}
where
\begin{equation}\label{GURNL-2} {\mathcal{Q}}(\Omega):=(\Omega\times\Omega)\cup(\Omega^c\times\Omega)\cup
(\Omega\times\Omega^c),\end{equation}
being~$\Omega^c:=\R^n\setminus\Omega$.

The notation~${\mathcal{Q}}(\Omega)$ stands for a
``cross-shaped set'' (``${\mathcal{Q}}$'' stands for cross, since~``$C$'' is used for constants
and~``$K$'' for kernels!). The intuition behind this set is that
we are prescribing here the order parameter~$\eta$ outside the domain~$\Omega$,
which is the ``global'' counterpart of prescribing~$\eta$ along~$\partial\Omega$
in~\eqref{GUR}. Accordingly, the energy functional in~\eqref{GURNL} should account for
all the configurations which involve the values of the state parameter in~$\Omega$:
whatever piece of energy containing only the values
of the state parameter outside~$\Omega$ is prescribed, whence does not influence
energy minimization (interestingly, in this way, one considers the energy confined outside the domain
as ``constant'', even if this constant can actually be infinite!).
In this spirit, the cross-shaped set in~\eqref{GURNL-2}
accounts for all the phase interactions in which at least one of the sites is located in~$\Omega$.

Special cases of kernels~$K$ are the ones which are positively homogeneous of some degree~$d$,
that is~$K(|tz|)=t^d K(|z|)$ for all~$z\in\R^n\setminus\{0\}$ and~$t\in(0,+\infty)$.

Note that the degree~$d$ cannot be arbitrarily chosen in the reals, since to make sense of the interaction energy
in~\eqref{GURNL} it is desirable to have it finite at least when~$\eta\in C^\infty_0(\Omega,\,[0,2])$
with~$\eta(x)=2-|x-x_0|^2$ for all~$x\in B_r(x_0)$,
for some small~$r\in(0,1)$ such that~$B_{2r}(x_0)\subseteq\Omega\subseteq B_{1/r}(x_0)$.
Hence,
if~$x\in B_r(x_0)$, then~$\eta(x)\ge 2-r^2>1$.
This yields that, for every~$\theta\in(0,+\infty)$,
\begin{eqnarray*} +\infty&>&\iint_{B_r(x_0)\times(\R^n\setminus B_{1/r}(x_0))} (\eta(x)-\eta(y))^2\,K(|x-y|)\,dx\,dy\\&
\ge&\iint_{B_r(x_0)\times(\R^n\setminus B_{1/r}(x_0))} K(|x-y|)\,dx\,dy\\
&\ge&\iint_{B_r(x_0)\times(\R^n\setminus B_{2/r})} K(|z|)\,dx\,dz\\&=&
|B_r|\int_{\R^n\setminus B_{2/r}} K\left(\left|\frac{|z|}\theta\,\theta e_1\right|\right)\,dz\\&=&
|B_r|\int_{\R^n\setminus B_{2/r}} 
\frac{|z|^d}{\theta^d}
K\left(\left|\theta e_1\right|\right)\,dz\\&=&
\frac{|B_r|\,K(\theta)}{\theta^d}\int_{\R^n\setminus B_{2/r}} 
|z|^d\,dz
\end{eqnarray*}
and accordingly
\begin{equation}\label{APPL9weGUIS-1}
d<-n,
\end{equation}
unless of course~$K$ vanishes identically.

In a similar fashion, using the substitutions~$z:=x-y$ and~$w:=2(x-x_0)-z$,
\begin{eqnarray*} +\infty&>&\iint_{B_r(x_0)\times B_r(x_0)} (\eta(x)-\eta(y))^2\,K(|x-y|)\,dx\,dy\\&
=&\iint_{B_r(x_0)\times B_r(x_0)} \big(|x-x_0|^2-|y-x_0|^2\big)^2\,K(|x-y|)\,dx\,dy
\\&
=&\iint_{B_r(x_0)\times B_r(x_0)} \big((x+y-2x_0)\cdot(x-y)\big)^2\,K(|x-y|)\,dx\,dy\\&\ge&
\iint_{B_{r/2}(x_0)\times B_{r/2}} \big((2(x-x_0)-z)\cdot z\big)^2\,K(|z|)\,dx\,dz
\\&\ge&\frac1{2^n}
\int_{B_{r/8}}\left[\int_{B_{r/8}\cap\left\{w\cdot \frac{z}{|z|}\ge \frac12\right\}} \big(w\cdot z\big)^2\,K(|z|)\,dw\right]\,dz
\\&\ge&
\frac1{2^{n+2}}\int_{B_{r/8}}\left[\int_{B_{r/8}\cap\left\{w\cdot \frac{z}{|z|}\ge \frac12\right\}} 
|z|^2\,K\left(\left|\frac{|z|}\theta\,\theta e_1\right|\right)\,dw\right]\,dz\\&=&
\frac{\left|B_{r/8}\cap\left\{w_1\ge \frac12\right\}\right|}{2^{n+2}}\int_{B_{r/8}}
\frac{|z|^{2+d}}{\theta^d}\,
K\left(\left|\theta e_1\right|\right)\,dz\\&=&
\frac{\left|B_{r/8}\cap\left\{w_1\ge \frac12\right\}\right|\,K(\theta)}{2^{n+2}\theta^d}\int_{B_{r/8}}
|z|^{2+d}\,dz,
\end{eqnarray*}
whence
\begin{equation}\label{APPL9weGUIS-2}
d>-n-2,
\end{equation}
unless~$K$ vanishes identically.

{F}rom~\eqref{APPL9weGUIS-1} and~\eqref{APPL9weGUIS-2} (and using the normalization~$K(1):=1$), it follows that
necessarily
$$ K(|z|)=|z|^d\,K(|e_1|)=\frac1{|z|^{n+\alpha}},$$
for some~$\alpha\in(0,2)$. Thus in this case~\eqref{GURNL} boils down to
\begin{equation}\label{GURNL-sal}
\iint_{{\mathcal{Q}}(\Omega)} \frac{(\eta(x)-\eta(y))^2}{|x-y|^{n+\alpha}}\,dx\,dy+\int_\Omega W(\eta(x))\,dx,\end{equation}
which is the interaction energy presented in~\eqref{PRES}.

When~$\alpha\in(0,1)$, minimizers of~\eqref{GURNL-sal} can be easily related to a geometric
minimization problem, since if
$$\eta(x)=\chi_E(x)-\chi_{E^c}(x)=\begin{dcases}1 & {\mbox{ if }}x\in E,\\
-1 & {\mbox{ if }}x\in E^c,
\end{dcases}$$ for some~$E\subseteq\R^n$, then the energy functional in~\eqref{GURNL-sal}
boils down to
\begin{equation*}\begin{split}
&\iint_{{\mathcal{Q}}(\Omega)\cap(E\times E^c)} \frac{4}{|x-y|^{n+\alpha}}\,dx\,dy=
\iint_{(E\cap\Omega)\times (E^c\cap\Omega)} \frac{4}{|x-y|^{n+\alpha}}\,dx\,dy\\
&\qquad\qquad
+
\iint_{(E\cap\Omega^c)\times (E^c\cap\Omega)} \frac{4}{|x-y|^{n+\alpha}}\,dx\,dy
+
\iint_{(E\cap\Omega)\times (E^c\cap\Omega^c)} \frac{4}{|x-y|^{n+\alpha}}\,dx\,dy
\\&\qquad\qquad\qquad={\rm Per}_\alpha(E,\Omega),
\end{split}\end{equation*}
due to~\eqref{per-a}.

This observation will be better formalized in Section~\ref{GAMMACO2},
in which we will revisit the main results of~\cite{MR2948285}:
in this situation, the limit interfaces of~\eqref{GURNL-sal} will be rigorously related
to the minimizers of the $\alpha$-perimeter when~$\alpha\in(0,1)$.
Interestingly, when~$\alpha\in[1,2)$,
the limit interfaces of~\eqref{GURNL-sal} will be instead related
to the minimizers of the classical perimeter, showing a remarkable ``localization effect
for nonlocal energies'' when the interaction parameter~$\alpha$ is larger than or equal to~$1$.

For additional information on the phase coexistence modelization, see e.g.~\cite{zbMATH02132471, MR2162511, zbMATH05046576, dipierro2021elliptic} and the references therein.

\section{The Allen-Cahn equation}\label{Skdfe0rthX234rt}

Critical points of the energy functional~\eqref{WDW-CLASSi} give rise to the equation
\begin{equation}\label{ACHAW}
\Delta u(x)=W'(u(x))\end{equation}
for~$x\in\Omega$. The model case in which the potential takes the form in~\eqref{WDW}
reduces to
\begin{equation}\label{ALLEC}-\Delta u=u-u^3,\end{equation} which is known as the Allen-Cahn equation. This equation indeed produces the stationary states
of an evolution equation presented in~\cite{ALLEN1972423} with the specific goal of
describing the phase separation in multi-component alloy systems. Moreover, solutions of the Allen-Cahn equation
are also stationary states of the so-called Cahn-Hilliard equation
\begin{equation}\label{CAHIL}\partial_t u=\Delta\left(u^{3}-u-\Delta u\right),\end{equation}
which was introduced in~\cite{doi:10.1063/1.1744102}
to represent
the process of  spontaneous phase separation in a binary fluid.

\begin{figure}[h]
\includegraphics[width=0.19\textwidth]{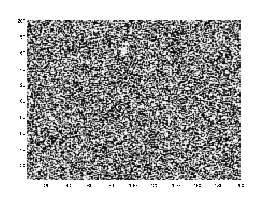}
\includegraphics[width=0.19\textwidth]{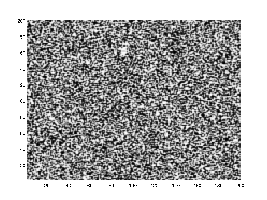}
\includegraphics[width=0.19\textwidth]{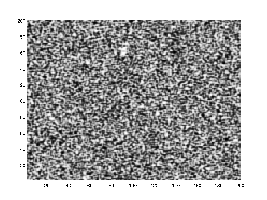}
\includegraphics[width=0.19\textwidth]{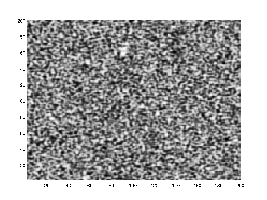}
\includegraphics[width=0.19\textwidth]{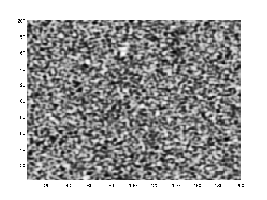}\\ \vskip0.1cm
\includegraphics[width=0.19\textwidth]{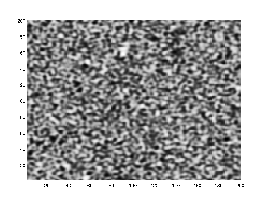}
\includegraphics[width=0.19\textwidth]{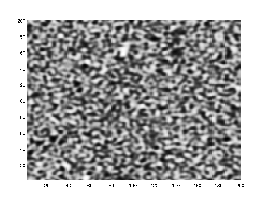}
\includegraphics[width=0.19\textwidth]{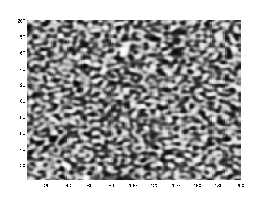} 
\includegraphics[width=0.19\textwidth]{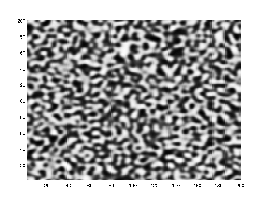}
\includegraphics[width=0.19\textwidth]{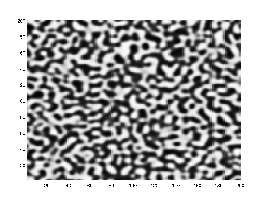}\\ \vskip0.1cm
\includegraphics[width=0.19\textwidth]{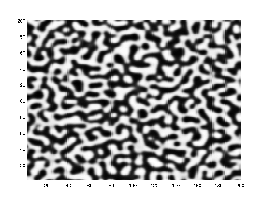}
\includegraphics[width=0.19\textwidth]{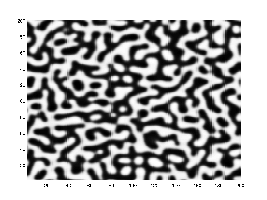}
\includegraphics[width=0.19\textwidth]{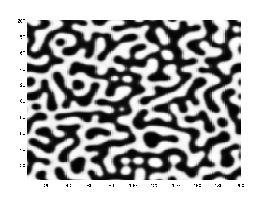}
\includegraphics[width=0.19\textwidth]{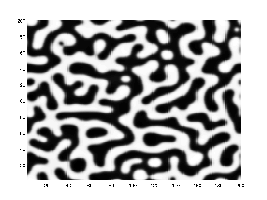}
\includegraphics[width=0.19\textwidth]{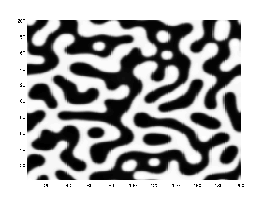}
\caption{Evolution of random initial data under the Cahn-Hilliard equation
(image from Wikipedia, by Yonatan Oren, public domain).}
        \label{19FIGUR236t4y54eu65ATBORDIbnNP-090PD}
\end{figure}

The success of equation~\eqref{CAHIL} in
describing spontaneous phase separation with
a tendency of similar phases to cluster together
is indeed quite perceptible, see Figure~\ref{19FIGUR236t4y54eu65ATBORDIbnNP-090PD}.

The above equations can also be modified to account for a heterogeneous material.
For instance, the potential in~\eqref{ACHAW} can be modulated by a spatially dependent function~$Q=Q(x)$,
ranging in~$[\underline{Q},\overline{Q}]$, for some~$\overline{Q}\ge \underline{Q}>0$. In this setting,
equation~\eqref{ACHAW} can be generalized to
\begin{equation}\label{ACHAW2}
\Delta u(x)=Q(x)\,W'(u(x))\end{equation}
and further generalization are also possible (e.g. by considering more general spatially dependent
double-well potentials or spatially dependent elliptic operators, say in divergence form
to have a consistent link with the arguments related to energy functionals in Section~\ref{LANDATHE}).

An interesting feature of the Allen-Cahn equation, as well as of the equations closely related to it,
is to be strongly cross-disciplinary. And this not only because its deep relation with phase coexistence
models makes it a standard tool of investigation for physicists, chemists, biologists and engineers too,
but also because it presents an immediate link with the theory of dynamical systems,
in which the space dimension~$n$ is taken to be~$1$, and the role space is replaced by time.
To clarify this point, one can consider equation~\eqref{ACHAW} when~$n=1$
and replace the name of the variable~$x$ by~$t$ (to physically represent time), thus obtaining
\begin{equation}\label{PEANH} \ddot u(t)=W'(u(t)).\end{equation}
For instance, when, for all~$u\in[-1,1]$, one has~$W(u):=1-\cos(\pi u)$, equation~\eqref{PEANH}
reduces to~$\ddot u(t)=\pi\sin(\pi u(t))$, which is the equation of the simple pendulum.
A similar procedure applied to~\eqref{ACHAW2} produces the equation~$\ddot u(t)=\pi Q(t)\sin(\pi u(t))$, which is the equation of a pendulum in a variable in time gravitational field.

This type of relations inspired a series of works in which methodologies typical of dynamical systems
(such as the so-called Aubry-Mather theory~\cite{MR670747, MR719634},
also related to classical works on geodesics~\cite{MR1501263}) are employed to
study phenomena strictly related to elliptic partial differential equations
and the Allen-Cahn equation, see e.g.~\cite{MR847308, MR991874, MR1782992, MR2354993, MR2416096, MR2505410, MR2809349} and the references therein.

For example, the theory of dynamical systems often focuses on the determination of orbits
with given ``rotation numbers'': for instance, especially in periodic settings, it is common to seek trajectories
of ordinary differential
equations exhibiting
a linear growth, up to bounded oscillations, that is such that
\begin{equation}\label{GEOTRAPL-S}
|u(t)-\rho t|\le K\qquad{\mbox{for all }}t\in\R,
\end{equation}
for some~$\rho\in\R$ and~$K\ge0$.

{F}rom a geometric viewpoint, one can restate~\eqref{GEOTRAPL-S} by saying that the graph of~$u$,
that is~${\mathcal{G}}_u:=\big\{(t,u(t)),$ $t\in\R\big\}$, is trapped in a slab perpendicular to~$\omega:=\frac{(-\rho,1)}{\sqrt{1+\rho^2}}$
and of width~$2M$, with~$M:=\frac{K}{\sqrt{1+\rho^2}}$,
namely
\begin{equation}\label{MS-PSLDibvAlcK1}
{\mathcal{G}}_u\subseteq\big\{
x\in\R^2{\mbox{ s.t. }}|x\cdot\omega|\le M
\big\}.
\end{equation}

With respect to this matter, a natural counterpart for phase transitions would be to construct
examples in which the interface
$$ {\mathcal{I}}_{u,\vartheta}:=\big\{
x\in\R^n {\mbox{ s.t. }}|u(x)| <\vartheta\big\}$$
is trapped in a slab (or equivalently, trapped within two hyperplanes).
This is indeed the content of the following result (which is a particular case of Theorem~8.1
in~\cite{MR2099113}):

\begin{theorem}\label{MS-PSLDibvAlcK2}
Let~$n\ge2$, $\vartheta\in(0,1)$ and assume that~$Q(x+m)=Q(x)$ for all~$m\in\Z^n$.

Then, there exists a constant~$M> 0$, depending only on~$\vartheta$, $n$, $Q$ and~$W$, such that, given any direction~$\omega\in\partial B_1$, we can construct a solution~$u$ of~\eqref{ACHAW2}
such that
\begin{equation}\label{MS-PSLDibvAlcK} {\mathcal{I}}_{u,\vartheta}\subseteq
\big\{
x\in\R^2{\mbox{ s.t. }}|x\cdot\omega|\le M
\big\}.
\end{equation}
\end{theorem}

The reader can appreciate the similarity between~\eqref{MS-PSLDibvAlcK1}
and~\eqref{MS-PSLDibvAlcK}. For completeness, we also mention that
the solution constructed in Theorem~\ref{MS-PSLDibvAlcK2} enjoys several extra features,
such as it is a local minimizer of the associated energy functional,
it is periodic when~$\omega$ is rationally dependent, its interface fulfills a suitable non-self-intersecting property, etc.

\section{The limit interface and the theory of~$\Gamma$-convergence}\label{GAMMACO}

\subsection{Singular perturbations and $\Gamma$-convergence theory}
We now come back to the perturbation problem induced by the energy functional in~\eqref{GUR2}.
Compared with~\eqref{ACHAW}, this setting produces a ``rescaled version'' of the Allen-Cahn equation of the form
\begin{equation}\label{EPSACAN}
\e^2\Delta u(x)=W'(u(x)).\end{equation}
This can be seen as a ``singular perturbation'' problem, since
the main term in the equation would disappear when~$\e=0$
and therefore a delicate analysis is required to capture the essential features of the problem
for small values of the parameter~$\e$.

As a matter of fact,
this kind of singular perturbation analysis is one of the main fillip to the development
of the notion of~$\Gamma$-convergence, see~\cite{MR375037}.
To appreciate this theory, let us consider a functional of the form~${\mathcal{F}}_\e$
and let us try to discuss a convenient meaning for a suitable convergence of~${\mathcal{F}}_\e$
to some~${\mathcal{F}}$ as~$\e\searrow0$. Notice that a pointwise convergence could be out of reach,
because a singular perturbation problem may drastically change the structure of the limit functional
as well as its natural domain of definition, therefore a ``new'' notion of convergence is called for.
In particular, to make this notion practical and serviceable, it is desirable to keep the notion of local energy minimizers
in the limit: namely,
\begin{equation}\label{OKSPREFBINV}
\begin{split}&
{\mbox{if~$u_\e$ is a local minimizer for~${\mathcal{F}}_\e$}}\\
&{\mbox{and~$u_\e\to u$ in some topology~$X$ as~$\e\searrow0$,}}\\
&{\mbox{this functional notion of convergence should entail}}
\\ &{\mbox{that~$u$ is a local minimizer for~${\mathcal{F}}$.}}\end{split}\end{equation} To this extent, the limit functional~${\mathcal{F}}$
may be considered as an ``effective energy'' and the choice\footnote{For simplicity, we will always implicitly assume that
the topology of~$X$ is induced by a metric space, so that compactness and sequential compactness are the same.}
of the topology~$X$
can be possibly made ``loose enough'' to ensure compactness of the minimizers beforehand
(choosing a ``too strong'' topology~$X$ produces the pitfall that minimizers
may not converge!).

It is also desirable that
\begin{equation}\label{MSJNDNTIONDS}
{\mbox{the limit functional~${\mathcal{F}}$ is lower semicontinuous}}\end{equation}
in order to develop a solid existence theory for its minimizers.

With these remarks in mind, it is not too difficult to ``guess'' what an ``appropriate'' notion
of functional convergence should be. 
To this end, we distinguish between
the lower limit and the upper limit.
For the lower limit we take inspiration from the classical Fatou's Lemma
(after all, it is sensible that a good functional convergence turns out to be compatible with the classical scenarios)
in which one considers the very special case of~$u_\e$ being a sequence of nonnegative measurable functions converging pointwise to~$u$,
takes~${\mathcal{F}}_\e(v):={\mathcal{F}}(v):=\int_{\R^n}v(x)\,dx$ and writes that
$$ \liminf_{\e\searrow0} {\mathcal{F}}_\e(u_\e)=\liminf_{\e\searrow0}
\int_{\R^n}u_\e(x)\,dx\ge\int_{\R^n}u(x)\,dx
={\mathcal{F}}(u).$$
Hence, a natural requirement for a general notion of functional convergence is that
\begin{equation}\label{REQ:LIMINF}
{\mbox{whenever~$u_\e\to u$ in~$X$,}}\quad
\liminf_{\e\searrow0} {\mathcal{F}}_\e(u_\e)\ge{\mathcal{F}}(u).\end{equation}

Let us now consider an upper limit condition in the light of~\eqref{OKSPREFBINV}.
Namely, let us consider a minimizer~$u_\e^\star$ for~${\mathcal{F}}_\e$
(say, under suitable boundary or external conditions, or mass prescriptions, or so): then, it holds that,
for every competitor~$u_\e$ for~$u_\e^\star$,
\begin{equation}\label{0j24rhcalFeu}
{\mathcal{F}}_\e(u_\e)\ge{\mathcal{F}}_\e(u_\e^\star)
\end{equation}
and the recipe in~\eqref{OKSPREFBINV} would suggest to obtain the same structural inequality as~$\e\searrow0$.
To wit, for every competitor~$u$ for~$u^\star$, we aim at showing that
\begin{equation}\label{0j24rhcalFeu-GO}
{\mathcal{F}}(u)\ge{\mathcal{F}}(u^\star).
\end{equation}
For this, if~$u_\e$ is any sequence converging to~$u$ in~$X$, we know from~\eqref{REQ:LIMINF}
and~\eqref{0j24rhcalFeu} that
$$ \limsup_{\e\searrow0} {\mathcal{F}}_\e(u_\e)\ge\liminf_{\e\searrow0} {\mathcal{F}}_\e(u_\e)\ge\liminf_{\e\searrow0} {\mathcal{F}}_\e(u_\e^\star)\ge{\mathcal{F}}(u^\star).
$$
Therefore, to obtain~\eqref{0j24rhcalFeu-GO}, it suffices to find one, possibly very special sequence~$u_\e$
converging to~$u$ in~$X$ for which
\begin{equation}\label{REQ:LIMSUP-224} {\mathcal{F}}(u)\ge \limsup_{\e\searrow0} {\mathcal{F}}_\e(u_\e).\end{equation}
This special sequence making the job is sometimes called ``recovery sequence'', the name coming from the
observation that if~$u$ is as in~\eqref{REQ:LIMSUP-224}, by~\eqref{REQ:LIMINF} one in fact has that
$$ {\mathcal{F}}(u)= \lim_{\e\searrow0} {\mathcal{F}}_\e(u_\e).$$
Thus, a natural upper limit condition to complement~\eqref{REQ:LIMINF} consists in requiring that
\begin{equation}\label{REQ:LIMSUP}
\begin{split}&{\mbox{there exists a sequence~$u_\e\to u$ as~$\e\searrow0$ in~$X$ such that}}\\ &
\limsup_{\e\searrow0} {\mathcal{F}}_\e(u_\e)\le{\mathcal{F}}(u).\end{split}\end{equation}

Conditions~\eqref{REQ:LIMINF} and~\eqref{REQ:LIMSUP} are often accompanied by
a compactness assumption under a bounded energy requirement, such that
\begin{equation}\label{REQ:LIMCOMPA}
\begin{split}&{\mbox{if }}\sup_{\e\in(0,1)} {\mathcal{F}}_\e(u_\e)<+\infty,
\\ & {\mbox{then there exists a subsequence~$u_{\e'}$ converging in~$X$ as~$\e'\searrow0$.}}\end{split}\end{equation}
When conditions~\eqref{REQ:LIMINF}, \eqref{REQ:LIMSUP} and~\eqref{REQ:LIMCOMPA} are met, then we say that~${\mathcal{F}}_\e$ $\Gamma$-converges to~${\mathcal{F}}$.

We also observe that conditions~\eqref{REQ:LIMINF} and~\eqref{REQ:LIMSUP} entail the lower semicontinuity property mentioned in~\eqref{MSJNDNTIONDS},
because if~$w_k\to w$ in~$X$ as~$k\to+\infty$ we can take a subsequence~$w_{k_j}$ such that
$$ \liminf_{k\to+\infty} {\mathcal{F}}(w_{k})=\lim_{j\to+\infty} {\mathcal{F}}(w_{k_j}).$$
Also, for any given~$j$,
we can find a recovery sequence~$w_{j,m}$ that converges to~$w_{k_j}$ in~$X$ as~$m\to+\infty$ and
such that $$
\lim_{m\to+\infty} {\mathcal{F}}_{\e_m} (w_{j,m})={\mathcal{F}}(w_{k_j}).$$
Hence, given~$\ell\in\N$, we pick~$j_\ell\in\N$ such that if~$j\ge j_\ell$ then~$\| w_{k_j}-w\|_X\le\frac1\ell$
and then we pick~$m_\ell\in\N$ such that if~$m\ge m_\ell$ then~$\|w_{j_\ell,m}-w_{k_{j_\ell}}\|_X\le\frac1\ell$ and
$$ \big|{\mathcal{F}}_{\e_m} (w_{j_\ell,m})-{\mathcal{F}}(w_{k_{j_\ell}})\big|\le\frac1\ell.$$
This construction gives that~$\|w_{j_\ell,m_\ell}-w\|_X\le\frac2\ell$ and consequently~$w_{j_\ell,m_\ell}\to w$ in~$X$
as~$\ell\to+\infty$. As a result, by~\eqref{REQ:LIMINF},
\begin{eqnarray*}
{\mathcal{F}}(w)&\le&
\liminf_{\ell\to+\infty} {\mathcal{F}}_{\e_{m_\ell}}(w_{j_\ell,m_\ell})
\\&\le&\liminf_{\ell\to+\infty}\left( {\mathcal{F}}(w_{k_{j_\ell}})+\frac1\ell\right)
\\&=& \liminf_{k\to+\infty} {\mathcal{F}}(w_{k}),
\end{eqnarray*}
which proves the lower semicontinuity property mentioned in~\eqref{MSJNDNTIONDS}.\medskip

See~\cite{MR1968440} and the references therein for a thorough induction to~$\Gamma$-convergence
and related topics.

One of the chief achievements of the $\Gamma$-convergence theory consists precisely
in the correct limit assessment of the singular perturbation problem posed by the Allen-Cahn equation in~\eqref{EPSACAN}. Namely, as established in~\cite{MR473971}, we have that:

\begin{theorem}\label{MORTOLA}
The functional
\begin{equation}\label{FEPSI} {\mathcal{F}}_\e(u):=\int_\Omega \left(\frac{\e\,|\nabla u(x)|^2}2+\frac{W(u(x))}\e\right)\,dx\end{equation}
$\Gamma$-converges as~$\e\searrow 0$ to
\begin{equation}\label{ILPERO} {\mathcal{F}}(u):=\begin{dcases}
c\,\Per(E,\Omega) & \begin{matrix*}[l]&{\mbox{ if $u=\chi_E-\chi_{\R^n\setminus E}$}} \\ &{\mbox{ for some set~$E$ of finite perimeter,}}\end{matrix*}\\
\\
+\infty & \begin{matrix*}[l]&{\mbox{ otherwise,}}\end{matrix*}
\end{dcases}\end{equation}
where
$$ c:=\int_{-1}^1 \sqrt{2W(r)}\,dr.$$
\end{theorem}

The notion of perimeter in~\eqref{ILPERO} is the classical one induced by the functions of bounded variations,
see e.g.~\cite{MR775682}.

A deep variant of Theorem~\ref{MORTOLA} deals with the case of
solutions of~\eqref{EPSACAN} which are not necessarily local minimizers for the energy functional:
this analysis is carried out in~\cite{MR1803974},
in which it is shown that the phase interface converges to a suitably ``generalized'' minimal hypersurface
(possibly counting the ``multiplicity'' of the layers produced by non-minimal solutions).

\subsection{Density estimates and geometric convergence}
Another extremely useful variant of Theorem~\ref{MORTOLA} consists in a ``geometric''
convergence results for the level sets of the minimizers of~\eqref{FEPSI}, stating, roughly speaking,
that if~$u_\e$ is a minimizer of~\eqref{FEPSI}, then its level sets approach locally uniformly the limit
interface. To state this result, which was obtained in~\cite{MR1310848},
we make the notation more precise by setting~${\mathcal{F}}_\e(u,\Omega):={\mathcal{F}}_\e(u)$
whenever we want to emphasize the dependence on the domain~$\Omega$ of the functional in~\eqref{FEPSI},
and we have that:

\begin{theorem}\label{THM:CC}
Assume that~$u_\e$ is a local minimizer for the functional~${\mathcal{F}}_\e$
in~\eqref{FEPSI} in the ball~$B_{1+\e}$.

Then:
\begin{itemize}
\item There exists~$C>0$, depending only on~$n$ and~$W$, such that
\begin{equation}\label{DENESCC-0} {\mathcal{F}}_\e(u_\e,B_1)\le C.\end{equation}
\item Up to a subsequence,
\begin{equation}\label{DENESCC-00}
{\mbox{$u_\e\to\chi_E-\chi_{\R^n\setminus E}$ as~$\e\searrow0$ in~$L^1(B_1)$}}\end{equation}
and the set~$E$ has locally minimal perimeter in~$B_1$.
\item Given~$\vartheta_1$, $\vartheta_2\in(-1,1)$, if~$u_\e(0)>\vartheta_1$, then
\begin{equation}\label{DENESCC-1} \big|\{u_\e>\vartheta_2\}\cap B_r\big|\ge cr^n,\end{equation}
as long as~$r\in(0,1]$ and~$\e\in(0, c_\star\,r]$, where~$c>0$ depends only on~$n$ and~$W$
and~$c_\star>0$ depends only on~$n$, $W$, $\vartheta_1$ and~$\vartheta_2$.
\item Similarly, given~$\vartheta_1$, $\vartheta_2\in(-1,1)$, if~$u_\e(0)<\vartheta_1$, then
\begin{equation}\label{DENESCC-2} \big|\{u_\e<\vartheta_2\}\cap B_r\big|\ge cr^n,\end{equation}
as long as~$r\in(0,1]$ and~$\e\in(0, c_\star\,r]$, where~$c>0$ depends only on~$n$ and~$W$
and~$c_\star>0$ depends only on~$n$, $W$, $\vartheta_1$ and~$\vartheta_2$.
\item For every~$\vartheta\in(0,1)$, the set~$\{|u_\e|<\vartheta\}$ approaches~$\partial E$
locally uniformly as~$\e\searrow0$: more explicitly, given~$r_0\in(0,1)$ and~$\delta>0$
there exists~$\e_0>0$ such that if~$\e\in(0,\e_0)$ then
\begin{equation}\label{DENESCC-3} \{|u_\e|<\vartheta\}\cap B_{r_0}\subseteq \bigcup_{p\in\partial E} B_{\delta}(p).\end{equation}
\end{itemize}
\end{theorem}

A direct, but important, consequence of Theorem~\ref{THM:CC} is that the interface
of a phase transition behaves ``like a codimension one'' set, at least in terms of density estimates:
more specifically, given~$\vartheta\in(0,1)$, if~$u_\e(0)\in(-\vartheta,\vartheta)$, then,
when~$r\in(0,1]$ and~$\e\in(0, c_\star\,r]$,
\begin{eqnarray}
&& \label{CONSEG:1} \big|\{|u_\e|<\vartheta\}\cap B_{r}\big|\le C\e r^{n-1}\\
\label{CONSEG:2} {\mbox{and }}&&\min\Big\{
\big|\{u_\e>\vartheta\}\cap B_{r}\big|,\;\big|\{u_\e<-\vartheta\}\cap B_{r}\big|
\Big\}\ge cr^{n}.
\end{eqnarray}
To check~\eqref{CONSEG:1} and~\eqref{CONSEG:2}, one can argue as follows.
Taking~$\vartheta_1:=-\vartheta$ and~$\vartheta_2:=\vartheta$, one deduces from~\eqref{DENESCC-1} that~$\big|\{u_\e>\vartheta\}\cap B_r\big|\ge cr^n$. Similarly,
taking~$\vartheta_1:=\vartheta$ and~$\vartheta_2:=-\vartheta$, one deduces from~\eqref{DENESCC-2} that~$\big|\{u_\e<-\vartheta\}\cap B_r\big|\ge cr^n$. These observations lead to~\eqref{CONSEG:2}.

Additionally, setting~$\e':=\frac{\e}{r}$ and~$
v_{\e'}(x):=u_\e(rx)$, we have that~$v_{\e'}$ is a local minimizer of the functional~${\mathcal{F}}_{\e'}$
in the ball~$B_{\frac{1+\e}r}\supseteq B_{1+\e'}$. As a consequence, by~\eqref{DENESCC-0},
\begin{equation}\label{IPMAinAP234a4r3tJS} \begin{split}&
C\ge {\mathcal{F}}_{\e'}(v_{\e'},B_1)=
\int_{B_1}\left(\frac{\e'\,r^2\,|\nabla u_\e(rx)|^2}2+\frac{W(u_\e(rx))}{\e'}\right)\,dx\\&\qquad\quad
=\frac1{r^{n-1}}\int_{B_r}\left(\frac{\e\,|\nabla u_\e(y)|^2}2+\frac{W(u_\e(y))}\e\right)\,dy.
\end{split}
\end{equation}
In particular,
\begin{eqnarray*} C&\ge&
\frac{1}{\e r^{n-1}}\int_{\{|u_\e|<\vartheta\}\cap B_{r}}W(u_{\e}(y))\,dy
\\& \ge&\frac{1}{\e r^{n-1}}\,\min_{[-\vartheta,\vartheta]} W\;\big|
\{|u_\e|<\vartheta\}\cap B_{r}\big|,
\end{eqnarray*}
from which~\eqref{CONSEG:1} readily follows, up to renaming~$C$.\medskip

\begin{figure}[h]
\includegraphics[height=0.45\textwidth]{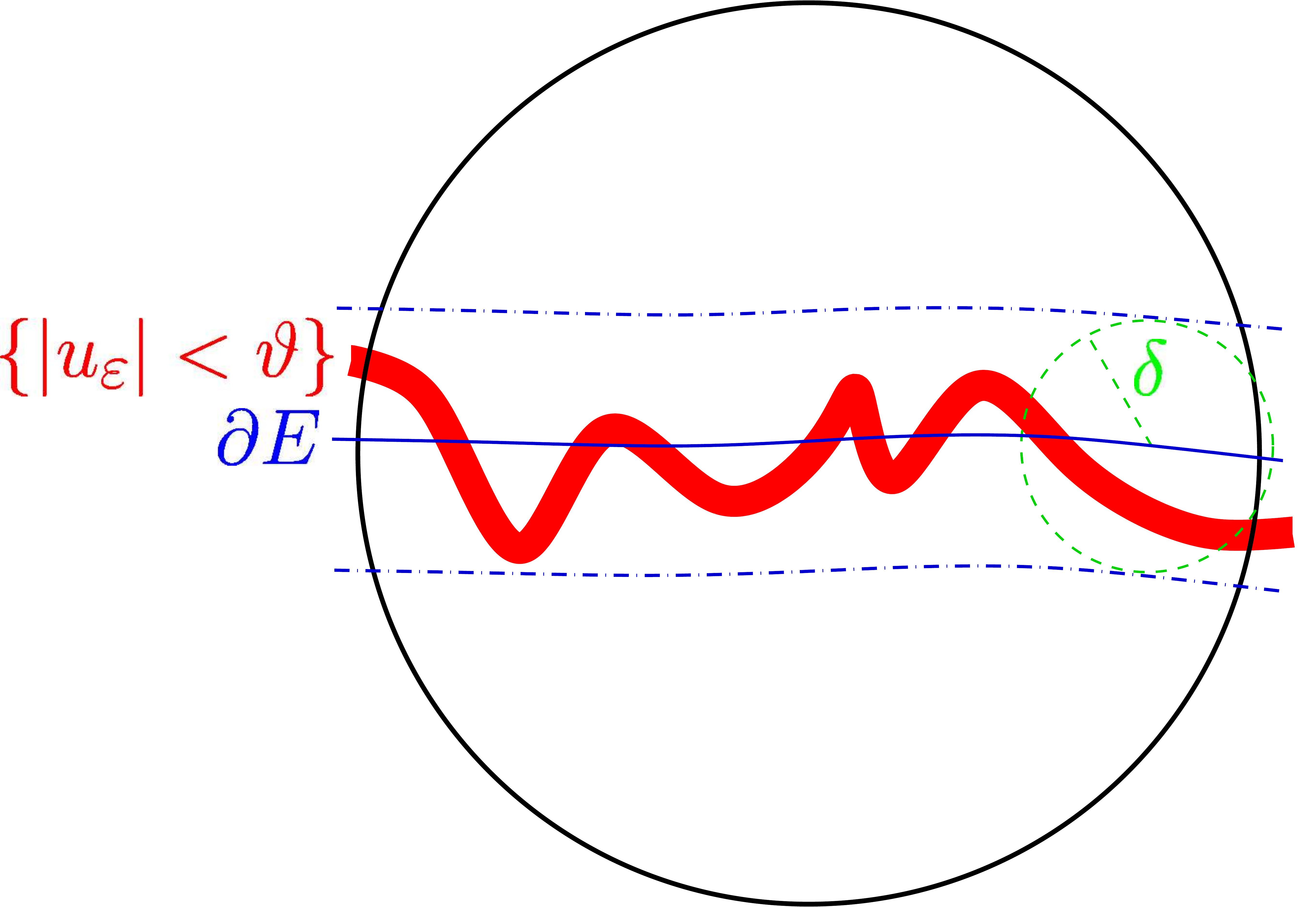}
\caption{Uniform convergence of level sets.}
        \label{d19FIGUR267utgbCLEa5ATBORDUNIgIbnNP-090}
\end{figure}

See also Figure~\ref{d19FIGUR267utgbCLEa5ATBORDUNIgIbnNP-090} for a sketch of the result stated in~\eqref{DENESCC-3}.
It is actually instructive to see the proof of~\eqref{DENESCC-3} as a consequence of~\eqref{DENESCC-00}
and~\eqref{CONSEG:2}.
For this, one can argue by contradiction, assuming that~\eqref{DENESCC-3} does not hold true,
thus finding a sequence of~$\e\searrow0$ for which there exists a sequence of points~$q_\e\in\{|u_\e|<\vartheta\}\cap B_{r_0}$ with~$B_{\delta/2}(q_\e)\cap(\partial E)=\varnothing$.
This gives that either~$B_{\delta/2}(q_\e)\subseteq E$ or~$B_{\delta/2}(q_\e)\subseteq\R^n\setminus E$.
We assume the latter to be true, up to swapping~$E$ and~$\R^n\setminus E$.

Moreover, by centering~\eqref{CONSEG:2} at~$q_\e$ we infer that
$$ \min\Big\{
\big|\{u_\e>\vartheta\}\cap B_{r}(q_\e)\big|,\;\big|\{u_\e<-\vartheta\}\cap B_{r}(q_\e)\big|
\Big\}\ge cr^{n},$$
as long as~$r+r_0<\frac{1+r_0}2$. In particular, we can take~$r:=\min\left\{\frac{1-r_0}4,\frac\delta4\right\}$ and deduce from~\eqref{DENESCC-00}
that
\begin{eqnarray*}
0&=&\lim_{\e\searrow0}\int_{B_1} \big|u_\e(x)-\chi_E(x)+\chi_{\R^n\setminus E}(x)\big|\,dx
\\&\ge&\lim_{\e\searrow0}
\int_{\{u_\e>\vartheta\}\cap B_{r}(q_\e)} \big|u_\e(x)-\chi_E(x)+\chi_{\R^n\setminus E}(x)\big|\,dx\\
&=&\lim_{\e\searrow0}
\int_{\{u_\e>\vartheta\}\cap B_{r}(q_\e)\cap(\R^n\setminus E)} \big|u_\e(x)-\chi_E(x)+\chi_{\R^n\setminus E}(x)\big|\,dx\\
&=&\lim_{\e\searrow0}
\int_{\{u_\e>\vartheta\}\cap B_{r}(q_\e)\cap(\R^n\setminus E)} (u_\e(x)+1)\,dx\\
&\ge&\lim_{\e\searrow0}
\int_{\{u_\e>\vartheta\}\cap B_{r}(q_\e)} (\vartheta+1)\,dx\\
&\ge& c(\vartheta+1)r^n.
\end{eqnarray*}
Since this is a contradiction, the proof of~\eqref{DENESCC-3} is complete.\medskip

It is also worth pointing out that~\eqref{CONSEG:1} and~\eqref{CONSEG:2} are essentially optimal,
since the inequalities presented there can be also reversed, up to changing constants.
Indeed, the inequality in~\eqref{CONSEG:2} can be of course reversed up to constants,
since
$$ \min\Big\{
\big|\{u_\e>\vartheta\}\cap B_{r}\big|,\;\big|\{u_\e<-\vartheta\}\cap B_{r}\big|
\Big\}\le |B_r|=|B_1|\,r^n.$$
As for reverting~\eqref{CONSEG:1}, we point out that, given~$\vartheta\in(0,1)$, if~$u_\e(0)\in(-\vartheta,\vartheta)$, then, when~$r\in(0,1]$ and~$\e\in(0, c_\star\,r]$,
\begin{equation}\label{AMSseddet823}
\big|\{|u_\e|<\vartheta\}\cap B_{r}\big|\ge c_o\,\e r^{n-1}\end{equation}
for some~$c_o>0$ depending only on~$n$, $W$ and~$\vartheta$.
To check this, we define
$$ \widetilde u_\e(x):=\begin{dcases}
u_\e(x) & {\mbox{ if }}u_\e(x)\in(-\vartheta,\vartheta),\\
\vartheta & {\mbox{ if }}u_\e(x)\in[\vartheta,+\infty),\\
-\vartheta & {\mbox{ if }}u_\e(x)\in(-\infty,-\vartheta]
\end{dcases}$$
and we let~$\mu$ be the average of~$\widetilde u_\e$ in~$B_r$.
Thus, recalling~\eqref{CONSEG:2}, we see that if~$\mu\le0$ then
\begin{eqnarray*}&&
\int_{B_{r}} |\widetilde u_\e(x)-\mu|\,dx\ge
\int_{\{\widetilde u_\e\ge\vartheta\}\cap B_r} (\widetilde u_\e(x)-\mu)\,dx\\&&\qquad
\ge \vartheta\,\big|\{\widetilde u_\e\ge\vartheta\}\cap B_r\big|=\vartheta\,\big|\{ u_\e\ge\vartheta\}\cap B_r\big|\ge
c \vartheta r^{n}
\end{eqnarray*}
and similarly if~$\mu>0$ then
\begin{eqnarray*}&&
\int_{B_r} |\widetilde u_\e(x)-\mu|\,dx\ge
\int_{\{\widetilde u_\e\le-\vartheta\}\cap B_r} (\mu-\widetilde u_\e(x))\,dx\\&&\qquad
\ge \vartheta\,\big|\{\widetilde u_\e\le-\vartheta\}\cap B_r\big|=\vartheta\,\big|\{u_\e\le-\vartheta\}\cap B_r\big|\ge
c \vartheta r^{n}.
\end{eqnarray*}
Either way, by Poincar\'e Inequality,
\begin{equation}\label{EWMiTFbia3omlrCA}
\int_{B_r} |\nabla \widetilde u_\e(x)|\,dx\ge \frac{c_1}{ r} \int_{B_r} |\widetilde u_\e(x)-\mu|\,dx\ge
c_2 r^{n-1},
\end{equation}
for suitable positive constants~$c_1$ and~$c_2$ depending only on~$n$, $W$ and~$\vartheta$.

Furthermore, using the Cauchy-Schwarz inequality and~\eqref{IPMAinAP234a4r3tJS}, for every~$\Lambda>0$,
\begin{equation*}
\begin{split}
\int_{B_r} |\nabla \widetilde u_\e(x)|\,dx&
=\int_{B_r} |\nabla u_\e(x)|\,\chi_{\{|u_\e|\le\vartheta\}}(x)\,dx\\&
\le\frac12\int_{B_r} \left(\frac{|\nabla u_\e(x)|^2}\Lambda+\Lambda\chi_{\{|u_\e|\le\vartheta\}}^2(x)\right)\,dx\\&
\le \frac{Cr^{n-1}}{\e\Lambda}+\frac{\Lambda}{2}\big| \{|u_\e|\le\vartheta\}\cap B_r\big|.
\end{split}\end{equation*}
{F}rom this and~\eqref{EWMiTFbia3omlrCA} we arrive at
\begin{eqnarray*}
c_2 r^{n-1}\le\frac{Cr^{n-1}}{\e\Lambda}+\frac{\Lambda}{2}\big| \{|u_\e|\le\vartheta\}\cap B_r\big|.
\end{eqnarray*}
Therefore, choosing~$\Lambda:=\frac{2C}{\e c_2}$,
\begin{eqnarray*}
\frac{c_2 r^{n-1}}2\le\frac{C}{\e c_2}\big| \{|u_\e|\le\vartheta\}\cap B_r\big|,
\end{eqnarray*}
from which~\eqref{AMSseddet823} follows, as desired.
\medskip

\begin{figure}[h]
\includegraphics[height=0.45\textwidth]{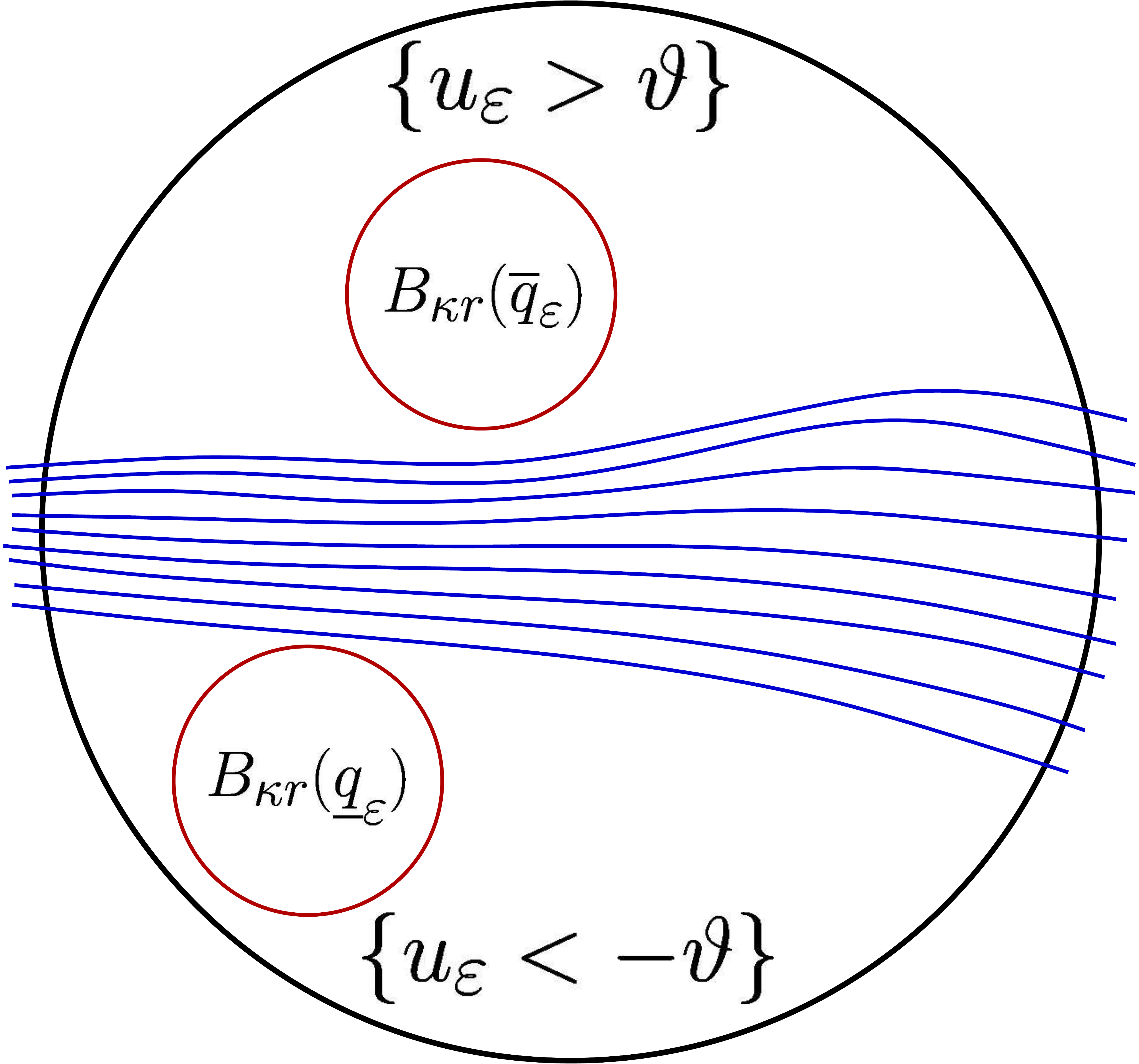}
\caption{The clean ball condition.}
        \label{d19FIGUR267utgbCLEa5ATBORDIbnNP-090}
\end{figure}

\subsection{The clean ball condition}
Another interesting consequence of the previous geometric constructions
is a ``clean ball condition'': namely, looking at a ball centered at the interface,
one can also find balls of comparable size in either side of the interface (hence the interface is not ``spread out''
here and there). The precise result in the spirit of Theorem~\ref{THM:CC}
is that if~$\vartheta\in(0,1)$, $r\in(0,1]$, $\e\in(0, c_\star\,r]$ and~$|u_\e(0)|<\vartheta$
then there exist~$\kappa\in(0,1)$, depending only on~$n$, $W$ and~$\vartheta$,
and points~$\underline{q}_\e$ and~$\overline{q}_\e$ such that
\begin{equation}\label{TYHJSisnSDcP-P}
B_{\kappa r}(\underline{q}_\e)\subseteq \{u_\e<-\vartheta\}\cap B_r\qquad{\mbox{and}}\qquad
B_{\kappa r}(\overline{q}_\e)\subseteq \{u_\e>\vartheta\}\cap B_r.
\end{equation}
See Figure~\ref{d19FIGUR267utgbCLEa5ATBORDIbnNP-090} for a sketch of this phenomenon.
To prove~\eqref{TYHJSisnSDcP-P} one can argue as follows.
Given~$\kappa\in\left(0,\frac1{100}\right)$, we have that\footnote{The choice \label{IJSMSECOfir93}
in~\eqref{TYHJSisnSDcP-P89} will lead to the proof of the first
claim in~\eqref{TYHJSisnSDcP-P}. The proof of the second claim in~\eqref{TYHJSisnSDcP-P},
would focus instead on~$\{u_\e>-\vartheta\}$ in the analog of~\eqref{TYHJSisnSDcP-P89}.}
\begin{equation}\label{TYHJSisnSDcP-P89} \{u_\e<\vartheta\}\cap B_{r/20}\subseteq \bigcup_{p \in\{u_\e<\vartheta\}\cap B_{r/20}} B_{2\kappa r}(p).\end{equation}
By the Vitali Covering Lemma, we can extract a family of disjoint balls~$\big\{B_{2\kappa r}(p_j)\big\}_{j\in{\mathcal{N}}
}$, for some at most countable set of indexes~${\mathcal{N}}$, such that
\begin{equation}\label{TYHJSisnSDcP-P899} \bigcup_{p \in \{u_\e<\vartheta\}\cap B_{r/20}} B_{2\kappa r}(p)\subseteq \bigcup_{j\in{\mathcal{N}}} B_{10\kappa r}(p_j).\end{equation}
By~\eqref{DENESCC-2}, we know that
\[ \big|\{u_\e<\vartheta\}\cap B_{r/20}\big|\ge c r^n\]
up to renaming~$c$,
and consequently, by~\eqref{TYHJSisnSDcP-P89} and~\eqref{TYHJSisnSDcP-P899},
\begin{eqnarray*}&&
cr^n\le\left|\bigcup_{p \in\{u_\e<\vartheta\}\cap B_{r/20}} B_{2\kappa r}(p)\right|\le
\left| \bigcup_{j\in{\mathcal{N}}} B_{10\kappa r}(p_j)\right|\le
\sum_{j\in{\mathcal{N}}} |B_{10\kappa r}(p_j)|=
\kappa^n r^n\,|B_{10}|\,\#{\mathcal{N}},
\end{eqnarray*}
yielding that
\begin{equation}\label{NAMnsikffuNS8iBemowek3}
\#{\mathcal{N}}\ge\frac{\widetilde{c}}{\kappa^n},\end{equation}
for some~${\widetilde{c}}>0$ depending only on~$n$.

Now, let~$\widetilde{\mathcal{N}}$ denote the indexes~$j\in{\mathcal{N}}$ for which~$B_{\kappa r}(p_j)\cap\{|u_\e| \le\vartheta\}\ne\varnothing$.
Accordingly, for each~$j\in\widetilde{\mathcal{N}}$, let us pick a point~$\zeta_j\in B_{\kappa r}(p_j)\cap\{|u_\e| \le\vartheta\}$. 
We stress that if~$x\in B_{\kappa r}(\zeta_j)$ then~$|x-p_j|\le|x-\zeta_j|+|\zeta_j-p_j|<2\kappa r$
and therefore
\begin{equation}\label{P01o2wkedm1masjcnanyhwedfbadfkerfcASGHBX}
B_{\kappa r}(\zeta_j)\subseteq B_{2\kappa r}(p_j).\end{equation}
We also note that, utilizing~\eqref{AMSseddet823},
$$ \big|\{|u_\e|<\vartheta\}\cap B_{\kappa r}(\zeta_j)\big|\ge c_o\,\e \kappa^{n-1} r^{n-1}.$$
Using this, \eqref{P01o2wkedm1masjcnanyhwedfbadfkerfcASGHBX} and the fact that the balls~$B_{2\kappa r}(p_j)$ are disjoint we find that
\begin{eqnarray*}
c_o\,\e \kappa^{n-1} r^{n-1}\;\#{\widetilde{\mathcal{N}}}&\le&\sum_{j\in\widetilde{\mathcal{N}}}\big|\{|u_\e|<\vartheta\}\cap B_{\kappa r}(\zeta_j)\big|\\&\le&\sum_{j\in\widetilde{\mathcal{N}}}\big|\{|u_\e|<\vartheta\}\cap B_{2\kappa r}(p_j)\big|\\&=&\left|\{|u_\e|<\vartheta\}\cap\left(\bigcup_{j\in\widetilde{\mathcal{N}}} B_{2\kappa r}(p_j)\right)\right|\\&\le&
\left|\{|u_\e|<\vartheta\}\cap B_r\right|.
\end{eqnarray*}
This and~\eqref{CONSEG:1} give that
$$ c_o\,\e \kappa^{n-1} r^{n-1}\;\#{\widetilde{\mathcal{N}}}\le C\e r^{n-1}$$ 
and, as a consequence,
$$ \#{\widetilde{\mathcal{N}}}\le\frac{C}{c_o \kappa^{n-1}}.$$
Comparing this with~\eqref{NAMnsikffuNS8iBemowek3} we conclude that, if~$\kappa$ is conveniently small,
$$ \#\big({\mathcal{N}}\setminus\widetilde{\mathcal{N}}\big)\ge\frac{\widetilde{c}}{2\kappa^n}>0.$$
In particular, we can pick~$j_\star\in {\mathcal{N}}\setminus\widetilde{\mathcal{N}}$, yielding that
$$ B_{\kappa r}(p_{j_\star})\cap\{|u_\e| \le\vartheta\}=\varnothing.$$
Since~$u_\e(p_{j_\star})\in \{u_\e<\vartheta\}$ due to~\eqref{TYHJSisnSDcP-P89},
we infer that~$B_{\kappa r}(p_{j_\star})\subseteq\{ u_\e \le-\vartheta\}$.
This establishes the first claim in~\eqref{TYHJSisnSDcP-P}
and the second can be proved similarly (recall footnote~\ref{IJSMSECOfir93}).

\section{Minimal surfaces and one-dimensional symmetry}\label{1243546578679ow3etg245tPOkcv}

The results of Theorems~\ref{MORTOLA} and~\ref{THM:CC}, linking the Allen-Cahn equation
to the theory of minimal surfaces (i.e., in our notation, of hypersurfaces which are boundary of sets
and minimize the perimeter functional under compactly supported perturbations in a given domain)
strongly suggest that the understanding of the locally minimal solutions of the Allen-Cahn equation
and that of perimeter minimizers are deeply related.
In general, solutions of semilinear equations, e.g. equations of the form~$\Delta u=f(u)$ for some~$f:\R\to\R$,
possess the remarkable property that their Laplacian is constant along level sets, thus again suggesting
a strong relation with geometric objects, such as hypersurfaces with constant mean curvature.

\subsection{Classical minimal surfaces}
\begin{figure}[h]
\includegraphics[height=0.45\textwidth]{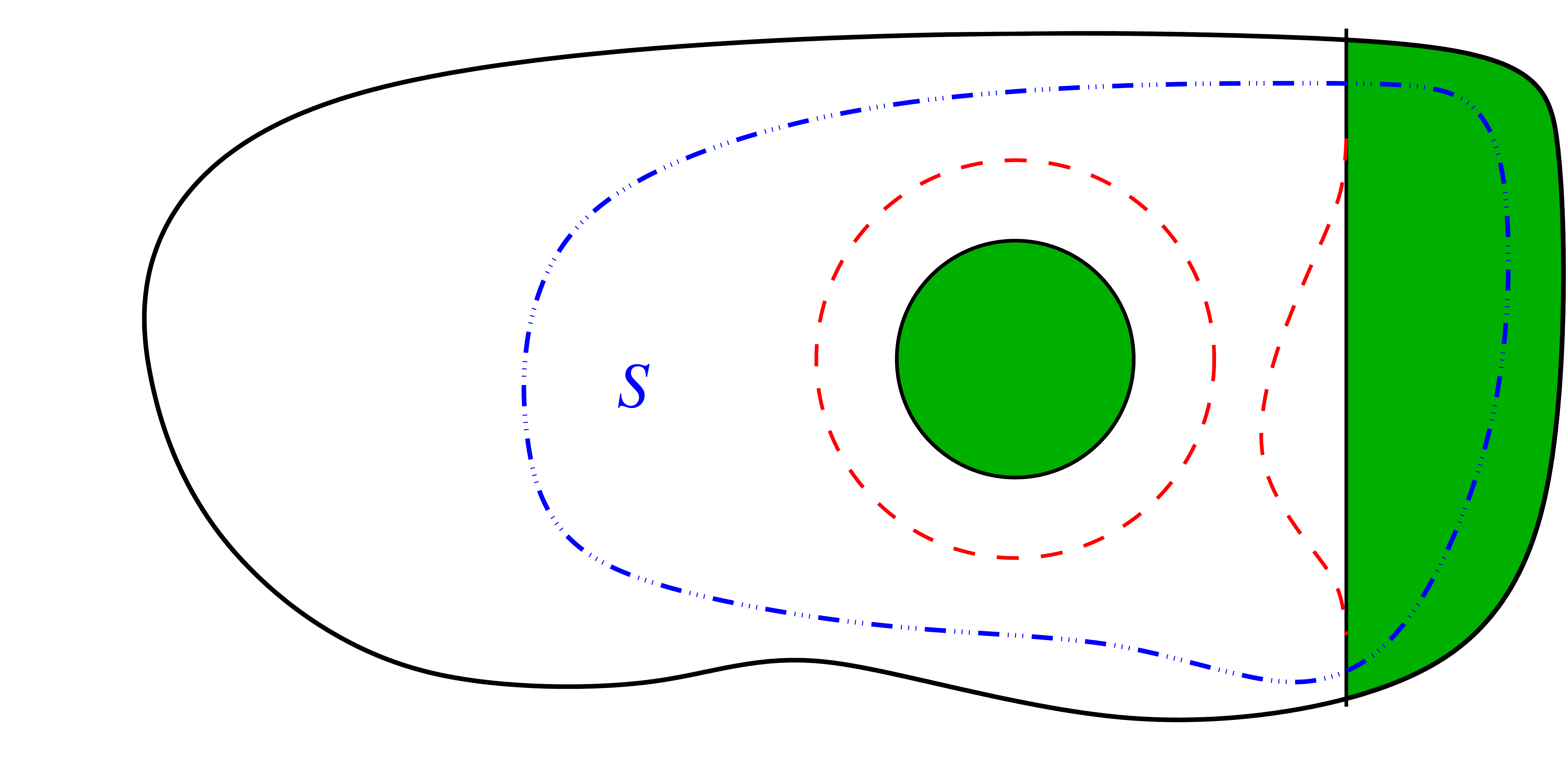}
\caption{Perturbation of a set~$E$ by a normal perturbation supported
in~$S\Subset\Omega$.}
        \label{dOPMperRDUNIgIbn2NP-090}
\end{figure}

To further appreciate the links between analytic and geometric results,
we briefly review some classical aspects of the theory of minimal surfaces.
First of all, the first and second variations of the perimeter functional can be explicitly
computed in terms of the mean curvature~$H$ and of the norm~$c$ of the second fundamental form.
More precisely (see e.g.~\cite[pages~115--120]{MR775682} for details), one can consider a domain~$\Omega\subset\R^n$,
a set~$E\subset\R^n$ and a function~$\phi\in C^\infty_0(\Omega)$
such that~$\partial E$ is a hypersurface of class~$C^2$ in the support~$S$ of~$\phi$.
Thus, we take the exterior normal vector~$\nu$ to~$E$ in~$S$, consider the vector field~$\phi\nu$
(extended to~$0$ outside~$S$). We denote by~$E_t$ the flow of the set~$E$ along this vector field, see Figure~\ref{dOPMperRDUNIgIbn2NP-090}.

Then (see~\cite[equations~(10.12) and~(10.13)]{MR775682}) it follows that, as~$t\to0$,
\begin{equation}\label{MS-odcPSKebad}\begin{split}
\Per(E_t,\Omega)=&\Per(E,\Omega)+t\int_{\partial E}H(x)\,\phi(x)\,d{\mathcal{H}}^{n-1}_x\\&\qquad +
\frac{t^2}2\int_{\partial E}\Big( |\nabla_T\phi(x)|^2-(c^2(x)-H^2(x))\phi^2(x)\Big)\,d{\mathcal{H}}^{n-1}_x+o(t^2),\end{split}
\end{equation}
where the ``tangential gradient'' is given by
\begin{equation}\label{UAMejnfEEthIA} \nabla_T\phi=\nabla\phi-(\nabla\phi\cdot\nu)\nu.\end{equation}
Therefore, by~\eqref{MS-odcPSKebad}, a critical point for the perimeter functional
is (regularity allowing) a hypersurface with vanishing mean curvature
and a minimizer satisfies additionally that
\begin{equation}\label{MS-odcPSKebad-ST}
\int_{\partial E}\Big( |\nabla_T\phi(x)|^2- c^2(x)\,\phi^2(x)\Big)\,d{\mathcal{H}}^{n-1}_x\geq0
\end{equation}
for every test function~$\phi\in C^\infty_0(\Omega)$.\medskip

Condition~\eqref{MS-odcPSKebad-ST} found its own place in the mathematical literature,
and it is indeed customary to say that a vanishing mean curvature hypersurface is ``stable''
if~\eqref{MS-odcPSKebad-ST} holds true (in particular, local minimizers of the perimeter
are stable).

One of the chief results in the classical theory of minimal surfaces
is that \begin{equation}\label{0qwoifk0OSJHDN-OKSDND8-UJHS}
{\mbox{perimeter minimizers are smooth in dimension~$n\le7$.}}\end{equation}
We stress that the dimensional assumption in~\eqref{0qwoifk0OSJHDN-OKSDND8-UJHS}
is optimal, since 
\begin{equation}\label{0qwoifk0OSJHDN-OKSDND8-UJHS-x2}
{\mbox{singular minimal cones occur in dimension~$n\ge8$,}}\end{equation}
as constructed in~\cite{MR250205}, see also e.g.~\cite{MR308905, MR331197, MR1356726} for additional
examples.

The result in~\eqref{0qwoifk0OSJHDN-OKSDND8-UJHS} is proved by reducing, after a blow-up procedure (see~\cite[Theorem~9.3 and Corollary~9.5]{MR775682})
and a dimensional
reduction (see~\cite[Theorem~9.10]{MR775682}), to the case in which the minimal surface~$E$ under consideration is a cone
(namely, if~$p\in E$ then~$tp\in E$ for all~$t>0$)
and its only possible singularity is at the origin.
Thus, in this setting, the claim in~\eqref{0qwoifk0OSJHDN-OKSDND8-UJHS}
is a consequence of a celebrated result in~\cite{MR233295}, according to which
\begin{equation}\label{0qwoifk0OSJHDN-OKSDND8-UJHS-b}\begin{split}&
{\mbox{a vanishing mean curvature cone in~$\R^n$ with~$n\le7$}}\\ &{\mbox{which is stable, and smooth outside the origin,}}\\ &{\mbox{is necessarily a halfplane.}}\end{split}\end{equation}
To prove~\eqref{0qwoifk0OSJHDN-OKSDND8-UJHS-b} we will show that if~$n\le7$ then~$c=0$: correspondingly, all the principal curvatures
of the cone must vanish identically, and therefore the cone must necessarily be a halfplane.

More precisely, the proof of~\eqref{0qwoifk0OSJHDN-OKSDND8-UJHS-b} relies on a beautiful
inequality of geometric type valid for all cones~$E$ with zero mean curvature at regular points, stating that
\begin{equation}\label{0qwoifk0OSJHDN-OKSDND8-UJHS-c}
\frac{\Delta_T c^2}{2} \ge |\nabla_T c|^2+\frac{2c^2}{|x|^2}-c^4,
\end{equation}
where~$\Delta_T$
is the Laplace-Beltrami operator, which can be defined, for instance, in the distributional sense
via the tangential gradient, for smooth and compactly supported functions~$f$ and~$g$, by
\begin{equation}\label{KSLABELTRA} \int_{\partial E} \Delta_T f(x)\,g(x)\,d{\mathcal{H}}^{n-1}_x=
-\int_{\partial E} \nabla_T f(x)\cdot\nabla_T g(x)\,d{\mathcal{H}}^{n-1}_x.\end{equation}
See e.g.~\cite{dipierro2021elliptic} for further details on the Laplace-Beltrami operator.

We postpone the proof of~\eqref{0qwoifk0OSJHDN-OKSDND8-UJHS-c} and we use it
now to obtain~\eqref{0qwoifk0OSJHDN-OKSDND8-UJHS-b}.

\begin{proof}[Proof of the statement in~\eqref{0qwoifk0OSJHDN-OKSDND8-UJHS-b}]
We consider a test function~$\zeta$ and exploit the stability inequality in~\eqref{MS-odcPSKebad-ST} with~$\phi:=c\zeta$, finding that
\begin{eqnarray*}0&\le&
\int_{\partial E}\Big( |\zeta(x)\nabla_T c(x)+c(x)\nabla_T\zeta(x)|^2- c^4(x)\,\zeta^2(x)\Big)\,d{\mathcal{H}}^{n-1}_x\\&=&
\int_{\partial E}\left( |\zeta(x)\nabla_T c(x)|^2+|c(x)\nabla_T\zeta(x)|^2
+\frac12\nabla_T c^2(x)\cdot\nabla_T\zeta^2(x)
- c^4(x)\,\zeta^2(x)\right)\,d{\mathcal{H}}^{n-1}_x\\&=&
\int_{\partial E}\left( \zeta^2(x)|\nabla_T c(x)|^2+c^2(x)|\nabla_T\zeta(x)|^2
-\frac{\Delta_T c^2(x) }2\zeta^2(x)
- c^4(x)\,\zeta^2(x)\right)\,d{\mathcal{H}}^{n-1}_x.
\end{eqnarray*}
Combining this and~\eqref{0qwoifk0OSJHDN-OKSDND8-UJHS-c}, we infer that
\begin{equation}\label{trfsdgbST3P}
\int_{\partial E} \frac{2c^2(x)\,\zeta^2(x)}{|x|^2}\,d{\mathcal{H}}^{n-1}_x
\le \int_{\partial E}c^2(x)\,|\nabla_T\zeta(x)|^2\,d{\mathcal{H}}^{n-1}_x.
\end{equation}

Now, given~$\alpha$, $\beta\in\R$, and~$\e\in(0,1)$ to be taken as small as we wish\footnote{Roughly speaking, the idea now is to take as a test function in~\eqref{trfsdgbST3P} something with
two different behaviors near zero and at infinity, such as
$$ \begin{dcases}
|x|^{\alpha}&{\mbox{if }}x\in B_1,\\
|x|^{\alpha+\beta}&{\mbox{if }}x\in \R^n\setminus B_1.
\end{dcases}$$
As a matter of fact, this choice would formally lead to~\eqref{098u09uyt9GTYbIlaBO09okVOE}.
However, one needs an approximation and cutoff argument in order to make the
above test function rigorously admissible, which is the reason for introducing the parameter~$\e$ here.}
in what follows, we consider~$\tau_\e\in C^\infty_0\big(B_{2/\e}\setminus B_{\e},\,[0,1]\big)$ with~$\tau_\e=1$ in~$B_{1/\e}\setminus B_{2\e}$
with
\begin{equation}\label{123PLSKMskcijfghiOSIJHB}
|\nabla\tau_\e|\le\frac4\e\chi_{B_{2\e}\setminus B_\e}+4\e\chi_{B_{2/\e}\setminus B_{1/\e}}.\end{equation}
Let also~$\zeta_\e:=\tau_\e \varphi_\e $, with
$$ \varphi_\e(x):=\frac{|x|^\alpha}{2}\,\Big( 
\sqrt{\big(|x|^\beta-1\big)^2+\e}+|x|^\beta+1\Big).$$
The idea is to use~$\zeta_\e$ as a test function in~\eqref{trfsdgbST3P}
and pass to the limit as~$\e\searrow0$. For this  approximation method to work,
we will need to choose appropriately the parameters~$\alpha$ and~$\beta$, which, in turn, will be possible only under the dimensional restriction in~\eqref{0qwoifk0OSJHDN-OKSDND8-UJHS-b}. For this, we note that
\begin{eqnarray*}&&\lim_{\e\searrow0} \zeta_\e(x)=
\lim_{\e\searrow0} \varphi_\e(x)=\frac{|x|^\alpha}{2}\,\Big( 
\big||x|^\beta-1\big|+|x|^\beta+1\Big)=\begin{dcases}
|x|^{\alpha+\beta} & {\mbox{ if }}x\in\R^n\setminus B_1,\\
|x|^\alpha &{\mbox{ if }}x\in B_1.
\end{dcases}
\end{eqnarray*}
Thus, Fatou's Lemma entails that
\begin{equation}\label{098uytgujdYUHS-24023P}\begin{split}&
\liminf_{\e\searrow0}
\int_{\partial E} \frac{2c^2(x)\,\zeta^2_\e(x)}{|x|^2}\,d{\mathcal{H}}^{n-1}_x\\ &\qquad\ge
\int_{(\partial E)\cap B_1} 2c^2(x)\,|x|^{2\alpha-2}\,d{\mathcal{H}}^{n-1}_x+
\int_{(\partial E)\setminus B_1} 2c^2(x)\,|x|^{2\alpha+2\beta-2}\,d{\mathcal{H}}^{n-1}_x.
\end{split}
\end{equation}
Furthermore, since~$E$ is a cone, its curvatures are positively homogeneous of degree~$-1$
and therefore, for all~$x\in\R^n\setminus\{0\}$,
\begin{equation}\label{CEST023}
|c(x)|=\frac{\left|c\left(\frac{x}{|x|}\right)\right|}{|x|}\le\frac{M}{|x|}, \qquad{\mbox{where}}\qquad
M:=\displaystyle\max_{(\partial E)\cap(\partial B_1)}|c|.\end{equation}
As a result,
\begin{eqnarray*}
&&\int_{(\partial E)\cap B_1} 2c^2(x)\,|x|^{2\alpha-2}\,d{\mathcal{H}}^{n-1}_x+
\int_{(\partial E)\setminus B_1} 2c^2(x)\,|x|^{2\alpha+2\beta-2}\,d{\mathcal{H}}^{n-1}_x\\&\le&
2M^2
\left[
\int_{(\partial E)\cap B_1} |x|^{2\alpha-4}\,d{\mathcal{H}}^{n-1}_x+
\int_{(\partial E)\setminus B_1} |x|^{2\alpha+2\beta-4}\,d{\mathcal{H}}^{n-1}_x
\right]\\&=&
2M^2
\left[
\int_{0}^1 \rho^{2\alpha-4}\,
{\mathcal{H}}^{n-2}\big((\partial E)\cap (\partial B_\rho)\big)\,d\rho+
\int^{+\infty}_1 \rho^{2\alpha+2\beta-4}\,
{\mathcal{H}}^{n-2}\big((\partial E)\cap (\partial B_\rho)\big)\,d\rho\right]\\&=&
2M^2 {\mathcal{H}}^{n-2}\big((\partial E)\cap (\partial B_1)\big)
\left[
\int_{0}^1 \rho^{2\alpha+n-6}\,+
\int^{+\infty}_1 \rho^{2\alpha+2\beta+n-6}\,d\rho\right]\\&=&
2M^2 {\mathcal{H}}^{n-2}\big((\partial E)\cap (\partial B_1)\big)
\left[\frac{1}{2\alpha+n-5}-\frac1{2\alpha+2\beta+n-5}\right]\\&<&+\infty,
\end{eqnarray*}
as long as
\begin{equation}\label{ijn6789YSHINS-wef}
2\alpha+n-5>0 \qquad{\mbox{and}}\qquad2\alpha+2\beta+n-5<0.
\end{equation}
In particular, these conditions guarantee that the integrals in
the right hand side of~\eqref{098uytgujdYUHS-24023P} are finite.

Additionally, we claim that
\begin{equation}\label{COkmOOS00PKzx2345mc}
|\nabla\varphi_\e(x)|\le C\max\{|x|^{\alpha+\beta-1},|x|^{\alpha-1}\},
\end{equation}
for some~$C>0$ depending only on~$\alpha$ and~$\beta$, as long as~$\e>0$ is small enough. Indeed,
using the notation~$r:=|x|$, we have that
\begin{equation}\label{MSnmsL97VhGYiMOal}
\begin{split}&
|\nabla\varphi_\e(x)|\,=\,
\frac{r^{\alpha-1}}{2}
\Bigg|(\alpha+\beta)\left( \sqrt{(r^\beta - 1)^2 + \e } + r^\beta+1\right) \\&\qquad\qquad\qquad\qquad\qquad-\beta
\frac{\sqrt{(r^\beta - 1)^2 + \e}  -r^\beta+1+\e}{\sqrt{(r^\beta - 1)^2 + \e} }\Bigg|\end{split}
\end{equation}
and~\eqref{COkmOOS00PKzx2345mc} plainly follows when~$r\in\left[\frac12,2\right]$.
Also, if~$r\in\left(0,\frac12\right)$, we have that
$$ \frac{\sqrt{(r^\beta - 1)^2 + \e}}{1-r^\beta}=
\sqrt{1 +\frac{ \e}{(1-r^\beta)^2}}\in [1-C\e,1+C\e].
$$
Accordingly,
\begin{eqnarray*}
\sqrt{(r^\beta - 1)^2 + \e } + r^\beta+1 &\in&
\big[(1-C\e)(1-r^\beta)+ r^\beta+1,\,(1+C\e)(1-r^\beta)+ r^\beta+1\big]\\&
\subseteq&[2-C\e,2+C\e]
\end{eqnarray*}
and
\begin{eqnarray*}&&
\frac{\sqrt{(r^\beta - 1)^2 + \e}  -r^\beta+1+\e}{\sqrt{(r^\beta - 1)^2 + \e} }-1
=\frac{1-r^\beta+\e}{\sqrt{(r^\beta - 1)^2 + \e} }\\
&&\qquad\in\left[\frac{(1-C\e)(1-r^\beta+\e)}{1-r^\beta},\,\frac{(1+C\e)(1-r^\beta+\e)}{1-r^\beta}\right]\subseteq\left[ 1-C\e, 1+C\e\right],
\end{eqnarray*}
where, as customary,~$C$ has been renamed line after line.

Hence, we deduce from~\eqref{MSnmsL97VhGYiMOal} that
\begin{equation*}
|\nabla\varphi_\e(x)|\le
\frac{ r^{\alpha-1}}2
\left( \big| 2(\alpha+\beta)-2 \big|+C\e\right)\le 
\left( \big| \alpha+\beta-1 \big|+C\e\right) r^{\alpha-1},
\end{equation*}
which proves~\eqref{COkmOOS00PKzx2345mc} in this case.

Finally, if~$r>2$,
$$ \frac{\sqrt{(r^\beta - 1)^2 + \e}}{r^\beta-1}=
\sqrt{1 +\frac{ \e}{(r^\beta-1)^2}}\in [1-C\e,1+C\e].
$$
In this case, we have that
\begin{eqnarray*}
\sqrt{(r^\beta - 1)^2 + \e } + r^\beta+1 &\in&
\big[(1-C\e)(r^\beta-1)+ r^\beta+1,\,(1+C\e)(r^\beta-1)+ r^\beta+1\big]\\&
\subseteq&[(2-C\e)r^\beta,(2+C\e)r^\beta]
\end{eqnarray*}
and
\begin{eqnarray*}&&
\frac{\sqrt{(r^\beta - 1)^2 + \e}  -r^\beta+1+\e}{\sqrt{(r^\beta - 1)^2 + \e} }=
1-\frac{ r^\beta-1-\e}{\sqrt{(r^\beta - 1)^2 + \e} }\in[-C\e,C\e],
\end{eqnarray*}
yielding that
\begin{equation*}
|\nabla\varphi_\e(x)|\le\frac{r^{\alpha-1}}2\left(
\big|2(\alpha+\beta)r^\beta\big|+C\e\right)\le {r^{\alpha+\beta-1}}\left(
\big|(\alpha+\beta)\big|+C\e\right)
.\end{equation*}
Combining these pieces of information and~\eqref{MSnmsL97VhGYiMOal}
we obtain that~\eqref{COkmOOS00PKzx2345mc} holds true in this case as well.

Now, exploiting~\eqref{123PLSKMskcijfghiOSIJHB}, \eqref{CEST023}
and~\eqref{COkmOOS00PKzx2345mc},
and possibly renaming~$C$ in dependence also on~$E$, we see that
\begin{eqnarray*}&&
\left|\,
\int_{\partial E}c^2(x)\,|\nabla_T\zeta_\e(x)|^2\,d{\mathcal{H}}^{n-1}_x-
\int_{\partial E}c^2(x)\,|\nabla_T\varphi_\e(x)|^2\,d{\mathcal{H}}^{n-1}_x
\right|\\&\le&C\left(
\int_{(\partial E)\cap 
\big(B_{2\e}\cup(\R^n\setminus B_{1/\e})\big)
}c^2(x)\,|\nabla\varphi_\e(x)|^2\,d{\mathcal{H}}^{n-1}_x
+
\int_{\partial E}c^2(x)\,|\nabla \tau_\e(x)|^2\,|\varphi_\e(x)|^2\,d{\mathcal{H}}^{n-1}_x\right)\\&\le&
C\Bigg[
\int_{(\partial E)\cap 
\big(B_{2\e}\cup(\R^n\setminus B_{1/\e})\big)
} 
\max\{|x|^{2\alpha+2\beta-4},|x|^{2\alpha-4}\}
\,d{\mathcal{H}}^{n-1}_x
\\&&\qquad\qquad+
\int_{\partial E\cap B_{2\e}} \e^{-2}|x|^{2\alpha-2}\,d{\mathcal{H}}^{n-1}_x+
\int_{\partial E\setminus B_{1/\e}} \e^2|x|^{2\alpha+2\beta-2}\,d{\mathcal{H}}^{n-1}_x
\Bigg]
\\&\le& C\big(\e^{2\alpha+n-5}+ \e^{5-2\alpha-2\beta-n}\big)
,\end{eqnarray*}
which is infinitesimal as long as condition~\eqref{ijn6789YSHINS-wef} is fulfilled.
This, together with~\eqref{trfsdgbST3P}
and~\eqref{098uytgujdYUHS-24023P}, yields that
\begin{equation*}\begin{split}&
\int_{(\partial E)\cap B_1} 2c^2(x)\,|x|^{2\alpha-2}\,d{\mathcal{H}}^{n-1}_x+
\int_{(\partial E)\setminus B_1} 2c^2(x)\,|x|^{2\alpha+2\beta-2}\,d{\mathcal{H}}^{n-1}_x\\
\le\;&
\liminf_{\e\searrow0}
\int_{\partial E} \frac{2c^2(x)\,\zeta^2_\e(x)}{|x|^2}\,d{\mathcal{H}}^{n-1}_x\\
\le\;&\liminf_{\e\searrow0}\int_{\partial E}c^2(x)\,|\nabla_T\zeta_\e(x)|^2\,d{\mathcal{H}}^{n-1}_x
\\=\;&
\liminf_{\e\searrow0}\int_{\partial E}c^2(x)\,|\nabla_T\varphi_\e(x)|^2\,d{\mathcal{H}}^{n-1}_x
\\ \leq\;&
\liminf_{\e\searrow0}\int_{\partial E}c^2(x)\,|\nabla\varphi_\e(x)|^2\,d{\mathcal{H}}^{n-1}_x.\end{split}\end{equation*}
Owing to~\eqref{CEST023}, \eqref{ijn6789YSHINS-wef} and~\eqref{COkmOOS00PKzx2345mc}, we can also
pass the limit inside the integral sign in the latter term, using the Dominated Convergence Theorem, and, recalling~\eqref{MSnmsL97VhGYiMOal}, we thus find that
\begin{equation}\label{098u09uyt9GTYbIlaBO09okVOE}\begin{split}&
\int_{(\partial E)\cap B_1} 2c^2(x)\,|x|^{2\alpha-2}\,d{\mathcal{H}}^{n-1}_x+
\int_{(\partial E)\setminus B_1} 2c^2(x)\,|x|^{2\alpha+2\beta-2}\,d{\mathcal{H}}^{n-1}_x\\
\le\,&\int_{\partial E}c^2(x)\,
\lim_{\e\searrow0}|\nabla\varphi_\e(x)|^2\,d{\mathcal{H}}^{n-1}_x\\ =\,&
\int_{(\partial E)\cap B_1} \alpha^2 c^2(x)\,|x|^{2\alpha-2}\,d{\mathcal{H}}^{n-1}_x
+\int_{(\partial E)\setminus B_1} (\alpha+\beta)^2c^2(x)\,|x|^{2\alpha+2\beta-2}\,d{\mathcal{H}}^{n-1}_x,
\end{split}\end{equation}
that is
\begin{equation}\label{FYSHD:ILSNCt}
\int_{(\partial E)\cap B_1} \kappa_1\,c^2(x)\,|x|^{2\alpha-2}\,d{\mathcal{H}}^{n-1}_x+
\int_{(\partial E)\setminus B_1} \kappa_2\,c^2(x)\,|x|^{2\alpha+2\beta-2}\,d{\mathcal{H}}^{n-1}_x
\le0,\end{equation}
where~$\kappa_1:=2-\alpha^2$
and~$\kappa_2:=2-(\alpha+\beta)^2$.

Our goal is now to exploit the dimension assumption~$n\le7$ in~\eqref{0qwoifk0OSJHDN-OKSDND8-UJHS-b}
in order to fulfill~\eqref{ijn6789YSHINS-wef} and also make~$\kappa_1$ and~$\kappa_2$ strictly positive.
We stress that this would complete the proof of~\eqref{0qwoifk0OSJHDN-OKSDND8-UJHS-b},
since we would deduce from~\eqref{FYSHD:ILSNCt} that~$c$ vanishes identically, giving that~$\partial E$ is a hyperplane.

Thus, it remains to check that it is possible to choose our free parameters~$\alpha$ and~$\beta$ such that
\begin{equation}\label{fier64b785b4v986nbfroew}
\varsigma:=\alpha+\beta<\frac{5-n}2<\alpha,\qquad|\varsigma|<\sqrt2
\qquad{\mbox{and}}\qquad|\alpha|<\sqrt2.
\end{equation}
Notice that if~$3\le n\le7$, then in particular~$-\sqrt2<\frac{5-n}{2}<\sqrt2$, and therefore
it is possible to choose~$\alpha$ and~$\beta$ in such a way that the conditions in~\eqref{fier64b785b4v986nbfroew} are satisfied.
Namely, if~$3\le n\le7$ we can choose
$$ \alpha:=\frac{5-n}4+\frac{\sqrt2}2 \qquad{\mbox{and}}\qquad \beta:=-\sqrt2,
$$
obtaining~\eqref{0qwoifk0OSJHDN-OKSDND8-UJHS-b} in this case.

We remark that when~$n=2$ the statement in~\eqref{0qwoifk0OSJHDN-OKSDND8-UJHS-b}
is clearly satisfied, since to rule out the possibility that a minimal cone is nontrivial
in the plane it is enough to take a straight line which has less perimeter. 
This completes the proof of~\eqref{0qwoifk0OSJHDN-OKSDND8-UJHS-b}.
\end{proof}

Now we provide a proof of~\eqref{0qwoifk0OSJHDN-OKSDND8-UJHS-c}. 

\begin{proof}[Proof of~\eqref{0qwoifk0OSJHDN-OKSDND8-UJHS-c}]
We start by rewriting~\eqref{UAMejnfEEthIA}
in coordinates as
\begin{equation}\label{DEDAS} \partial_{T,j}\phi=\partial_j\phi-(\nabla\phi\cdot\nu)\nu_j,\qquad{\mbox{ for }}\;j\in\{1,\dots,n\}.\end{equation}
We stress that~$\partial_{T,j}$ plays the role of a ``tangential derivative'', since~$\tau_j:=e_j-\nu_j\,\nu$
is a tangent vector to~$\partial E$ (because~$\tau_j\cdot\nu=e_j\cdot\nu-\nu_j=0$) and~$\nabla\phi\cdot\tau_j=\partial_{T,j}\phi$.
Note also that~$\nabla_T\phi$ plays the role of a tangent vector, since
\begin{equation} \label{G1014}\nu\cdot\nabla_T\phi=\sum_{\ell=1}^n \nu_\ell \,\partial_{T,\ell}\phi=\sum_{\ell=1}^n
\nu_\ell\Big(\partial_\ell\phi-(\nabla\phi\cdot\nu)\nu_\ell\Big)=\sum_{\ell=1}^n
\nu_\ell\partial_\ell\phi-(\nabla\phi\cdot\nu)=0.\end{equation}

Furthermore, we extend~$\nu$ to a neighborhood of~$\partial E$,
say by a normal extension according to which,
for all~$x\in\partial E$ and~$t\in\R$ with~$|t|$ close enough,
\begin{equation}\label{NE55}
\nu(x+t\nu(x))=\nu(x).\end{equation}
In this way, for each~$i\in\{1,\dots,n\}$,
\begin{equation}\label{NE56}\begin{split}&0=\left.\frac{d}{dt} \nu_i(x)\right|_{t=0}
\left.=\frac{d}{dt} \nu_i(x+t\nu(x))\right|_{t=0}=\sum_{j=1}^n\partial_j\nu_i(x)\,\nu_j(x).\end{split}
\end{equation}

Hence, since~$c$ is the norm of the second fundamental form~$\nabla_T\nu$, we have that
\begin{equation}\label{REGVD:LL} c^2=\sum_{i,j=1}^n(\partial_{T,j}\nu_i)^2.\end{equation}
At this stage, it is also convenient to observe that the Laplace-Beltrami operator
can be reconstructed from the tangential derivatives through the formula
\begin{equation}\label{BELL7}
\Delta_T f=\sum_{i=1}^n\partial_{T,i}\partial_{T,i}f.
\end{equation}
Not to interrupt the flow of the argument, we defer the proof of~\eqref{BELL7}
to Appendix~\ref{appeadd10}.

Moreover, recalling~\eqref{KSLABELTRA},
\begin{eqnarray*}
&&\int_{\partial E} \Big(\Delta_T f^2(x)-2|\nabla_Tf(x)|^2-2f(x)\Delta_T f(x)\Big)\,g(x)\,d{\mathcal{H}}^{n-1}_x\\&=&
-\int_{\partial E} \nabla_T f^2(x)\cdot\nabla_T g(x)\,d{\mathcal{H}}^{n-1}_x
-2\int_{\partial E}|\nabla_Tf(x)|^2\,g(x)\,d{\mathcal{H}}^{n-1}_x\\&&\qquad
+2\int_{\partial E} \nabla_T f(x)\cdot\nabla_T(f(x)g(x))\,d{\mathcal{H}}^{n-1}_x\\&=&0,\end{eqnarray*}
that is
$$ \frac{\Delta_T f^2}2=|\nabla_Tf|^2+f\Delta_T f.$$
{F}rom these observations, we arrive at
\begin{equation}\label{0-0-32edxCOMMUTATA}
\begin{split}&
\frac{\Delta_T c^2}2=\frac12\sum_{i,j=1}^n\Delta_T(\partial_{T,j}\nu_i)^2
=\sum_{i,j=1}^n\Big[|\nabla_T \partial_{T,j}\nu_i|^2+
\partial_{T,j}\nu_i
\Delta_T(\partial_{T,j}\nu_i)
\Big]\\&\qquad\qquad\qquad=\sum_{i,j,m=1}^n(\partial_{T,m} \partial_{T,j}\nu_i)^2+
\sum_{i,j=1}^n \partial_{T,j}\nu_i\Delta_T(\partial_{T,j}\nu_i).
\end{split}\end{equation}

We now recall two useful commutator identities for tangential derivatives, namely, for each~$i$, $j\in\{1,\dots,n\}$,
\begin{equation}\label{COMMUTATA0}
\partial_{T,i}\nu_j=\partial_{T,j}\nu_i
\end{equation}
and
\begin{equation}\label{COMMUTATA}
\partial_{T,i}\partial_{T,j}\phi-
\partial_{T,j}\partial_{T,i}\phi=\sum_{k=1}^n
\Big[\nu_i\partial_{T,j}\nu_k-\nu_j\partial_{T,i}\nu_k\Big]\partial_{T,k}\phi.
\end{equation}
The proofs of~\eqref{COMMUTATA0} and~\eqref{COMMUTATA} are deferred to Appendix~\ref{appeadd10}.

It is also useful to observe that,
for all~$i\in\{1,\dots,n\}$,
\begin{equation}\label{CABNSM:IAKMSbbMNs-1}
\sum_{k=1}^n\nu_k\partial_{T,i}\nu_k
=\frac12\sum_{k=1}^n\partial_{T,i}\nu_k^2
=\frac12\partial_{T,i}|\nu|^2=\frac12\partial_{T,i}1=0.
\end{equation}
Moreover, we point out that
$$ H=\div_T\nu,
$$
where~$\div_T$ denotes the tangential divergence of a field, see formula~\eqref{diuewyt8tuy9854674967698708}.

Hence, using~\eqref{BELL7},
\eqref{COMMUTATA0}, \eqref{COMMUTATA} and~\eqref{CABNSM:IAKMSbbMNs-1},
for each~$j\in\{1,\dots,n\}$ we have that
\begin{equation}\label{G1017}\begin{split}&
\Delta_T\nu_j+c^2\nu_j-\partial_{T,j} H
=\sum_{i=1}^n\partial_{T,i}\partial_{T,i}\nu_j+c^2\nu_j-\sum_{i=1}^n\partial_{T,j} \partial_{T,i}\nu_i\\
&\qquad=\sum_{i=1}^n\partial_{T,i} \partial_{T,j}\nu_i +c^2\nu_j-\sum_{i=1}^n\partial_{T,j} \partial_{T,i}\nu_i\\
&\qquad=\sum_{i=1}^n\partial_{T,j} \partial_{T,i}\nu_i 
+\sum_{i,k=1}^n
\Big[\nu_i\partial_{T,j}\nu_k-\nu_j\partial_{T,i}\nu_k\Big]\partial_{T,k}\nu_i
+c^2\nu_j-\sum_{i=1}^n\partial_{T,j} \partial_{T,i}\nu_i\\&\qquad=
\sum_{i,k=1}^n
\nu_i\partial_{T,j}\nu_k\,\partial_{T,k}\nu_i\\&\qquad=0.
\end{split}\end{equation}

Now, we consider another useful commutator identity, claiming that, for each~$k\in\{1,\dots,n\}$, if~$H=0$ along~$\partial E$ then
\begin{equation}\label{G1018}
\partial_{T,k}(\Delta_T\phi)-\Delta_T(\partial_{T,k}\phi)=
2\sum_{i,j=1}^n \Big[
\nu_k(\partial_{T,i}\nu_j)(\partial_{T,i}\partial_{T,j}\phi)+
(\partial_{T,k}\nu_j)(\partial_{T,j}\nu_i)(\partial_{T,i}\phi)
\Big].\end{equation}
The proof of~\eqref{G1018} is deferred to Appendix~\ref{appeadd10}.

In the light of~\eqref{COMMUTATA0}, \eqref{CABNSM:IAKMSbbMNs-1},
\eqref{G1017} and~\eqref{G1018} we infer that, if~$H=0$ along~$\partial E$,
\begin{equation}\label{5.30}
\begin{split}&
\sum_{i,j=1}^n (\partial_{T,i}\nu_j)\Delta_T(\partial_{T,i}\nu_j)\\=\;&
\sum_{i,j=1}^n (\partial_{T,i}\nu_j)\partial_{T,i}(\Delta_T\nu_j)\\&\quad
-2\sum_{i,j,h,m=1}^n (\partial_{T,i}\nu_j)\Big[
\nu_i(\partial_{T,h}\nu_m)(\partial_{T,h}\partial_{T,m}\nu_j)+
(\partial_{T,i}\nu_m)(\partial_{T,m}\nu_h)(\partial_{T,h}\nu_j)
\Big]\\=\;&-\sum_{i,j=1}^n (\partial_{T,i}\nu_j)\partial_{T,i}(c^2\nu_j)\\&\quad
-2\sum_{i,j,h,m=1}^n (\partial_{T,j}\nu_i)\Big[
\nu_i(\partial_{T,h}\nu_m)(\partial_{T,h}\partial_{T,m}\nu_j)+
(\partial_{T,i}\nu_m)(\partial_{T,m}\nu_h)(\partial_{T,h}\nu_j)
\Big]\\=\;&-\sum_{i,j=1}^n (\partial_{T,i}\nu_j)\partial_{T,i}(c^2\nu_j)
-2\sum_{i,j,h,m=1}^n (\partial_{T,j}\nu_i)(\partial_{T,i}\nu_m)(\partial_{T,m}\nu_h)(\partial_{T,h}\nu_j)\\=\;&-c^4-\sum_{i,j=1}^n \partial_{T,i} c^2\,\nu_j(\partial_{T,i}\nu_j)
-2\sum_{i,j,h,m=1}^n (\partial_{T,j}\nu_i)(\partial_{T,i}\nu_m)(\partial_{T,m}\nu_h)(\partial_{T,h}\nu_j)\\=\;&-c^4
-2\sum_{i,j,h,m=1}^n (\partial_{T,j}\nu_i)(\partial_{T,i}\nu_m)(\partial_{T,m}\nu_h)(\partial_{T,h}\nu_j)
.\end{split}\end{equation}
Furthermore, by~\eqref{COMMUTATA} we know that
\begin{eqnarray*}&&
\partial_{T,j}\partial_{T,h}\nu_m=\partial_{T,h}\partial_{T,j}\nu_m-\sum_{k=1}^n
\Big[\nu_h\partial_{T,j}\nu_k-\nu_j\partial_{T,h}\nu_k\Big]\partial_{T,k}\nu_m
\\ {\mbox{and }}&&
\partial_{T,m}\partial_{T,i}\nu_j=\partial_{T,i}\partial_{T,m}\nu_j-\sum_{\ell=1}^n
\Big[\nu_i\partial_{T,m}\nu_\ell-\nu_m\partial_{T,i}\nu_\ell\Big]\partial_{T,\ell}\nu_j.\end{eqnarray*}
Thus, by~\eqref{G1014},
\begin{eqnarray*}&&
\sum_{h=1}^n\nu_h\,\partial_{T,j}\partial_{T,h}\nu_m=-\sum_{h,k=1}^n
\nu_h^2\,\partial_{T,j}\nu_k\,\partial_{T,k}\nu_m=-\sum_{k=1}^n\partial_{T,j}\nu_k\,\partial_{T,k}\nu_m
\\ {\mbox{and }}&&
\sum_{i=1}^n\nu_i\,\partial_{T,m}\partial_{T,i}\nu_j=-\sum_{i,\ell=1}^n
\nu_i^2\,\partial_{T,m}\nu_\ell\,\partial_{T,\ell}\nu_j=-\sum_{\ell=1}^n
\partial_{T,m}\nu_\ell\,\partial_{T,\ell}\nu_j.\end{eqnarray*}
Taking the product of these two identities, we conclude that
\begin{eqnarray*}&&
\sum_{i,h=1}^n
\nu_i\nu_h(\partial_{T,j}\partial_{T,h}\nu_m)(\partial_{T,m}\partial_{T,i}\nu_j)=
\sum_{k,\ell=1}^n
(\partial_{T,m}\nu_\ell)(\partial_{T,\ell}\nu_j)(\partial_{T,j}\nu_k)(\partial_{T,k}\nu_m).
\end{eqnarray*}
Therefore, in view of~\eqref{5.30}, if~$H=0$ along~$\partial E$,
\begin{eqnarray*}&&
\sum_{i,j=1}^n (\partial_{T,i}\nu_j)\Delta_T(\partial_{T,i}\nu_j)=
-c^4
-2\sum_{i,j,h,m=1}^n
\nu_i\nu_h(\partial_{T,j}\partial_{T,h}\nu_m)(\partial_{T,m}\partial_{T,i}\nu_j)
.\end{eqnarray*}
{F}rom this and~\eqref{0-0-32edxCOMMUTATA} we infer that, if~$H=0$ along~$\partial E$,
\begin{equation}\label{P123}
\frac{\Delta_T c^2}2=\sum_{i,j,m=1}^n(\partial_{T,m} \partial_{T,j}\nu_i)^2-c^4
-2\sum_{i,j,h,m=1}^n
\nu_i\nu_h(\partial_{T,j}\partial_{T,h}\nu_m)(\partial_{T,m}\partial_{T,i}\nu_j).
\end{equation}

Now, to complete the proof of~\eqref{0qwoifk0OSJHDN-OKSDND8-UJHS-c},
we suppose that~$H=0$ along~$\partial E$, we pick a regular point~$x_0\in\partial E$ and we assume, up to a rotation, that
\begin{equation}\label{nuen}
\nu(x_0)=e_n.\end{equation}
Notice that~\eqref{DEDAS} evaluated at~$x_0$ gives that
\begin{equation}\label{PNMRmnAQSECEDDR-1} \partial_{T,n}\phi(x_0)=\partial_n\phi(x_0)-(\nabla\phi(x_0)\cdot e_n) =0.\end{equation}
As a result, by~\eqref{DEDAS}, \eqref{NE55}, \eqref{NE56}, \eqref{COMMUTATA0} and~\eqref{COMMUTATA},
\begin{equation}\label{PNMRmnAQSECEDDR-2}\begin{split}&
\partial_{T,i}\partial_{T,n}\nu_n(x_0)\\ =\,&\partial_{T,n}\partial_{T,i}\nu_n(x_0)+\sum_{k=1}^n
\Big[\nu_i(x_0)\partial_{T,n}\nu_k(x_0)-\nu_n(x_0)\partial_{T,i}\nu_k(x_0)\Big]\partial_{T,k}\nu_n(x_0)
\\=\,&0+\sum_{k=1}^n
\Big[0-\partial_{T,i}\nu_k(x_0)\Big]\partial_{T,k}\nu_n(x_0)\\=\,&
-\sum_{k=1}^n \partial_{T,i}\nu_k(x_0)\,\partial_{T,k}\nu_n(x_0)\\
=\,&-\sum_{k=1}^n
\Big( \partial_i\nu_k(x_0)-(\nabla\nu_k(x_0)\cdot\nu(x_0))\nu_i(x_0)\Big)
\Big( \partial_k\nu_n(x_0)-(\nabla\nu_n(x_0)\cdot\nu(x_0))\nu_k(x_0)\Big)\\
=\,&-\sum_{k=1}^n\partial_i\nu_k(x_0)\, \partial_k\nu_n(x_0)\\
=\,&-\sum_{k=1}^n\partial_i\nu_k(x_0)\, \partial_n\nu_k(x_0)\\
=\,&-\left.\frac12\partial_i\partial_n\left(\sum_{k=1}^n\nu_k^2(x)\right)\right|_{x=x_0}+\sum_{k=1}^n\nu_k(x_0)\partial_n\partial_i\nu_k(x_0)\\=\,&0+\partial_n\partial_i\nu_n(x_0)\\=\,&\partial_n\partial_n\nu_i(x_0)\\
=\,&\lim_{t\to0}\frac{\nu_i(x_0+t\nu(x_0))+\nu_i(x_0-t\nu(x_0))-2\nu_i(x_0)}{t^2}\\=\,&0.
\end{split}\end{equation}
 
Also, owing to~\eqref{P123}, 
\begin{equation}\label{FRTGSDB-00-0234-1}
\frac{\Delta_T c^2(x_0)}2 =
\sum_{i,j,m=1}^n(\partial_{T,m} \partial_{T,j}\nu_i(x_0))^2-c^4(x_0)
-2\sum_{j,m=1}^n
(\partial_{T,j}\partial_{T,n}\nu_m(x_0))(\partial_{T,m}\partial_{T,n}\nu_j(x_0)).\end{equation}
By~\eqref{PNMRmnAQSECEDDR-1} and~\eqref{PNMRmnAQSECEDDR-2}, we know that
$$ \partial_{T,n}\partial_{T,n}\nu_j(x_0)=0=\partial_{T,m}\partial_{T,n}\nu_n(x_0)$$
and therefore, exploiting~\eqref{COMMUTATA0},
\begin{equation}\label{FRTGSDB-00-0234-2a}
\sum_{j,m=1}^n
(\partial_{T,j}\partial_{T,n}\nu_m(x_0))(\partial_{T,m}\partial_{T,n}\nu_j(x_0))=
\sum_{j,m=1}^{n-1}
(\partial_{T,j}\partial_{T,m}\nu_n(x_0))(\partial_{T,m}\partial_{T,j}\nu_n(x_0)).\end{equation}

Besides, we infer from~\eqref{COMMUTATA} that, for all~$j$, $m\in\{1,\dots,n-1\}$,
$$ \partial_{T,j}\partial_{T,m}\phi(x_0)=\partial_{T,m}\partial_{T,j}\phi(x_0),$$
since in this case~$\nu_j(x_0)=\nu_m(x_0)=0$.

On this account, we can write~\eqref{FRTGSDB-00-0234-2a} in the form
\begin{equation}\label{FRTGSDB-00-0234-2}
\sum_{j,m=1}^n
(\partial_{T,j}\partial_{T,n}\nu_m(x_0))(\partial_{T,m}\partial_{T,n}\nu_j(x_0))=
\sum_{j,m=1}^{n-1}
(\partial_{T,m}\partial_{T,j}\nu_n(x_0))^2.\end{equation}

Moreover, using again~\eqref{PNMRmnAQSECEDDR-1},
\begin{eqnarray*}
&& \sum_{i,j,m=1}^n(\partial_{T,m} \partial_{T,j}\nu_i(x_0))^2\\
&=& \sum_{i,j,m=1}^{n-1}(\partial_{T,m} \partial_{T,j}\nu_i(x_0))^2+
\sum_{i,j=1}^n(\partial_{T,n} \partial_{T,j}\nu_i(x_0))^2\\&&\qquad+
\sum_{i,m=1}^n(\partial_{T,m} \partial_{T,n}\nu_i(x_0))^2
+\sum_{j,m=1}^n(\partial_{T,m} \partial_{T,j}\nu_n(x_0))^2\\&=& \sum_{i,j,m=1}^{n-1}(\partial_{T,m} \partial_{T,j}\nu_i(x_0))^2+0+
\sum_{i,m=1}^n(\partial_{T,m} \partial_{T,n}\nu_i(x_0))^2
+\sum_{j,m=1}^n(\partial_{T,m} \partial_{T,j}\nu_n(x_0))^2.
\end{eqnarray*}
Accordingly, by means of~\eqref{COMMUTATA0},
\begin{eqnarray*}
&& \sum_{i,j,m=1}^n(\partial_{T,m} \partial_{T,j}\nu_i(x_0))^2\\
&=& \sum_{i,j,m=1}^{n-1}(\partial_{T,m} \partial_{T,j}\nu_i(x_0))^2+
\sum_{i,m=1}^n(\partial_{T,m} \partial_{T,i}\nu_n(x_0))^2
+\sum_{j,m=1}^n(\partial_{T,m} \partial_{T,j}\nu_n(x_0))^2
\\&=& \sum_{i,j,m=1}^{n-1}(\partial_{T,m} \partial_{T,j}\nu_i(x_0))^2+2
\sum_{j,m=1}^n(\partial_{T,m} \partial_{T,j}\nu_n(x_0))^2.
\end{eqnarray*}
For this reason, recalling~\eqref{PNMRmnAQSECEDDR-1} and~\eqref{PNMRmnAQSECEDDR-2},
\begin{eqnarray*}
\sum_{i,j,m=1}^n(\partial_{T,m} \partial_{T,j}\nu_i(x_0))^2
&=& \sum_{i,j,m=1}^{n-1}(\partial_{T,m} \partial_{T,j}\nu_i(x_0))^2+2
\sum_{j,m=1}^{n-1}(\partial_{T,m} \partial_{T,j}\nu_n(x_0))^2.
\end{eqnarray*}

Combining this, \eqref{FRTGSDB-00-0234-1} and~\eqref{FRTGSDB-00-0234-2}, we arrive at
\begin{equation}\label{REGVD:LL2}\begin{split}&
\frac{\Delta_T c^2(x_0)}2 +c^4(x_0)\\&\qquad=
\sum_{i,j,m=1}^n(\partial_{T,m} \partial_{T,j}\nu_i(x_0))^2
-2\sum_{j,m=1}^n
(\partial_{T,j}\partial_{T,n}\nu_m(x_0))(\partial_{T,m}\partial_{T,n}\nu_j(x_0))\\&\qquad=
\sum_{i,j,m=1}^{n-1}(\partial_{T,m} \partial_{T,j}\nu_i(x_0))^2+2
\sum_{j,m=1}^{n-1}(\partial_{T,m} \partial_{T,j}\nu_n(x_0))^2-2\sum_{j,m=1}^{n-1}
(\partial_{T,m}\partial_{T,j}\nu_n(x_0))^2\\&\qquad=
\sum_{i,j,m=1}^{n-1}(\partial_{T,m} \partial_{T,j}\nu_i(x_0))^2.
\end{split}\end{equation}

Recalling that our goal is to prove~\eqref{0qwoifk0OSJHDN-OKSDND8-UJHS-c}, it is now
useful to give a closer look to the quantity~$|\nabla_Tc|^2$. For this, by~\eqref{REGVD:LL},
\begin{equation}\label{cvhjKJSiddssaPK1MS35}\begin{split}&
c^2|\nabla_Tc|^2=|c\nabla_T c|^2=\left|\frac{\nabla_T c^2}{2}\right|^2=
\frac14\sum_{m=1}^n (\partial_{T,m}c^2)^2\\&\quad=\frac14\sum_{m=1}^n\left( \partial_{T,m}
\sum_{i,j=1}^n(\partial_{T,j}\nu_i)^2
\right)^2=
\sum_{m=1}^n\left(
\sum_{i,j=1}^n(\partial_{T,j}\nu_i)(\partial_{T,m}\partial_{T,j}\nu_i)
\right)^2\\&\quad=
\sum_{i,j,h,k,m=1}^n(\partial_{T,j}\nu_i)(\partial_{T,k}\nu_h)(\partial_{T,m}\partial_{T,j}\nu_i)(\partial_{T,m}\partial_{T,k}\nu_h)\\&\quad=
\sum_{i,j,h,k,m=1}^{n-1}(\partial_{T,j}\nu_i)(\partial_{T,k}\nu_h)(\partial_{T,m}\partial_{T,j}\nu_i)(\partial_{T,m}\partial_{T,k}\nu_h).
\end{split}\end{equation}
We stress that the last step in this calculation
ows to~\eqref{COMMUTATA0}, \eqref{PNMRmnAQSECEDDR-1} and~\eqref{PNMRmnAQSECEDDR-2}.

{F}rom~\eqref{REGVD:LL2} and~\eqref{cvhjKJSiddssaPK1MS35} we arrive at
\begin{equation}\label{LS-sskYmerPLI}
\begin{split}\Upsilon\,:=\,&
\frac{\Delta_T c^2}2+c^4-|\nabla_Tc|^2\\=\,&
\sum_{i,j,m=1}^{n-1}(\partial_{T,m} \partial_{T,j}\nu_i)^2-\frac1{c^2}\sum_{i,j,h,k,m=1}^{n-1}(\partial_{T,j}\nu_i)(\partial_{T,k}\nu_h)(\partial_{T,m}\partial_{T,j}\nu_i)(\partial_{T,m}\partial_{T,k}\nu_h)\\=\,&
\frac1{c^2}
\sum_{i,j,m=1}^{n-1}
\left[
c^2(\partial_{T,m} \partial_{T,j}\nu_i)^2-\sum_{h,k=1}^{n-1}(\partial_{T,j}\nu_i)(\partial_{T,k}\nu_h)(\partial_{T,m}\partial_{T,j}\nu_i)(\partial_{T,m}\partial_{T,k}\nu_h)\right]\\=\,&
\frac1{c^2}
\sum_{i,j,h,k,m=1}^{n-1}
\Big[
(\partial_{T,k}\nu_h)^2(\partial_{T,m} \partial_{T,j}\nu_i)^2-(\partial_{T,j}\nu_i)(\partial_{T,k}\nu_h)(\partial_{T,m}\partial_{T,j}\nu_i)(\partial_{T,m}\partial_{T,k}\nu_h)\Big]
,\end{split}\end{equation}
where all the quantities are computed at the point~$x_0$ as above
(and, from now on, in the calculations needed to establish~\eqref{0qwoifk0OSJHDN-OKSDND8-UJHS-c},
the fact that the quantities involved are computed at~$x_0$ will be omitted, to ease the notation).

To (slightly) simplify the notation in~\eqref{LS-sskYmerPLI}, one can observe that, for tensors~$a_{kh}$ and~$b_{mji}$,
\begin{eqnarray*}&&
\sum_{i,j,h,k,m=1}^{n-1}\Big(a_{kh} b_{mji}-a_{ji} b_{mkh}\Big)^2=\sum_{i,j,h,k,m=1}^{n-1}\Big(
a_{kh}^2 b_{mji}^2+a_{ji}^2 b_{mkh}^2
-2a_{kh} a_{ji}b_{mji} b_{mkh}\Big)\\
&&\qquad=\sum_{i,j,h,k,m=1}^{n-1}\Big(
a_{kh}^2 b_{mji}^2-a_{kh} a_{ji}b_{mji} b_{mkh}\Big)+
\sum_{i,j,h,k,m=1}^{n-1}\Big(
a_{ji}^2 b_{mkh}^2
-a_{kh} a_{ji}b_{mji} b_{mkh}\Big)\\&&\qquad=2
\sum_{i,j,h,k,m=1}^{n-1}\Big(
a_{kh}^2 b_{mji}^2
-a_{ji} a_{kh}b_{mji}b_{mkh} \Big).
\end{eqnarray*}
Using this observation with~$a_{kh}:=\partial_{T,k}\nu_h$ and~$b_{mji}:=
\partial_{T,m} \partial_{T,j}\nu_i$, we write~\eqref{LS-sskYmerPLI} in the form
\begin{equation}\label{P124-001}
\Upsilon=\frac1{2c^2}
\sum_{i,j,h,k,m=1}^{n-1}
\Big((\partial_{T,k}\nu_h)(\partial_{T,m} \partial_{T,j}\nu_i)-(\partial_{T,j}\nu_i)(\partial_{T,m}\partial_{T,k}\nu_h)\Big)^2.
\end{equation}

Now we observe that rotations fixing~$e_n$
preserve the normalization condition~\eqref{nuen}. Accordingly, we may suppose without loss of generality that~$x_0=|x_0|\,e_{n-1}$. Additionally, to complete the proof of~\eqref{0qwoifk0OSJHDN-OKSDND8-UJHS-c} we are now exploiting the cone structure of~$E$ and we obtain that, using the signed distance function~$d$, for all~$t>0$ and~$x\in\partial E$ we have that~$d(tx)=0$ and therefore
$$ 0=\left.\frac{d}{dt}d(tx)\right|_{t=1}=\nabla d(x)\cdot x=\nu(x)\cdot x.$$
This observation and~\eqref{G1014} give that, for each~$x\in\partial E$ and~$i\in\{1,\dots,n-1\}$,
$$ 0=\partial_{T,i}(\nu(x)\cdot x)=\partial_{T,i}\nu(x)\cdot x+\nu(x)\cdot (\partial_{T,i} x)=\partial_{T,i}\nu(x)\cdot x.$$
Evaluating this at the point~$x_0$, we deduce that~$ 0=|x_0|\, \partial_{T,i}\nu(x_0)\cdot e_{n-1}$, whence, by~\eqref{COMMUTATA0},
\begin{equation}\label{nuqw8eighbkns24e3r63e8tyjghn} \partial_{T,n-1}\nu_{i}(x_0)=\partial_{T,i}\nu_{n-1}(x_0)=0.\end{equation}
Consequently, by~\eqref{PNMRmnAQSECEDDR-1} and~\eqref{P124-001},
\begin{equation}\label{234ty5-23tyh-24T460m}
\begin{split}
\Upsilon\,&=\frac1{2c^2}\Bigg[
\sum_{{1\le i,j,h,m\le n-1}\atop{1\le k\le n-2}}
\Big((\partial_{T,k}\nu_h)(\partial_{T,m} \partial_{T,j}\nu_i)-(\partial_{T,j}\nu_i)(\partial_{T,m}\partial_{T,k}\nu_h)\Big)^2
\\&\qquad\qquad+
\sum_{{1\le i,j,h,m\le n-1}}(\partial_{T,j}\nu_i)^2(\partial_{T,m}\partial_{T,n-1}\nu_h)^2\Bigg]\\
&=\frac1{2c^2}\Bigg[
\sum_{{1\le i,h,m\le n-1}\atop{1\le j,k\le n-2}}
\Big((\partial_{T,k}\nu_h)(\partial_{T,m} \partial_{T,j}\nu_i)-(\partial_{T,j}\nu_i)(\partial_{T,m}\partial_{T,k}\nu_h)\Big)^2
\\&\qquad\qquad+
\sum_{{1\le i,h,m\le n-1}\atop{1\le k\le n-2}}
(\partial_{T,k}\nu_h)^2(\partial_{T,m} \partial_{T,n-1}\nu_i)^2
+c^2\sum_{{1\le h,m\le n-1}}(\partial_{T,m}\partial_{T,n-1}\nu_h)^2\Bigg]\\
&=\frac1{2c^2}\Bigg[
\sum_{{1\le i,m\le n-1}\atop{1\le j,h,k\le n-2}}
\Big((\partial_{T,k}\nu_h)(\partial_{T,m} \partial_{T,j}\nu_i)-(\partial_{T,j}\nu_i)(\partial_{T,m}\partial_{T,k}\nu_h)\Big)^2
\\&\qquad\qquad+
\sum_{{1\le i,m\le n-1}\atop{1\le j,k\le n-2}}(\partial_{T,j}\nu_i)^2(\partial_{T,m}\partial_{T,k}\nu_{n-1})^2
+2c^2\sum_{{1\le h,m\le n-1}}(\partial_{T,m}\partial_{T,n-1}\nu_h)^2\Bigg]\\
&=\frac1{2c^2}\Bigg[
\sum_{{1\le m\le n-1}\atop{1\le i,j,h,k\le n-2}}
\Big((\partial_{T,k}\nu_h)(\partial_{T,m} \partial_{T,j}\nu_i)-(\partial_{T,j}\nu_i)(\partial_{T,m}\partial_{T,k}\nu_h)\Big)^2
\\&\qquad\qquad+
\sum_{{1\le m\le n-1}\atop{1\le j,h,k\le n-2}}
(\partial_{T,k}\nu_h)^2(\partial_{T,m} \partial_{T,j}\nu_{n-1})^2
+3c^2\sum_{{1\le h,m\le n-1}}(\partial_{T,m}\partial_{T,n-1}\nu_h)^2\Bigg]\\&=\frac1{2c^2}\Bigg[
\sum_{{1\le m\le n-1}\atop{1\le i,j,h,k\le n-2}}
\Big((\partial_{T,k}\nu_h)(\partial_{T,m} \partial_{T,j}\nu_i)-(\partial_{T,j}\nu_i)(\partial_{T,m}\partial_{T,k}\nu_h)\Big)^2
\\&\qquad\qquad+4c^2\sum_{{1\le h,m\le n-1}}(\partial_{T,m}\partial_{T,n-1}\nu_h)^2\Bigg]\\&\ge2\sum_{{1\le h,m\le n-1}}(\partial_{T,m}\partial_{T,n-1}\nu_h)^2.
\end{split}\end{equation}

Moreover, by~\eqref{COMMUTATA}, \eqref{nuen}
and~\eqref{nuqw8eighbkns24e3r63e8tyjghn}, at the point~$x_0$ we have that
\begin{equation}\label{NSD-SNCJK-0121} \partial_{T,i}\partial_{T,n-1}\phi-
\partial_{T,n-1}\partial_{T,i}\phi=\sum_{k=1}^n
\Big[\nu_i\partial_{T,n-1}\nu_k-\nu_{n-1}\partial_{T,i}\nu_k\Big]\partial_{T,k}\phi=0.\end{equation}
Similarly, at the point~$x_0$ we have that
\begin{equation}\label{NSD-SNCJK-0122} \partial_{T,n-1}\phi=
\partial_{n-1}\phi-(\nabla\phi\cdot\nu)\nu_{n-1}=\partial_{n-1}\phi.\end{equation}

Furthermore, since~$E$ is a cone we have that~$\nu(t x_0)=\nu(x_0)$ for all~$t>0$
and therefore~$\partial_{T,m}\nu_h$ is a homogeneous function of degree~$-1$.
This yields that
$$ \partial_{T,m}\nu_h(x_0)=-\nabla \partial_{T,m}\nu_h(x_0)\cdot x_0=
-|x_0|\,\partial_{n-1}\partial_{T,m}\nu_h(x_0).$$
Hence, by~\eqref{NSD-SNCJK-0121} and~\eqref{NSD-SNCJK-0122},
$$ \partial_{T,m}\partial_{T,n-1}\nu_h(x_0)=\partial_{T,n-1}\partial_{T,m}\nu_h(x_0)=
\partial_{n-1}\partial_{T,m}\nu_h(x_0)=-\frac{\partial_{T,m}\nu_h(x_0)}{|x_0|}.
$$
Combining this information with~\eqref{234ty5-23tyh-24T460m} we conclude that
$$ \Upsilon
\ge2\sum_{{1\le h,m\le n-1}}\frac{\big(\partial_{T,m}\nu_h(x_0)\big)^2}{|x_0|^2}=\frac{2c^2}{|x_0|^2}.$$
The proof of~\eqref{0qwoifk0OSJHDN-OKSDND8-UJHS-c} is thereby complete.\medskip
\end{proof}

\subsection{Bernstein's problem and a conjecture by E. De Giorgi}
We also recall that the regularity of minimal surfaces is strictly linked to the so-called
Bernstein's problem (named after S. Bernstein who first solved the case in dimension~3,
see~\cite{MR1544873}). 
Bernstein's problem asks whether or not a minimal graph in~$\R^n$ (i.e., a minimal
surface which possesses a global graphical structure of the form~$x_n=u(x')$ with~$x'\in\R^{n-1}$) is necessarily affine.
The answer to this problem is affirmative in \label{0ojlrf4-234rot}
dimension~$n\le8$
(due to the works of~\cite{MR157263, MR178385, MR1556840, MR200816, MR233295}), and negative when~$ n\ge 9$ (see~\cite{MR250205}).

As a matter of fact, the connection between Bernstein's problem and the regularity
of minimal surfaces in~\eqref{0qwoifk0OSJHDN-OKSDND8-UJHS} is clearly showcased in~\cite{MR178385, MR1556840} by showing that
\begin{equation}\label{DGIoBErdfv}
\begin{split}&{\mbox{if all minimal cones in~$\R^{n-1}$ are halfplanes,}}\\ &{\mbox{then Bernstein's problem has an affirmative answer in~$\R^n$.}}\end{split}\end{equation}\medskip

The link between Bernstein's problem and the limit interfaces of phase transition models
(as described by the $\Gamma$-convergence theory in Theorem~\ref{MORTOLA})
was possibly an inspiring motivation for Ennio De Giorgi to state one of his most famous conjectures~\cite{MR533166}.

The gist of this conjecture could be as follows: given that, at a large scale,
the level sets of ``good'' solutions of the Allen-Cahn equation approach perimeter minimizing surfaces (as made precise by the $\Gamma$-convergence theory in Theorem~\ref{MORTOLA}
and by the geometric convergence of Theorem~\ref{THM:CC})
and given that minimal graphs reduce to hyperplanes in dimension~$n\le8$ (according to 
Bernstein's problem), would it be possible that level sets of ``good'' global solutions
of the Allen-Cahn equation are already hyperplanes? Since level sets corresponding to different values of the solution cannot intersect, this would say that all the level sets are in fact parallel hyperplanes and therefore the solution only depends on the distance to one of these hyperplanes
(in particular, the solution would be a function depending only on one Euclidean variable).

In all this heuristic discussion, we have been vague about what a ``good'' solution precisely is:
in a sense, besides boundedness and regularity assumptions, in view of Theorems~\ref{MORTOLA} and~\ref{THM:CC} a natural hypothesis would be to require that the solution is a local minimizer; furthermore, to fall within the range of application of
Bernstein's problem, it would be desirable to know that the limit minimal surface has a graphical structure and for this some monotonicity assumption on the solution could be helpful (since, at least locally, it would entail a graphical structure of the level set via Implicit Function Theorem).

It would be however desirable to keep the number of assumptions to the minimum and possibly to confine them to assumptions of ``geometric'' type: in this spirit, one may be tempted to remove the minimality assumption (which is instead of ``variational'' and ``energetic'' type)
and focus mainly on a monotonicity assumption
(roughly speaking, after all, maybe monotonicity is already an indication of some ``weak'' form of minimality\footnote{A more precise link between monotonicity
and this weak notion of minimality, and more precisely stability, will be given in~\eqref{STABBI}.} since it avoids oscillations that increase energy).

These arguments (and likely many others of much deeper nature) have probably inspired De Giorgi for the precise formulation of this conjecture, which goes as follows:

\begin{conjecture} \label{DG:CO}
Let~$n\le8$ and~$u\in C^2(\R^n)\cap L^\infty(\R^n)$ be a global solution of the Allen-Cahn equation
$$ -\Delta u=u-u^3$$ such that
$$ \partial_n u(x)>0\quad{\mbox{ for every }}\,x\in\R^n.$$
Is it true that~$u$ is one-dimensional, i.e. that there exist~$u_0:\R\to\R$ and~$\omega\in\partial B_1$ such that~$u(x)=u_0(\omega\cdot x)$ for all~$x\in\R^n$?
\end{conjecture}

Conjecture~\ref{DG:CO} has been proven for~$n\in\{2,3\}$
but it is still open for~$n\in\{4,\dots,8\}$, see~\cite{MR1637919, MR1655510, MR1775735, MR1843784}. For~$n\ge9$,
an example of global, bounded and monotone solution of the Allen-Cahn equation which is not one-dimensional has been constructed in~\cite{MR2846486}, thus confirming that the dimensional restriction in Conjecture~\ref{DG:CO} cannot be removed.

In dimension~$n\in\{4,\dots,8\}$ Conjecture~\ref{DG:CO} is known to hold under an additional assumption on the ``profiles of the solution at infinity''. Namely, since~$u$
in Conjecture~\ref{DG:CO} is bounded and monotone in the direction of~$e_n$, one can define, for all~$x'\in\R^{n-1}$,
\begin{equation*} \overline{u}(x'):=\lim_{x_n\to+\infty}u(x',x_n)\qquad{\mbox{and}}\qquad
\underline{u}(x'):=\lim_{x_n\to-\infty}u(x',x_n).\end{equation*}
In this setting, it has been proved in~\cite{MR2480601} that Conjecture~\ref{DG:CO}
holds true under the additional assumption
\begin{equation}\label{omUjPoOP1-23e}
\overline{u}(x')=-\underline{u}(x')=1 \quad{\mbox{ for every }}\,x'\in\R^{n-1}.
\end{equation}
For further results establishing Conjecture~\ref{DG:CO} under suitable assumptions
on the limit profiles see~\cite{MR2483642, MR2728579, MR3488250}. See also~\cite{MR2528756}
for an overview of Conjecture~\ref{DG:CO} and of related problems.\medskip

We also mention that a related problem, posed by theoretical physicist Gary William Gibbons 
consisted in determining whether a bounded global solution of the Allen-Cahn equation
is necessarily one dimensional under the uniform limit assumption
\begin{equation*} \lim_{T\to+\infty}\sup_{x_n\ge T} |u(x',x_n)-1|=0\qquad{\mbox{and}}\qquad
\lim_{T\to+\infty}\sup_{x_n\le -T} |u(x',x_n)+1|=0.\end{equation*}
Note that this condition is stronger than~\eqref{omUjPoOP1-23e}.
The answer to Gibbons' problem is known to be positive for every dimension~$n$, see~\cite{MR1755949, MR1763653, MR1765681}.\medskip

It is now worth coming back to the relation between monotonicity and some weak form of minimality which was raised before the statement of Conjecture~\ref{DG:CO}.
A precise notion of this is given by the observation that monotonicity implies stability:
namely, if~$u$ is a solution of
$$ \Delta u=W'(u)$$
such that~$\partial_n u>0$ in some domain~$\Omega\subseteq\R^n$, then, for all~$\phi\in C^\infty_0(\Omega)$, we have that
\begin{equation}\label{STABBI}
\int_\Omega \Big( |\nabla\phi(x)|^2+W''(u(x))\,\phi^2(x)\Big)\,dx\ge0.
\end{equation}
Indeed, under the monotonicity assumption it is fair to define~$\psi:=\frac{\phi^2}{\partial_n u}$
and infer that
\begin{eqnarray*}
&&\int_\Omega \Big( |\nabla\phi(x)|^2+W''(u(x))\,\phi^2(x)\Big)\,dx
=\int_\Omega \left( |\nabla\phi(x)|^2+\partial_n\Big( W'(u(x))\Big)\,\frac{\phi^2(x)}{\partial_nu(x)}\right)\,dx\\&&\qquad=
\int_\Omega \left( |\nabla\phi(x)|^2+\partial_n\Big( \Delta u(x)\Big)\,\psi(x)\right)\,dx\\&&\qquad=
\int_\Omega \Big( \big|\nabla\big(\sqrt{\psi(x)}\,\sqrt{\partial_nu(x)}\big)\big|^2-\nabla\partial_n u(x)\cdot\nabla\psi(x)\Big)\,dx\\&&\qquad=
\int_\Omega \left( \left|
\frac{\sqrt{\partial_nu(x)}\,\nabla\psi(x)}{2\sqrt{\psi(x)}} 
+\frac{\sqrt{\psi(x)}\,\nabla\partial_nu(x)}{2\sqrt{\partial_nu(x)}}\right|^2-\nabla\partial_n u(x)\cdot\nabla\psi(x)\right)\,dx\\&&\qquad=
\int_\Omega \left( 
\frac{{\partial_nu(x)}\,|\nabla\psi(x)|^2}{4 {\psi(x)}} 
+\frac{{\psi(x)}\,|\nabla\partial_nu(x)|^2}{4{\partial_nu(x)}}-\frac12\nabla\partial_n u(x)\cdot\nabla\psi(x)\right)\,dx\\&&\qquad=\int_\Omega  \left|
\frac{\sqrt{\partial_nu(x)}\,\nabla\psi(x)}{2\sqrt{\psi(x)}} 
-\frac{\sqrt{\psi(x)}\,\nabla\partial_nu(x)}{2\sqrt{\partial_nu(x)}}\right|^2 \,dx
\ge0,
\end{eqnarray*}
which is~\eqref{STABBI}.

We also note that~\eqref{STABBI} states that the second derivative of the corresponding
energy functional (say, \eqref{FEPSI} with~$\e:=1$) is nonnegative. In particular,
minimizers satisfy the stability inequality in~\eqref{STABBI}.

\section{Long-range interactions and the nonlocal Allen-Cahn equation}\label{x-09i8uytf-0iuytfdfghjdP}

In view of the discussion on page~\pageref{GURNL-sal}, it is of interest to consider minimizers,
and more generally critical points,
of the long-range energy functional
\begin{equation*} \frac14 \iint_{{\mathcal{Q}}(\Omega)} \frac{(u(x)-u(y))^2}{|x-y|^{n+\alpha}}\,dx\,dy+\int_\Omega W(u(x))\,dx,\end{equation*}
where~${\mathcal{Q}}(\Omega)$
has been defined in~\eqref{GURNL-2}, for a given~$\alpha\in(0,2)$,
and this not only in view of natural generalizations of the classical setting to more complicated ones,
but also due to the structural formulation of the phase coexistence problem which is intrinsically
long-range (being the short-range case a very handy and important simplification).
 
The corresponding critical points in this setting are solutions of
the fractional counterpart of the Allen-Cahn equation given by
\begin{equation*}
(-\Delta)^{\frac\alpha2}u(x)+W'(u(x))=0
\end{equation*}
for~$x\in\Omega$ (to be compared with the classical case in~\eqref{ACHAW}).

The case of heterogeneous materials (to be compared with~\eqref{ACHAW2})
can also be taken into account, via the more general equation
\begin{equation*}
(-\Delta)^{\frac\alpha2} u(x)+Q(x)\,W'(u(x))=0,\end{equation*}
for~$Q$
ranging in~$[\underline{Q},\overline{Q}]$, for some~$\overline{Q}\ge \underline{Q}>0$.

Here above, we are using the fractional Laplacian operator, defined (up to a normalizing constant that we omit) as
$$ (-\Delta)^{\frac\alpha2} u(x)=
\int_{\R^n} \frac{u(x)-u(y)}{|x-y|^{n+\alpha}}\,dx\,dy.$$
Notice that the above integral is singular, hence, if needed, one has to consider it in the Cauchy principal value sense,
to allow for cancellations (see e.g.~\cite{MR2707618, MR3967804}
and the references therein for the basics on the fractional Laplacian).

Solutions of the fractional heterogeneous Allen-Cahn equation in a periodic medium
whose interface is trapped within two hyperplanes have been constructed in~\cite{MR3816747}
(this can be considered as a nonlocal counterpart of Theorem~\ref{MS-PSLDibvAlcK2}).
See also~\cite{MR4410572} for a general discussion about the
fractional Allen-Cahn equation.

\section{The nonlocal limit interface and the theory of nonlocal~$\Gamma$-convergence}\label{GAMMACO2}

The long-range interaction energies present a richer $\Gamma$-convergence theory than their
classical counterpart. First of all, the short-range functional in~\eqref{WDW-CLASSi} needs to be replaced by its fractional analog
in~\eqref{PRES}, but also the fractional counterpart of the
rescaled functional in~\eqref{FEPSI} requires some care in determining the appropriate penalization parameters.
Moreover, the result obtained in this case deeply depends on the fractional exponent~$\alpha$:
as we will now clarify, for~$\alpha\in(0,1)$ the $\Gamma$-limit is related to nonlocal minimal surfaces and
for~$\alpha\in[1,2)$ to classical ones. This is especially interesting since it suggests that, while
for small values of the fractional exponent the phase transition problem always maintains a clearly distinctive long-range feature,
for large values of the fractional exponent, at a large scale, the phase transition problem
tends to resemble a local one (the threshold between ``small'' and ``large'' fractional parameters being given by
the specific value~$\alpha=1$).

As for the ``appropriate'' rescaled functional, one can try to get inspired by~\eqref{EPSACAN}
and look for a perturbative fractional Allen-Cahn equation of the form
\begin{equation*}
\e^\alpha (-\Delta)^\alpha u(x)+W'(u(x))=0,\end{equation*}
for a small parameter~$\e$. This would correspond to an energy functional (up to normalization constants that we disregard) of the type
\begin{equation} \label{PROKFba}\e^{\alpha}\iint_{{\mathcal{Q}}(\Omega)} \frac{(u(x)-u(y))^2}{|x-y|^{n+\alpha}}\,dx\,dy+\int_\Omega W(u(x))\,dx.\end{equation}
But to develop a solid theory of $\Gamma$-convergence one would like to have an energy functional that has the tendency to remain
bounded as the perturbative parameter vanishes (and note that this would have been a point to raise
even in the classical case, when shifting from the equation in~\eqref{EPSACAN}
to the energy functional in~\eqref{FEPSI}): for this, one exploits the freedom of multiplying the energy by a scalar,
which does not change the critical points, and then chooses appropriately this gauge to have a control of the energy
as the perturbative parameter vanishes.
That is, without affecting the minimization problem, we replace~\eqref{PROKFba} with the functional
\begin{equation}\label{PROKFba2}
\varsigma_\e
\left[\e^\alpha\iint_{{\mathcal{Q}}(\Omega)} \frac{(u(x)-u(y))^2}{|x-y|^{n+\alpha}}\,dx\,dy+\int_\Omega W(u(x))\,dx\right]\end{equation}
and we choose~$\varsigma_\e>0$ such that the energy of a ``typical'' transition
remains bounded, and nontrivial, as~$\e\searrow0$.

As a model transition for this calculation, one can suppose that~$\Omega=(-1,1)^n$ and take, for instance, $u_\e(x):=u_\star\left(\frac{x_1}\e\right)$ for a given smooth, odd function~$u_\star:\R\to\R$ with~$u_\star(0)=0$, increasing in~$(-1,1)$ and with~$u_\star=1$ in~$[1,+\infty)$.

We note that, on the one hand,
\begin{equation*}
\int_\Omega W(u_\e(x))\,dx=\int_{(-\e,\e)\times(-1,1)^{n-1}} W\left(u_\star\left(\frac{x_1}\e\right)\right)\,dx=
\e^n\int_{(-1,1)\times\left(-\frac1\e,\frac1\e\right)^{n-1}} W (u_\star(y_1) )\,dy
\simeq \e.
\end{equation*}
On the other hand,
\begin{equation}\label{jyhtgrfu7y6tg5rfet8569b7c5tnxm4utxc46}
\e^{\alpha-n}
\iint_{{\mathcal{Q}}(\Omega)} \frac{(u_\e(x)-u_\e(y))^2}{|x-y|^{n+\alpha}}\,dx\,dy \simeq
\begin{dcases}
\e^{\alpha-n}&{\mbox{ if }}\alpha\in(0,1),\\
\e^{1-n}|\ln\e|&{\mbox{ if }}\alpha=1,\\
\e^{1-n}&{\mbox{ if }}\alpha\in(1,2).
\end{dcases}\end{equation}
We defer the proof of~\eqref{jyhtgrfu7y6tg5rfet8569b7c5tnxm4utxc46}
to Appendix~\ref{sec:jyhtgrfu7y6tg5rfet8569b7c5tnxm4utxc46} not to interrupt the flow of the argument.

As a result,
\begin{eqnarray*}
&&\e^\alpha\iint_{{\mathcal{Q}}(\Omega)} \frac{(u(x)-u(y))^2}{|x-y|^{n+\alpha}}\,dx\,dy+\int_\Omega W(u(x))\,dx=\begin{dcases}
\e^\alpha&{\mbox{ if }}\alpha\in(0,1),\\
\e|\ln\e|&{\mbox{ if }}\alpha=1,\\
\e&{\mbox{ if }}\alpha\in(1,2).
\end{dcases}
\end{eqnarray*}
Thus, in light of~\eqref{PROKFba2}, it is convenient to consider, as perturbed long-range energy functional
(to be compared with~\eqref{FEPSI}) the quantity
\begin{equation}\label{PROKFba3}\begin{split}&{\mathcal{F}}_\e^\alpha(u):=
\varsigma_\e
\left[\e^\alpha\iint_{{\mathcal{Q}}(\Omega)} \frac{(u(x)-u(y))^2}{|x-y|^{n+\alpha}}\,dx\,dy+\int_\Omega W(u(x))\,dx\right],\\
&\qquad\qquad {\mbox{with }}\;\varsigma_\e:=\begin{dcases}
\e^{-\alpha} & {\mbox{ if }}\alpha\in(0,1),\\
(\e|\ln\e|)^{-1} & {\mbox{ if }}\alpha=1,\\
\e^{-1} &{\mbox{ if }}\alpha\in(1,2).
\end{dcases}\end{split}\end{equation}
The $\Gamma$-convergence theory for this object has been established in
Theorems~1.4 and~1.5 in~\cite{MR2948285}, thus providing the long-range counterpart of Theorem~\ref{MORTOLA}:

\begin{theorem}\label{oqjdlwn0304-53-tgol}
The functional~${\mathcal{F}}_\e^\alpha$ in~\eqref{PROKFba3}
$\Gamma$-converges as~$\e\searrow 0$ to
\begin{equation*} {\mathcal{F}}^\alpha(u):=\begin{dcases}
c\,\Per_\alpha(E,\Omega) & \begin{matrix*}[l]&{\mbox{ if $\alpha\in(0,1)$ and $u=\chi_E-\chi_{\R^n\setminus E}$}} \\ &{\mbox{ for some set~$E$ of finite $\alpha$-perimeter,}}\end{matrix*}\\
\\
c\,\Per(E,\Omega) & \begin{matrix*}[l]&{\mbox{ if $\alpha\in[1,2)$ and $u=\chi_E-\chi_{\R^n\setminus E}$}} \\ &{\mbox{ for some set~$E$ of finite perimeter,}}\end{matrix*}\\
\\
+\infty & \begin{matrix*}[l]&{\mbox{ otherwise,}}\end{matrix*}
\end{dcases}\end{equation*}
where~$c>0$ depends on~$n$, $\alpha$ and~$W$.
\end{theorem}

The proof of Theorem~\ref{oqjdlwn0304-53-tgol} is relatively straightforward when~$\alpha\in(0,1)$
because in this case step functions are admissible not only for the limit functional~${\mathcal{F}}^\alpha$
but also for the original functional~${\mathcal{F}}_\e^\alpha$ (e.g., this provides a recovery sequence straight away).
But the case~$\alpha\in[1,2)$ is trickier, since one needs to relate a long-range functional
to the classical perimeter in the limit, and for this a careful analysis of different integral contributions
is mandatory and some cancellations must be suitably spotted. 

For example,
to understand which integral contributions survive in the classical setting obtained
in the limit, it is useful to have a lower bound on the nonlocal interaction in~\eqref{per-I}.
For this, if~${\mathcal{A}}$ and~${\mathcal{D}}$ are, say, disjoint subsets of the unit cube~${\mathcal{Q}}$,
when~$\alpha\in(0,1)$ we know from~\eqref{per-IO}
that~$ {\mathcal{I}}_\alpha({\mathcal{A}},{\mathcal{D}})\le
{\rm Per}_\alpha({\mathcal{A}})$ which is finite for smooth and bounded sets~${\mathcal{A}}$.
Instead, when~$\alpha\in[1,2)$ contributions of the type~${\mathcal{I}}_\alpha({\mathcal{A}},{\mathcal{Q}}\setminus{\mathcal{A}})$ are infinite (unless the sets involved in the interactions become of null measure).
One can indeed quantify this feature (see Proposition~3.1 in~\cite{MR2948285}) and find that,
if~$\min\{|{\mathcal{A}}|,|{\mathcal{D}}|\}\ge c>0$, then
\begin{equation}\label{jdnDAVSsrytt902P5} {\mathcal{I}}_\alpha({\mathcal{A}},{\mathcal{D}})\ge\begin{dcases}
c_0\,\big|\ln\big|{\mathcal{Q}}\setminus({\mathcal{A}}\cup{\mathcal{D}})\big|\big| & {\mbox{ if }}\alpha=1,\\
c_0\,\big|{\mathcal{Q}}\setminus({\mathcal{A}}\cup{\mathcal{D}})\big|^{1-\alpha}& {\mbox{ if }}\alpha\in(1,2),
\end{dcases}\end{equation}
for some~$c_0 > 0$ depending only on~$c$, $n$ and~$\alpha$. Estimates of this sort are useful
since they entail that order parameters exhibiting a positive measure of states close to both the pure phases
necessarily have energy bounded from below (as we can expect, some energy has to be spent to
produce two different phases and these kinds of estimates provide bounds on the energy of the limit
interface).

Additional complications for the proof of Theorem~\ref{oqjdlwn0304-53-tgol}
surface in the construction of recovery sequences, since some fine energy estimate is needed
to control the interpolation of two functions across a given domain, and also the existence and basic
properties of one-dimensional transition layers play a significant role (see~\cite{MR3081641, MR3165278}).

For further results about $\Gamma$-convergence theories related to nonlocal problems,
see~\cite{MR1612250, MR2546026}.\medskip

Having settled the $\Gamma$-convergence theory in the nonlocal framework, we now point
out that the geometric convergence in Theorem~\ref{THM:CC} has also a nonlocal counterpart,
as established in Theorems~1.2, 1.3 and~1.4, and Corollary~1.7,
of~\cite{MR3133422} (see also Theorem~1.6 there for a convenient extension of~\eqref{jdnDAVSsrytt902P5}):

\begin{theorem}\label{THM:CC:E}
Assume that~$u_\e$ is a local minimizer for the functional~${\mathcal{F}}_\e^\alpha$
in~\eqref{PROKFba3} in the ball~$B_{1+\e}$.

Then:
\begin{itemize}
\item There exists~$C>0$, depending only on~$n$, $\alpha$ and~$W$, such that
\begin{equation*} {\mathcal{F}}_\e^\alpha(u_\e,B_1)\le C.\end{equation*}
\item Up to a subsequence,
\begin{equation*}
{\mbox{$u_\e\to\chi_E-\chi_{\R^n\setminus E}$ as~$\e\searrow0$ in~$L^1(B_1)$}}\end{equation*}
and the set~$E$ has locally minimal perimeter in~$B_1$ when~$\alpha\in[1,2)$
and locally minimal $\alpha$-perimeter when~$\alpha\in(0,1)$.
\item Given~$\vartheta_1$, $\vartheta_2\in(-1,1)$, if~$u_\e(0)>\vartheta_1$, then
\begin{equation*} \big|\{u_\e>\vartheta_2\}\cap B_r\big|\ge cr^n,\end{equation*}
as long as~$r\in(0,1]$ and~$\e\in(0, c_\star\,r]$, where~$c>0$ depends only on~$n$, $\alpha$ and~$W$
and~$c_\star>0$ depends only on~$n$, $\alpha$, $W$, $\vartheta_1$ and~$\vartheta_2$.
\item Similarly, given~$\vartheta_1$, $\vartheta_2\in(-1,1)$, if~$u_\e(0)<\vartheta_1$, then
\begin{equation*} \big|\{u_\e<\vartheta_2\}\cap B_r\big|\ge cr^n,\end{equation*}
as long as~$r\in(0,1]$ and~$\e\in(0, c_\star\,r]$, where~$c>0$ depends only on~$n$, $\alpha$ and~$W$
and~$c_\star>0$ depends only on~$n$, $\alpha$, $W$, $\vartheta_1$ and~$\vartheta_2$.
\item For every~$\vartheta\in(0,1)$, the set~$\{|u_\e|<\vartheta\}$ approaches~$\partial E$
locally uniformly as~$\e\searrow0$: more explicitly, given~$r_0\in(0,1)$ and~$\delta>0$
there exists~$\e_0>0$ such that if~$\e\in(0,\e_0)$ then
\begin{equation*} \{|u_\e|<\vartheta\}\cap B_{r_0}\subseteq \bigcup_{p\in\partial E} B_{\delta}(p).\end{equation*}
\end{itemize}
\end{theorem}

We stress that Theorem~\ref{THM:CC:E} is a perfect nonlocal counterpart of Theorem~\ref{THM:CC}:
indeed, being of measure theoretic type, the exponents involved do not depend on~$\alpha$
(but the structural constants may).

\section{Nonlocal minimal surfaces and one-dimensional symmetry}\label{09i2w-e23rt5-5PKM78}

\subsection{Nonlocal minimal surfaces}
Given the strong connection between long-range phase transitions in the ``genuinely nonlocal'' range~$\alpha\in(0,1)$ and the nonlocal minimal surfaces (as showcased in Theorems~\ref{oqjdlwn0304-53-tgol}
and~\ref{THM:CC:E}), it is desirable to understand better the regularity and flatness properties of
the minimizers of the $\alpha$-perimeter. While the classical situation is fully understood,
in light of~\eqref{0qwoifk0OSJHDN-OKSDND8-UJHS}
and~\eqref{0qwoifk0OSJHDN-OKSDND8-UJHS-x2}, the nonlocal counterpart of this regularity theory is broadly open. To the best of our knowledge, no example of singular nonlocal minimal surface is known,
and an interior regularity theory for $\alpha$-perimeter minimizers
in~$\R^n$ has been established
\begin{eqnarray}
\label{AGbhgVjwehe-1}&&{\mbox{when~$n=2$, for all~$\alpha\in(0,1)$,}}\\
\label{AGbhgVjwehe-3}&&{\mbox{when~$n\le7$, for all~$\alpha\in(\alpha_0,1)$, for a suitable~$\alpha_0\in(0,1)$,}}\\
\label{AGbhgVjwehe-5}&&{\mbox{when the set has a graphical structure, for all~$\alpha\in(0,1)$ and~$n\in\N$.}}
\end{eqnarray}
Indeed, \eqref{AGbhgVjwehe-1} has been established in Corollary~1 of~\cite{MR3090533}
(with smooth regularity coming from Theorem~1.1 in~\cite{MR3331523}),
\eqref{AGbhgVjwehe-3} in Theorem~3 of~\cite{MR3107529}, and~\eqref{AGbhgVjwehe-5} in Theorem~1.1 of~\cite{MR3934589}.

Notice in particular that~\eqref{AGbhgVjwehe-3} carries the classical minimal surfaces regularity theory in~\eqref{0qwoifk0OSJHDN-OKSDND8-UJHS}
over to the nonlocal setting, provided that the fractional exponent is ``large enough''.\medskip

As for the nonlocal version of Bernstein's problem, the classical result by De Giorgi in~\eqref{DGIoBErdfv}
has a full counterpart for nonlocal minimal graphs, as proved in
Theorem~1.2 of~\cite{MR3680376} (see also~\cite{MR4050198} for a different proof and
for related results): therefore one can infer from~\eqref{AGbhgVjwehe-1}
and~\eqref{AGbhgVjwehe-3}
that
\begin{equation}\label{BEJDrn92o3esefrP}
\begin{split}&
{\mbox{Bernstein's problem for $\alpha$-minimal graphs has}}\\&{\mbox{an affirmative answer in~$\R^n$
when~$n\le3$ for all~$\alpha\in(0,1)$,}}\\&{\mbox{as well as
when~$n\le8$ for all~$\alpha\in(\alpha_0,1)$, for a suitable~$\alpha_0\in(0,1)$}}.\end{split}\end{equation}
Once again, this transfers into the nonlocal world the classical affirmative answer to Bernstein's problem
discussed on page~\pageref{0ojlrf4-234rot}, provided that the fractional exponent is ``large enough''
(no nonlocal counterexamples to the Bernstein's problem being available so far).
\medskip

\begin{figure}[h]
\includegraphics[width=0.65\textwidth]{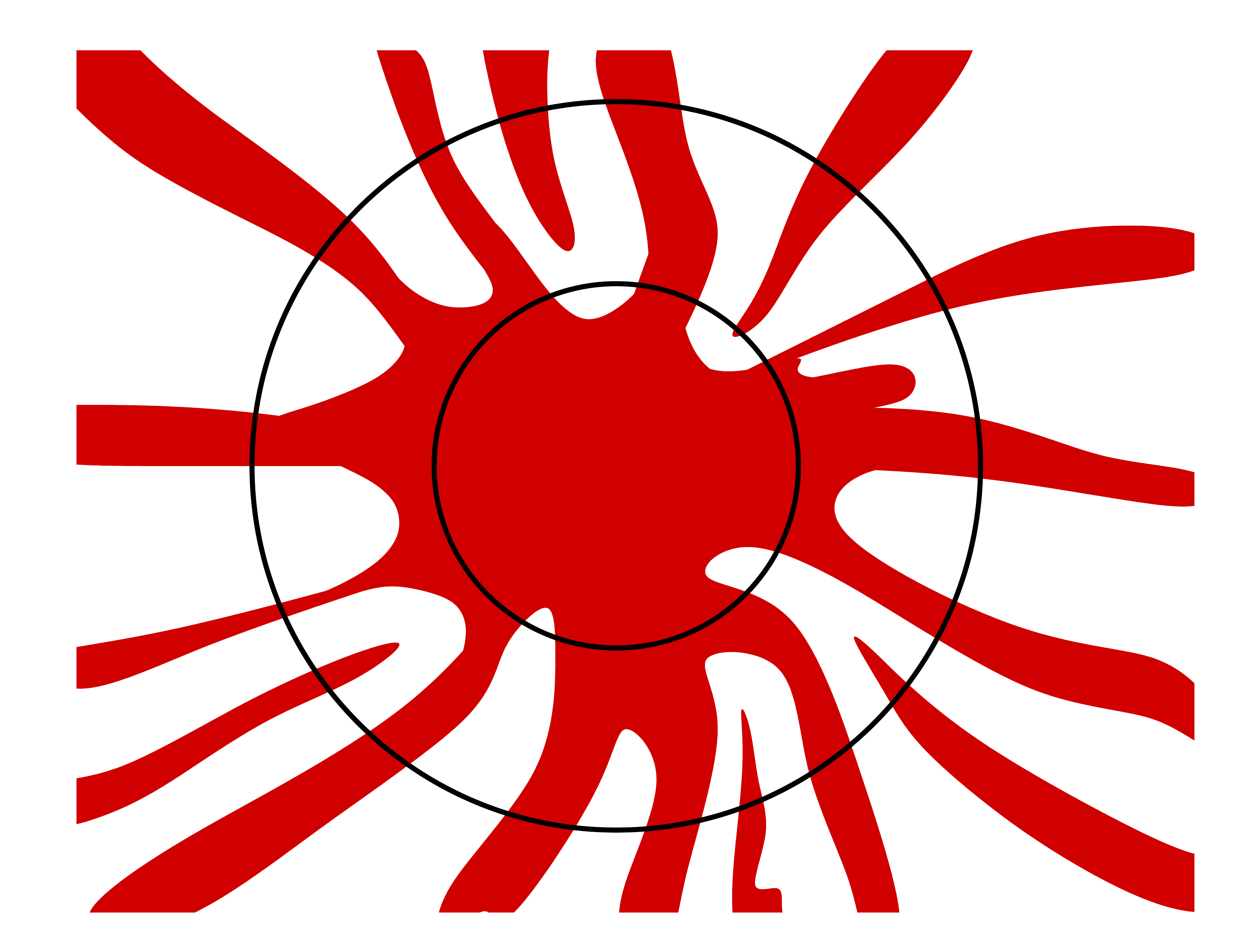}
\caption{Fingers either merging or getting extinguished.}
        \label{cheb3erAV}
\end{figure}

It is open, and very interesting, to understand whether the results in~\eqref{AGbhgVjwehe-1}
and~\eqref{AGbhgVjwehe-3} can be carried over to higher dimensions.
In such a generality, an interesting feature is that, even without smooth regularity results,
one can obtain a perimeter estimate for $\alpha$-minimal surfaces, and, even more generally,
for $\alpha$-stable ones. As a matter of fact, as proved in Theorem~1.1
of~\cite{MR3981295}, for any dimension~$n$ and any fractional exponent~$\alpha\in(0,1)$,
given any~$\alpha$-stable set in the unit ball~$B_1$, its classical perimeter in~$B_{1/2}$
is bounded by a constant depending only on~$n$ and~$\alpha$.
This is interesting not only because it provides a regularity result
(and note that even for $\alpha$-minimizers, the perimeter bound is sharper than an~$\alpha$-perimeter
bound which can follow by comparison with a given competitor), but also because
the interior bound does not depend on the shape of the set outside~$B_1$.
In particular, the given $\alpha$-stable set in the unit ball~$B_1$ may have arbitrarily many ``fingers'' reaching out towards infinity, but their contribution to the perimeter in~$B_{1/2}$
is bounded universally, independently on the number of these fingers: that is, most of the fingers have
the tendency to either merge or get extinguished before getting to the center, to avoid the formation
of extra perimeter, see Figure~\ref{cheb3erAV}. Notice that this is a purely nonlocal phenomenon,
since in the classical case, an arbitrarily large number of parallel hyperplanes would describe
the case of a stable zero mean curvature set with arbitrarily large perimeter in~$B_{1/2}$.\medskip

For flatness results about $\alpha$-stable case, see also~\cite{MR3981295, MR4116635, NHARDY}.
\medskip

To understand model cases for the regularity theory of nonlocal minimal surfaces,
to start with one can also focus specifically on cones possessing some special symmetry:
in particular, inspired by Lawson's construction~\cite{MR308905}, we may be willing to
fix our attention to cones of the form
\begin{equation}\label{LWDE}
\{(y, z) \in\R^m\times\R^{n-m} {\mbox{ s.t. }} |z|=\delta\,|y|\},\end{equation}
for some~$\delta>0$.

This case already reveals a very intriguing property, specific of the nonlocal setting. As a matter of fact,
as proved in Theorems~3 and~4 of~\cite{MR3798717}, on the one hand one can always find a unique~$\delta$
for which the cone in~\eqref{LWDE} has vanishing nonlocal mean curvature
(and this is already in contrast with the classical case, in which some dimensional restriction is needed).
On the other hand, there exists~$\alpha_\star\in(0,1)$ such that,
\begin{equation}\label{0qwoifk0OSJHDN-OKSDND8-UJHS-b-c}
{\mbox{for~$\alpha\in(0,\alpha_\star)$
these cones are unstable when~$n\le6$ and stable when~$n=7$. }}\end{equation}
That is, stability for this type
of cones occurs already in dimension~$7$ when the fractional exponent is ``small enough''
(this has to be compared instead with the different
classical scenario in~\eqref{0qwoifk0OSJHDN-OKSDND8-UJHS-b}).

The result in~\eqref{0qwoifk0OSJHDN-OKSDND8-UJHS-b-c} is extremely interesting because it reveals that the nonlocal case in dimension~$7$ must exhibit one of these remarkable situations:
\begin{itemize}
\item either the stable cones in~\eqref{0qwoifk0OSJHDN-OKSDND8-UJHS-b-c}
are minimizers (which would affirm that the restriction in~\eqref{AGbhgVjwehe-3}
on the fractional exponent to be ``large enough'' cannot be dropped in this case,
and would also establish a different regularity theory with respect to the classical case in~\eqref{0qwoifk0OSJHDN-OKSDND8-UJHS}),
\item or the stable cones in~\eqref{0qwoifk0OSJHDN-OKSDND8-UJHS-b-c}
are not minimizers (which would provide an example of stable cone which is not a minimizer,
while no example of this type is known in the classical setting).
\end{itemize}
Therefore, the result in~\eqref{0qwoifk0OSJHDN-OKSDND8-UJHS-b-c} entails that something
beautiful is happening in the nonlocal world (though we do not understand what is really happening!).
\medskip

Though we do not go into details here, we mention that the boundary regularity
for $\alpha$-minimal surfaces is also completely different from the classical case,
since in the nonlocal world these surfaces have the possibility of producing boundary jumps
(and they do so generically), according to a phenomenon that we christened
with the name of ``stickiness'': see~\cite{MR3596708, MR3824212, MR3926519, MR4104542, MR4178752, MR4392355, 2021arXiv211209299B}
for further details. See also~\cite{MR4188155}
for further comparisons between classical and nonlocal minimal surfaces.

\subsection{A nonlocal version of De Giorgi's conjecture}
The regularity of nonlocal minimal surfaces has also a solid link with the fractional counterpart
of the question posed by De Giorgi about the one-dimensional symmetry of phase transitions
(recall Conjecture~\ref{DG:CO}, and see~\cite{MR3035063} for a unified version
of these problems and methodologies). More specific, the
long-range version of Conjecture~\ref{DG:CO} can be stated as follows:

\begin{conjecture} \label{DG:CO:frac}
Let~$\alpha\in(0,2)$ and~$u\in C^2(\R^n)\cap L^\infty(\R^n)$ be a global solution of the fractional Allen-Cahn equation
\begin{equation}\label{dhffRAsdaokhHnPasdf} -(-\Delta)^{\frac\alpha2} u=u-u^3\end{equation} such that
$$ \partial_n u(x)>0\quad{\mbox{ for every }}\,x\in\R^n.$$
Is it true that~$u$ is one-dimensional (at least under suitable restrictions on~$\alpha$ and~$n$)?
\end{conjecture}

Note that the dimensional restriction in Conjecture~\ref{DG:CO:frac} could be, in principle,
different from that in Conjecture~\ref{DG:CO}, since the fractional version may be intertwined with the Bernstein property
for $\alpha$-minimal graphs, which, in its full generality, is still under investigation
(recall~\eqref{BEJDrn92o3esefrP} and~\eqref{0qwoifk0OSJHDN-OKSDND8-UJHS-b-c}).\medskip

To the best of our knowledge, the state of the art about Conjecture~\ref{DG:CO:frac} is that it is known to hold true:
\begin{itemize}
\item when~$n=2$ and~$\alpha=1$, as proved in~\cite{MR2177165},
\item when~$n=2$, for all~$\alpha\in(0,2)$, as proved in~\cite{MR2498561, MR3280032, MR3469920},
\item when~$n=3$ and~$\alpha=1$, as proved in~\cite{MR2644786},
\item when~$n=3$ and~$\alpha\in(1,2)$, as proved in~\cite{MR3148114},
\item when~$n=3$ and~$\alpha\in(0,1)$, as proved in~\cite{MR3740395},
\item when~$n=4$ and~$\alpha=1$, as proved in~\cite{MR4050103},
\end{itemize}
the other cases remaining open.\medskip

Additionally, under the additional limit assumption~\eqref{omUjPoOP1-23e},
Conjecture~\ref{DG:CO:frac} when~$n\le8$ has been established
when~$\alpha\in(1,2)$ in~\cite{MR3812860},
when~$\alpha=1$ in~\cite{MR3939768}
and when~$\alpha\in(\alpha_0,1)$ for a suitable~$\alpha_0\in(0,1)$ in~\cite{MR4124116}.\medskip

Interestingly, the case in~\cite{MR3740395} also relies on some of the results introduced in~\cite{MR4124116},
hence we now briefly comment on these results and on their connection with
Conjecture~\ref{DG:CO:frac}.

Though~\cite{MR4124116}
actually deals with a more general form of nonlocal Allen-Cahn equation, both
in terms of interaction kernels and of bistable nonlinearities allowed, we stick here to the model case
investigated in Conjecture~\ref{DG:CO:frac} for the sake of simplicity. Then, one picks~$\theta\in\left(
\frac{\sqrt{3}}3,1\right)$ and introduces the notion of ``interface'' for a function~$u$ taking values in~$[-1,1]$ as the set
$$ I(u):=\left\{x\in\R^n{\mbox{ s.t. }}|u(x)|\le\theta\right\}.$$
Note that, by the intuition arising from the phase coexistence model, $I(u)$ collects precisely the points
in the space which are not ``close'' to the pure phases (and the explicit value of~$\frac{\sqrt{3}}3$
is taken to ensure the monotonicity of the nonlinearity in~\eqref{dhffRAsdaokhHnPasdf} for~$u\in[-1,-\theta]\cup[\theta,1]$).

\begin{figure}[h]
\includegraphics[height=0.45\textwidth]{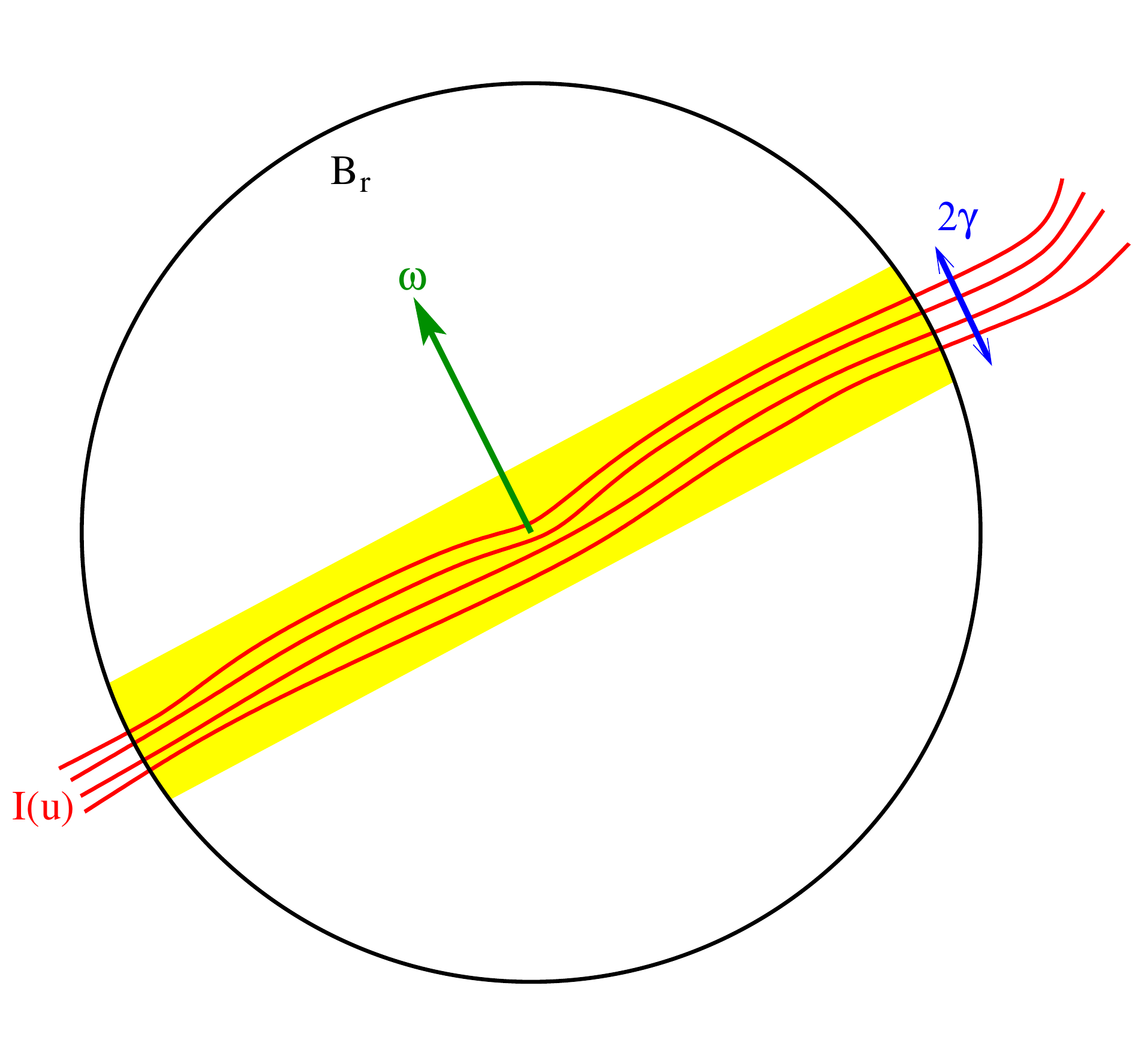}
\caption{The trapping notion for the interface.}
        \label{d19TRAUNIgIbnNP-090}
\end{figure}

Given~$r>0$, $\gamma>0$ and $\omega\in \partial B_1$,
one says that the interface~$I(u)$ is $\gamma$-trapped in~$B_r$ in direction~$\omega$
if
\begin{eqnarray*}
&& \{\omega\cdot x\le-\gamma\}\cap B_r \subseteq
\{u\le-\theta\}\cap B_r 
\\{\mbox{and }}&&
\{u\le\theta\}\cap B_r 
\subseteq
\{\omega\cdot x\le\gamma\}\cap B_r ,\end{eqnarray*}
see Figure~\ref{d19TRAUNIgIbnNP-090}.

We stress that when~$\gamma\geq r$, $I(u)$ is always
$\gamma$-trapped in~$B_r$ in direction~$\omega$, for every possible choices of~$r$ and~$\gamma$,
hence this condition becomes meaningful for~$\gamma<r$, and in fact it is the flatness ratio~$a:=\frac\gamma{r}$ which plays a role in the geometric description of the interface.
\medskip

In this setting, in Theorem~1.2 of~\cite{MR4124116}
one obtains the one-dimensional symmetry of the solution as a consequence of a suitable
``flatness at infinity'', i.e. from the sublinear growth of the interface~$I(u)$ at infinity. The precise result
is here below:

\begin{theorem}\label{NS:ocla234}
Let~$\alpha\in(0,1)$ and $u$ be a solution of the fractional Allen-Cahn equation~\eqref{dhffRAsdaokhHnPasdf}
in~$\R^n$. 

Assume that there exists~$a:(1,+\infty)\to(0,1]$
such that
$$ \lim_{r\to+\infty} a(r)=0$$
and
$I(u)$ is~$a(r)\,r$-trapped in~$B_r$ in some direction~$\omega_r$.

Then, $u$ is one-dimensional.
\end{theorem}

\begin{figure}[h]
\includegraphics[width=0.42\textwidth]{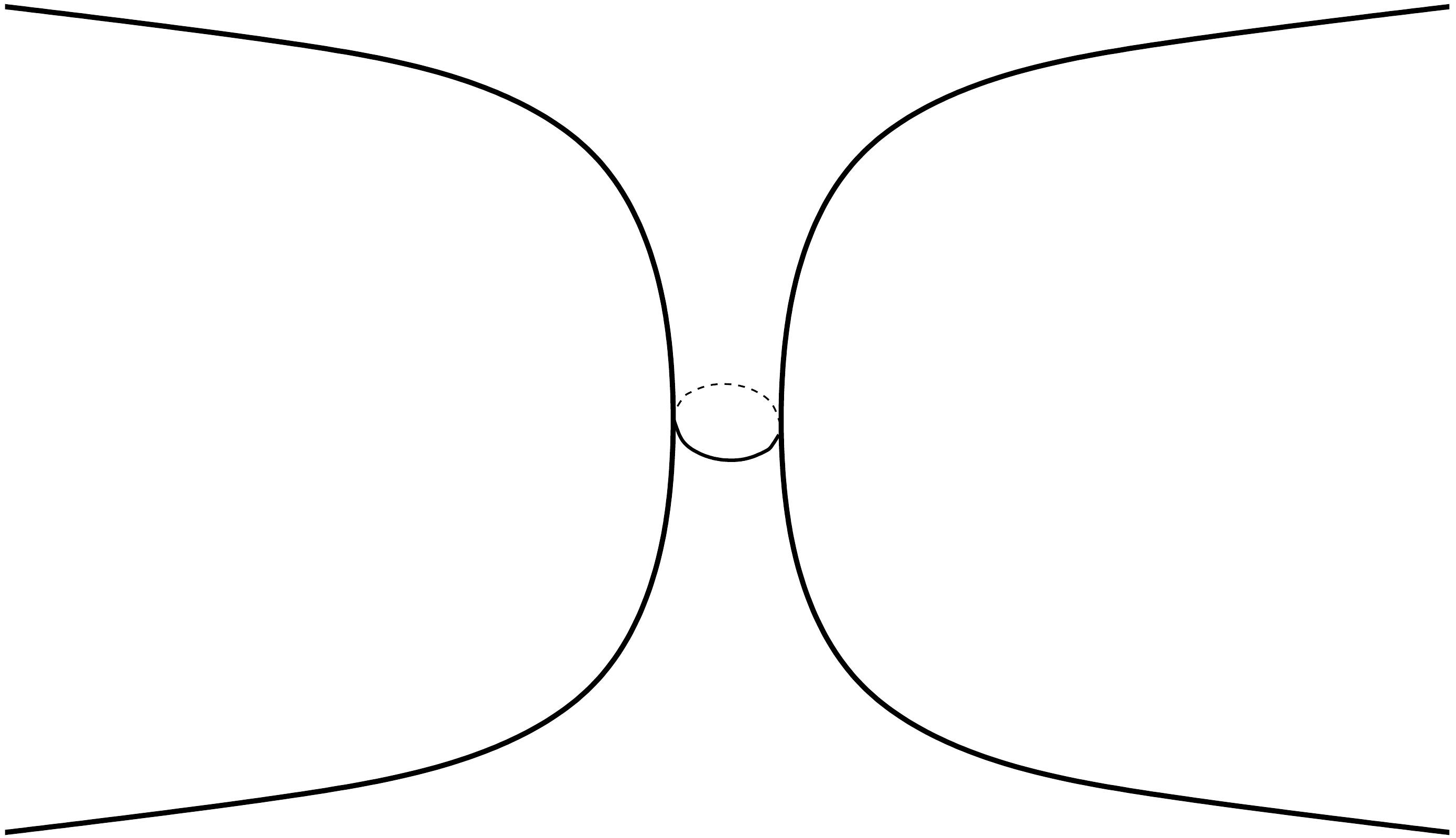}$\qquad\qquad$\includegraphics[width=0.42\textwidth]{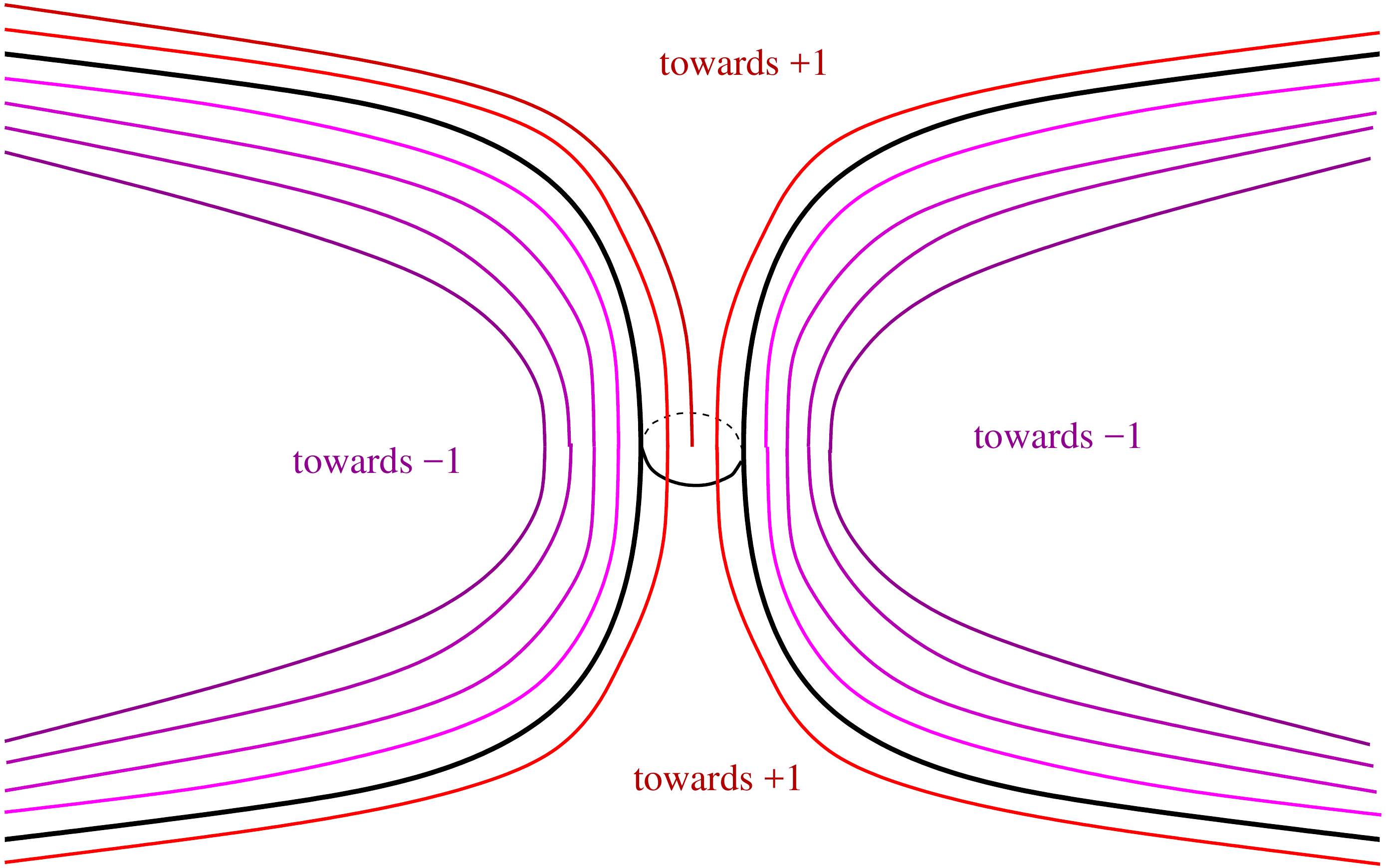}
\caption{Level sets of classical phase transition modeled on a catenoid.}
        \label{d19TRACAT090}
\end{figure}

We observe that this result is also purely nonlocal, in the sense that in the classical case
an analogous result is false. Indeed, in~\cite{MR3019512} an entire solution of the classical Allen-Cahn
equation has been constructed whose interface is modeled on a catenoid
(see Figure~\ref{d19TRACAT090} for a sketch): then, since the catenoid grows logarithmically at infinity,
the interface in this case is~$\gamma$-trapped in~$B_r$ in the vertical direction with~$\gamma\simeq\ln r=a(r)\,r$, where~$a(r):=\frac{\ln r}r\to0$ as~$r\to+\infty$ (but~$u$ is not one-dimensional,
showing that a classical counterpart of Theorem~\ref{NS:ocla234} does not hold true).\medskip

It is also worth pointing out that Theorem~\ref{NS:ocla234} is actually the consequence
of a more general trapping result (namely Theorem~1.1 of~\cite{MR4124116})
which does not require the solution to be a global one, but rather a quantitative decay of the interface
at infinity, via an improvement of flatness.\medskip

Then, the strategy adopted in~\cite{MR3740395} to prove 
Conjecture~\ref{DG:CO:frac} when~$n=3$ and~$\alpha\in(0,1)$
is to first check that in this case the solution is a local minimizer for the energy
(for this, extension, energy renormalization and sliding methods turn in handy,
combined with the observation
that the profiles at infinity are stable and two-dimensional solutions,
hence necessarily one dimensional, thanks to~\cite{MR2498561, MR3280032}).

Hence, one considers a rescaled version~$u_\e(x):=u\left(\frac{x}\e\right)$ of the original solution and employs the $\Gamma$-convergence
theory in Theorem~\ref{oqjdlwn0304-53-tgol} to deduce that~$u_\e$ approaches a step function of the type~$\chi_E-\chi_{E^c}$,
where~$E$ is a fractional minimal surface in~$\R^3$.

A Bernstein type result in the spirit of~\cite{MR3680376} then gives that~$\partial E$ is a hyperplane.
But one knows from the geometric convergence in Theorem~\ref{THM:CC:E} that actually the level sets
of~$u_\e$ must converge to this hyperplane and consequently, scaling back, one obtains
an interface satisfying the assumptions of Theorem~\ref{NS:ocla234}, thus allowing one to conclude
that~$u$ is necessarily one-dimensional.

\appendix

\section{Establishing some claims in the proof of~\eqref{0qwoifk0OSJHDN-OKSDND8-UJHS-c}}\label{appeadd10}

In this appendix we collect the proofs of some claims stated through
the proof of~\eqref{0qwoifk0OSJHDN-OKSDND8-UJHS-c}.

\begin{proof}[Proof of~\eqref{BELL7}]
We denote by~$d$ the signed distance function (positive outside~$E$) and observe that,
in the vicinity of~$\partial E$,
\begin{equation}\label{ODR-3r-2} \nu(X)-\nabla d(X)=o(d(X)),\end{equation}
for some~$C>0$, because if~$X=x+d(X)\nu(x)$ for some~$x\in\partial E$, using~\eqref{NE55} we have that
\begin{eqnarray*}
\nu_i(X)-\partial_i d(X)&=&\nu_i(x)-\partial_i d(x+d(X)\nu(x))\\&=&
\nu_i(x)-\partial_i d(x)-d(X)\nabla\partial_id(x)\cdot \nu(x)+o(d(X))\\&=&0
-d(X)\nabla\partial_id(x)\cdot \nabla d(x)+o(d(X))\\
&=&-\frac{d(X)}2\partial_i|\nabla d(x)|^2+o(d(X))\\
&=&-\frac{d(X)}2\partial_i1+o(d(X))\\&=&o(d(X)),
\end{eqnarray*}
which demonstrates~\eqref{ODR-3r-2}.

{F}rom~\eqref{ODR-3r-2} it also follows that, for all~$x\in\partial E$,
\begin{equation*}
\partial_i\nu_j(x)=\lim_{\e\to0}\frac{\nu_j(x+\e e_i)-\nu_j(x)}\e
=\lim_{\e\to0}\frac{\partial_jd(x+\e e_i)+o(\e)-\partial_jd(x)}\e=\partial_{ij}^2d(x)
\end{equation*}
and, in particular,
\begin{equation}\label{334-ipkjd-3}
\partial_i\nu_j(x)=\partial_j\nu_i(x).
\end{equation}

Besides, we recall the notion of tangential divergence of a vector field~$F$
\begin{equation}\label{diuewyt8tuy9854674967698708}
\div_T F:=\div F-\nabla(F\cdot\nu)\cdot\nu,\end{equation}
see e.g.~\cite{dipierro2021elliptic}.
We observe that
\begin{eqnarray*}&&
\sum_{i=1}^n (\partial_{T,i} f e_i)\cdot\nu=
\sum_{i=1}^n \Big( \partial_i f-(\nabla f\cdot\nu)\nu_i\Big)e_i\cdot\nu\\&&\qquad=
\sum_{i=1}^n \Big( \partial_i f-(\nabla f\cdot\nu)\nu_i\Big)\nu_i
=(\nabla f\cdot\nu)-(\nabla f\cdot\nu)\sum_{i=1}^n\nu_i^2=0
\end{eqnarray*}
and, as a result, by the Tangential Divergence Theorem (see e.g.~\cite{dipierro2021elliptic}) applied
to the vector field~$F:=\partial_{T,i} f e_i$, we have that, for all~$g\in C^1(\partial E)$,
\begin{equation}\label{NE10}\begin{split}
&\sum_{i=1}^n\int_{\partial E} \div_T \big(\partial_{T,i} f (x)\,e_i\big)\,g(x)\,d{\mathcal{H}}^{n-1}_x\\=\;&\sum_{i=1}^n\int_{\partial E} (\partial_{T,i} f (x)\,e_i\big)\cdot\Big(
H(x)\,\nu(x)\,g(x)-\nabla_T g(x)
\Big)\,d{\mathcal{H}}^{n-1}_x\\
=\;&-\sum_{i=1}^n\int_{\partial E} (\partial_{T,i} f (x)\,e_i\big)
\cdot\nabla_T g(x)\,d{\mathcal{H}}^{n-1}_x\\
=\;&-\sum_{i=1}^n\int_{\partial E} \partial_{T,i} f (x)\,\partial_{T,i} g(x)\,d{\mathcal{H}}^{n-1}_x.
\end{split}\end{equation}
Moreover, recalling~\eqref{DEDAS}, in view of~\eqref{NE56},
\begin{eqnarray*}
&&\sum_{i=1}^n \div_T (\partial_{T,i} f e_i)=\sum_{i=1}^n\Big[
\div(\partial_{T,i} f e_i)-\nabla\big((\partial_{T,i} f e_i)\cdot\nu\big)\cdot\nu\Big]\\&&\qquad=\sum_{i=1}^n\Big[
\partial_i(\partial_{T,i} f)-\nabla\big(\partial_{T,i} f \nu_i\big)\cdot\nu\Big]\\&&\qquad=\sum_{i=1}^n\Big[
\partial_{T,i}(\partial_{T,i} f)
+\big( \nabla(\partial_{T,i} f)\cdot\nu\big)\nu_i
-\nabla\big(\partial_{T,i} f \nu_i\big)\cdot\nu\Big]\\&&\qquad=\sum_{i=1}^n
\partial_{T,i}\partial_{T,i} f
+\sum_{i,j=1}^n\Big[\partial_j(\partial_{T,i} f)\nu_j\nu_i
-\partial_j\big(\partial_{T,i} f \nu_i\big)\nu_j\Big]\\&&\qquad=\sum_{i=1}^n
\partial_{T,i}\partial_{T,i} f
-\sum_{i,j=1}^n \partial_{T,i} f \,\partial_j\nu_i \,\nu_j\\&&\qquad=\sum_{i=1}^n
\partial_{T,i}\partial_{T,i} f
\end{eqnarray*}
and consequently, by~\eqref{NE10},
\begin{eqnarray*}
\sum_{i=1}^n\int_{\partial E} \partial_{T,i}\partial_{T,i} f(x)\,g(x)\,d{\mathcal{H}}^{n-1}_x&=&-\sum_{i=1}^n\int_{\partial E} \partial_{T,i} f (x)\,\partial_{T,i} g(x)\,d{\mathcal{H}}^{n-1}_x\\&=&
- \int_{\partial E} \nabla_T f (x)\cdot\nabla_T g(x)\,d{\mathcal{H}}^{n-1}_x.
\end{eqnarray*}
Owing to~\eqref{KSLABELTRA}, the proof of~\eqref{BELL7} is thereby complete.
\end{proof}

\begin{proof}[Proof of~\eqref{COMMUTATA0} and~\eqref{COMMUTATA}]
We first recall~\eqref{334-ipkjd-3} to see that
\begin{eqnarray*}
&&\partial_{T,i}\nu_j-\partial_{T,j}\nu_i=
\partial_i\nu_j-(\nabla\nu_j\cdot\nu)\nu_i
-\partial_j\nu_i+(\nabla\nu_i\cdot\nu)\nu_j=
\partial_i\nu_j-\partial_j\nu_i=0,
\end{eqnarray*}
which demonstrates~\eqref{COMMUTATA0}.

Now, to prove~\eqref{COMMUTATA}, we use~\eqref{DEDAS} and~\eqref{COMMUTATA0}
to see that
\begin{eqnarray*}&&
\partial_{T,i}\partial_{T,j}\phi\\&=&
\partial_{T,i}\Big(\partial_j\phi-(\nabla\phi\cdot\nu)\nu_j\Big)\\&=&
\partial_{T,i}\partial_j\phi-\sum_{k=1}^n\Big[ (\partial_{T,i}\partial_k\phi)\nu_k\nu_j
+\partial_k\phi \Big(\partial_{T,i}\nu_k\nu_j+\nu_k\partial_{T,i}\nu_j\Big)\Big]\\
&=&\partial_{ij}^2\phi-\sum_{k=1}^n\Big[ 
\partial^2_{jk}\phi \,\nu_k\nu_i+
(\partial_{T,i}\partial_k\phi)\nu_k\nu_j
+\partial_k\phi \Big(\partial_{T,i}\nu_k\nu_j+\nu_k\partial_{T,i}\nu_j\Big)\Big]\\&=&\partial_{ij}^2\phi-\sum_{k=1}^n\Big[ 
\partial^2_{jk}\phi \,\nu_k\nu_i+
\partial_{ik}^2\phi\,\nu_k\nu_j
+\partial_k\phi \Big(\partial_{T,i}\nu_k\nu_j+\nu_k\partial_{T,i}\nu_j\Big)\Big]\\&&\qquad\qquad\qquad
+\sum_{k,m=1}^n \partial_{mk}^2\phi\,\nu_m\nu_k\nu_j\nu_i
\\&=&
-\sum_{k=1}^n\partial_k\phi \Big(\partial_{T,i}\nu_k\nu_j\Big)
+{\mathcal{S}},
\end{eqnarray*}
where~${\mathcal{S}}$ denotes a quantity
which is symmetric in the indexes~$i$ and~$j$.

As a result,
\begin{equation}\label{CABNSM:IAKMSbbMNs-2}
\partial_{T,i}\partial_{T,j}\phi-\partial_{T,j}\partial_{T,i}\phi=
\sum_{k=1}^n\Big[
\partial_{T,j}\nu_k\nu_i-\partial_{T,i}\nu_k\nu_j\Big] \partial_k\phi.\end{equation}

Now we notice that
\begin{equation*}\begin{split}&
\sum_{k=1}^n\Big[
\partial_{T,j}\nu_k\nu_i-\partial_{T,i}\nu_k\nu_j\Big] \big(\partial_{T,k}\phi-\partial_k\phi\big)=-
\sum_{k=1}^n\Big[
\partial_{T,j}\nu_k\nu_i-\partial_{T,i}\nu_k\nu_j\Big] (\nabla\phi\cdot\nu)\nu_k=0,\end{split}
\end{equation*}
owing to~\eqref{DEDAS} and~\eqref{CABNSM:IAKMSbbMNs-1}.

By inserting this information into~\eqref{CABNSM:IAKMSbbMNs-2}, we arrive at~\eqref{COMMUTATA}, as desired.
\end{proof}

\begin{proof}[Proof of~\eqref{G1018}]
Exploiting~\eqref{BELL7} and~\eqref{COMMUTATA}, we have that
\begin{eqnarray*}&&
\partial_{T,k}(\Delta_T\phi)-\Delta_T(\partial_{T,k}\phi)\\&=&
\sum_{i=1}^n \partial_{T,k}\partial_{T,i}\partial_{T,i}\phi-\sum_{i=1}^n \partial_{T,i}\partial_{T,i}\partial_{T,k}\phi
\\&=&
\sum_{i=1}^n \partial_{T,k}\partial_{T,i}\partial_{T,i}\phi-\sum_{i=1}^n \partial_{T,i}\Bigg[
\partial_{T,k}\partial_{T,i}\phi+\sum_{m=1}^n
\Big(\nu_i\partial_{T,k}\nu_m-\nu_k\partial_{T,i}\nu_m\Big)\partial_{T,m}\phi
\Bigg]\\&=&\sum_{i=1}^n \partial_{T,k}\partial_{T,i}\partial_{T,i}\phi-\sum_{i=1}^n \partial_{T,i}
\partial_{T,k}\partial_{T,i}\phi-\sum_{i,m=1}^n \partial_{T,i}\Bigg[
\Big(\nu_i\partial_{T,k}\nu_m-\nu_k\partial_{T,i}\nu_m\Big)\partial_{T,m}\phi\Bigg]
\\&=&\sum_{i=1}^n \partial_{T,k}\partial_{T,i}\partial_{T,i}\phi-\sum_{i=1}^n\Bigg[ \partial_{T,k}
\partial_{T,i}\partial_{T,i}\phi+\sum_{m=1}^n \Big(\nu_i\partial_{T,k}\nu_m-\nu_k\partial_{T,i}\nu_m
\Big)\partial_{T,m}\partial_{T,i}\phi
\Bigg]\\&&\qquad
-\sum_{i,m=1}^n 
\Big((\partial_{T,i}\nu_i)(\partial_{T,k}\nu_m)+\nu_i\partial_{T,i}\partial_{T,k}\nu_m
-(\partial_{T,i}\nu_k)(\partial_{T,i}\nu_m)-\nu_k\partial_{T,i}\partial_{T,i}\nu_m\Big)\partial_{T,m}\phi
\\&&\qquad-\sum_{i,m=1}^n 
\Big(\nu_i\partial_{T,k}\nu_m-\nu_k\partial_{T,i}\nu_m\Big)\partial_{T,i}\partial_{T,m}\phi
\\&=&-\sum_{i,m=1}^n \Big(\nu_i\partial_{T,k}\nu_m-\nu_k\partial_{T,i}\nu_m
\Big)\partial_{T,m}\partial_{T,i}\phi
\\&&\qquad
-\sum_{i,m=1}^n 
\Big((\partial_{T,i}\nu_i)(\partial_{T,k}\nu_m)+\nu_i\partial_{T,i}\partial_{T,k}\nu_m
-(\partial_{T,i}\nu_k)(\partial_{T,i}\nu_m)-\nu_k\partial_{T,i}\partial_{T,i}\nu_m\Big)\partial_{T,m}\phi
\\&&\qquad-\sum_{i,m=1}^n 
\Big(\nu_i\partial_{T,k}\nu_m-\nu_k\partial_{T,i}\nu_m\Big)\partial_{T,i}\partial_{T,m}\phi
\\&=&-\sum_{i,m=1}^n \Big(\nu_i\partial_{T,k}\nu_m-\nu_k\partial_{T,i}\nu_m
\Big)\Bigg[\partial_{T,i}\partial_{T,m}\phi+\sum_{\ell=1}^n \Big(\nu_m\partial_{T,i}\nu_\ell-\nu_i\partial_{T,m}\nu_\ell
\Big)\partial_{T,\ell}\phi\Bigg]
\\&&\qquad
-\sum_{i,m=1}^n 
\Big((\partial_{T,i}\nu_i)(\partial_{T,k}\nu_m)+\nu_i\partial_{T,i}\partial_{T,k}\nu_m
-(\partial_{T,i}\nu_k)(\partial_{T,i}\nu_m)-\nu_k\partial_{T,i}\partial_{T,i}\nu_m\Big)\partial_{T,m}\phi
\\&&\qquad-\sum_{i,m=1}^n 
\Big(\nu_i\partial_{T,k}\nu_m-\nu_k\partial_{T,i}\nu_m\Big)\partial_{T,i}\partial_{T,m}\phi
\\&=&-\sum_{i,m,\ell=1}^n \Big(\nu_i\partial_{T,k}\nu_m-\nu_k\partial_{T,i}\nu_m
\Big) \Big(\nu_m\partial_{T,i}\nu_\ell-\nu_i\partial_{T,m}\nu_\ell
\Big)\partial_{T,\ell}\phi
\\&&\qquad
-\sum_{i,m=1}^n 
\Big((\partial_{T,i}\nu_i)(\partial_{T,k}\nu_m)+\nu_i\partial_{T,i}\partial_{T,k}\nu_m
-(\partial_{T,i}\nu_k)(\partial_{T,i}\nu_m)-\nu_k\partial_{T,i}\partial_{T,i}\nu_m\Big)\partial_{T,m}\phi
\\&&\qquad-2\sum_{i,m=1}^n 
\Big(\nu_i\partial_{T,k}\nu_m-\nu_k\partial_{T,i}\nu_m\Big)\partial_{T,i}\partial_{T,m}\phi.
\end{eqnarray*}
Now we observe that the fact that~$H=0$ gives that
$$\sum_{i=1}^n\partial_{T,i}\nu_i=0.$$
Also, by~\eqref{G1014} we deduce that
$$ \sum_{i=1}^n\nu_i\partial_{T,i}\partial_{T,k}\nu_m=0.$$
Furthermore, recalling~\eqref{G1017} and using again the fact that~$H=0$ and~\eqref{G1014},
$$ \sum_{i,m=1}^n (\partial_{T,i}\partial_{T,i}\nu_m)(\partial_{T,m}\phi)=\sum_{m=1}^n(\Delta_T\nu_m)(\partial_{T,m}\phi)
=-c^2\sum_{m=1}^n\nu_m\partial_{T,m}\phi=0
.$$
{F}rom these considerations, we conclude that
\begin{equation}\label{ao3ircbo9tu45icy56vuj785v}\begin{split}&
\partial_{T,k}(\Delta_T\phi)-\Delta_T(\partial_{T,k}\phi)\\=\;&
-\sum_{i,m,\ell=1}^n \Big(\nu_i\partial_{T,k}\nu_m-\nu_k\partial_{T,i}\nu_m
\Big) \Big(\nu_m\partial_{T,i}\nu_\ell-\nu_i\partial_{T,m}\nu_\ell
\Big)\partial_{T,\ell}\phi
\\&\qquad
+\sum_{i,m=1}^n (\partial_{T,i}\nu_k)(\partial_{T,i}\nu_m)(\partial_{T,m}\phi)
-2\sum_{i,m=1}^n 
\Big(\nu_i\partial_{T,k}\nu_m-\nu_k\partial_{T,i}\nu_m\Big)\partial_{T,i}\partial_{T,m}\phi\\=\;&
-\sum_{i,m,\ell=1}^n \Bigg[\nu_i\nu_m (\partial_{T,k}\nu_m)(\partial_{T,i}\nu_\ell)
-\nu_i^2(\partial_{T,k}\nu_m)(\partial_{T,m}\nu_\ell)\\&\qquad\qquad\qquad
-\nu_k\nu_m(\partial_{T,i}\nu_m)(\partial_{T,i}\nu_\ell)
+\nu_k\nu_i(\partial_{T,i}\nu_m)(\partial_{T,m}\nu_\ell)\Bigg]\partial_{T,\ell}\phi
\\&\qquad
+\sum_{i,m=1}^n (\partial_{T,i}\nu_k)(\partial_{T,i}\nu_m)(\partial_{T,m}\phi)
-2\sum_{i,m=1}^n 
\Big(\nu_i\partial_{T,k}\nu_m-\nu_k\partial_{T,i}\nu_m\Big)\partial_{T,i}\partial_{T,m}\phi.
\end{split}\end{equation}
Notice that
\begin{eqnarray*}&& \sum_{m=1}^n\nu_m \partial_{T,k}\nu_m=\sum_{m=1}^n\nu_m \partial_{T,m}\nu_k=0,\\
&&\sum_{m=1}^n\nu_m \partial_{T,i}\nu_m=\sum_{m=1}^n\nu_m \partial_{T,m}\nu_i=0\\
{\mbox{and }}&&\sum_{i=1}^n\nu_i\partial_{T,i}\nu_m=0,
\end{eqnarray*}
in light of~\eqref{G1014} and~\eqref{COMMUTATA0}.

Similarly,
$$ \sum_{i=1}^n\nu_i \partial_{T,i}\partial_{T,m}\phi=0.
$$
Thus, plugging these pieces of information into~\eqref{ao3ircbo9tu45icy56vuj785v}
and using~\eqref{COMMUTATA0}, we obtain that
\begin{eqnarray*}&&
\partial_{T,k}(\Delta_T\phi)-\Delta_T(\partial_{T,k}\phi)\\&=&
\sum_{i,m,\ell=1}^n \nu_i^2(\partial_{T,k}\nu_m)(\partial_{T,m}\nu_\ell)( \partial_{T,\ell}\phi)
+\sum_{i,m=1}^n (\partial_{T,i}\nu_k)(\partial_{T,i}\nu_m)(\partial_{T,m}\phi)
\\&&\qquad+2\sum_{i,m=1}^n 
\nu_k(\partial_{T,i}\nu_m)(\partial_{T,i}\partial_{T,m}\phi)\\
&=&\sum_{m,\ell=1}^n (\partial_{T,k}\nu_m)(\partial_{T,m}\nu_\ell)( \partial_{T,\ell}\phi)
+\sum_{i,m=1}^n (\partial_{T,k}\nu_i)(\partial_{T,i}\nu_m)(\partial_{T,m}\phi)
\\&&\qquad+2\sum_{i,m=1}^n 
\nu_k(\partial_{T,i}\nu_m)(\partial_{T,i}\partial_{T,m}\phi)\\&=&
2\sum_{i,m=1}^n (\partial_{T,k}\nu_i)(\partial_{T,i}\nu_m)(\partial_{T,m}\phi)
+2\sum_{i,m=1}^n 
\nu_k(\partial_{T,i}\nu_m)(\partial_{T,i}\partial_{T,m}\phi),
\end{eqnarray*}
which proves~\eqref{G1018}.
\end{proof}

\section{Proof of formula~\eqref{jyhtgrfu7y6tg5rfet8569b7c5tnxm4utxc46}}\label{sec:jyhtgrfu7y6tg5rfet8569b7c5tnxm4utxc46}

By~\eqref{GURNL-2}, we have that
$$ {\mathcal{Q}}(\Omega)=(-1,1)^{2n}\cup\big(((-1,1)^n)^c\times(-1,1)^n\big)\cup
\big((-1,1)^n\times((-1,1)^n)^c\big),$$
and consequently, disregarding multiplicative constants,
\begin{equation}\begin{split}\label{diwebutr4587685bvc7x489}
&\e^{\alpha-n}
\iint_{{\mathcal{Q}}(\Omega)} \frac{(u_\e(x)-u_\e(y))^2}{|x-y|^{n+\alpha}}\,dx\,dy\\=\;&\e^{\alpha-n}
\iint_{(-1,1)^{2n}} \frac{(u_\star(x_1/\e)-u_\star(y_1/\e))^2}{|x-y|^{n+\alpha}}\,dx\,dy\\&\qquad
+2\e^{\alpha-n}
\iint_{(-1,1)^n\times((-1,1)^n)^c} \frac{(u_\star(x_1/\e)-u_\star(y_1/\e))^2}{|x-y|^{n+\alpha}}\,dx\,dy\\=\;&
\iint_{\left(-\frac1\e,\frac1\e\right)^{2n}} \frac{(u_\star(X_1)-u_\star(Y_1))^2}{|X-Y|^{n+\alpha}}\,dX\,dY
+2
\iint_{\left(-\frac1\e,\frac1\e\right)^n\times\left(\left(-\frac1\e,\frac1\e\right)^n\right)^c} \frac{(u_\star(X_1)-u_\star(Y_1))^2}{|X-Y|^{n+\alpha}}\,dX\,dY\\ \simeq\;&
\iint_{\left(-\frac1\e,\frac1\e\right)^n\times\R^n} \frac{(u_\star(X_1)-u_\star(Y_1))^2}{|X-Y|^{n+\alpha}}\,dX\,dY
\\=\;&
\iint_{\left(-\frac1\e,\frac1\e\right)^n\times\R^n} \frac{(u_\star(X_1)-u_\star(Y_1))^2}{(|X_1-Y_1|^2+\zeta^2)^{\frac{n+\alpha}2}}\,dX\,dY_1\,d\zeta\\=\;&
\iint_{\left(-\frac1\e,\frac1\e\right)^n\times\R^n} \frac{(u_\star(X_1)-u_\star(Y_1))^2}{|X_1-Y_1|^{1+\alpha}(1+\eta^2)^{\frac{n+\alpha}2}}\,dX\,dY_1\,d\eta\\
\simeq\;&\frac1{\e^{n-1}}
\iint_{{\left(-\frac1\e,\frac1\e\right)\times\R}} \frac{(u_\star(X_1)-u_\star(Y_1))^2}{|X_1-Y_1|^{1+\alpha}}\,dX_1\,dY_1
.
\end{split}\end{equation}
We now use the odd symmetry of~$u_\star$ to write that
\begin{eqnarray*}&&
\iint_{{\left(-\frac1\e,\frac1\e\right)\times\R}} \frac{(u_\star(X_1)-u_\star(Y_1))^2}{|X_1-Y_1|^{1+\alpha}}\,dX_1\,dY_1\\&=&
\iint_{{\left(-\frac1\e,\frac1\e\right)\times(0,+\infty)}} \frac{(u_\star(X_1)-u_\star(Y_1))^2}{|X_1-Y_1|^{1+\alpha}}\,dX_1\,dY_1
+
\iint_{{\left(-\frac1\e,\frac1\e\right)\times(-\infty,0)}} \frac{(u_\star(X_1)-u_\star(Y_1))^2}{|X_1-Y_1|^{1+\alpha}}\,dX_1\,dY_1
\\&=&
\iint_{{\left(-\frac1\e,\frac1\e\right)\times(0,+\infty)}} \frac{(u_\star(X_1)-u_\star(Y_1))^2}{|X_1-Y_1|^{1+\alpha}}\,dX_1\,dY_1
+
\iint_{{\left(-\frac1\e,\frac1\e\right)\times(0,+\infty)}} \frac{(u_\star(-X_1)-u_\star(-Y_1))^2}{|X_1-Y_1|^{1+\alpha}}\,dX_1\,dY_1
\\&=&2
\iint_{{\left(-\frac1\e,\frac1\e\right)\times(0,+\infty)}} \frac{(u_\star(X_1)-u_\star(Y_1))^2}{|X_1-Y_1|^{1+\alpha}}\,dX_1\,dY_1
.\end{eqnarray*}
Using this into~\eqref{diwebutr4587685bvc7x489}, we obtain that
\begin{equation}\label{dweiout59698765434yujkjhgfdsdfghytrgwsqabhtgrfedyhtgrfed}
\e^{\alpha-n}
\iint_{{\mathcal{Q}}(\Omega)} \frac{(u_\e(x)-u_\e(y))^2}{|x-y|^{n+\alpha}}\,dx\,dy
\simeq\frac1{\e^{n-1}}
\iint_{{\left(-\frac1\e,\frac1\e\right)\times(0,+\infty)}} \frac{(u_\star(X_1)-u_\star(Y_1))^2}{|X_1-Y_1|^{1+\alpha}}\,dX_1\,dY_1
.\end{equation}

Now we notice that
\begin{eqnarray*}&&
\iint_{{\left(-\frac1\e,\frac1\e\right)\times(0,2)}\atop{\{|X_1-Y_1|\le1\}}} \frac{(u_\star(X_1)-u_\star(Y_1))^2}{|X_1-Y_1|^{1+\alpha}}\,dX_1\,dY_1=
\iint_{{\left(-1,3\right)\times(0,2)}\atop{\{|X_1-Y_1|\le1\}}} \frac{(u_\star'(\xi))^2|X_1-Y_1|^2}{|X_1-Y_1|^{1+\alpha}}\,dX_1\,dY_1\\&&\qquad\qquad=(u_\star'(\xi))^2
\iint_{{\left(-1,3\right)\times(-1,1)} } |Z|^{1-\alpha}\,dX_1\,dZ=4(u_\star'(\xi))^2\int_{-1}^1|Z|^{1-\alpha}\,dZ
=\frac{8(u_\star'(\xi))^2}{2-\alpha},
\end{eqnarray*}
for some~$\xi$ lying on the segment joining~$X_1$ and~$Y_1$.

Furthermore,
\begin{eqnarray*}&&
\iint_{{\left(-\frac1\e,\frac1\e\right)\times(0,2)}\atop{\{|X_1-Y_1|>1\}}} \frac{(u_\star(X_1)-u_\star(Y_1))^2}{|X_1-Y_1|^{1+\alpha}}\,dX_1\,dY_1=\iint_{{\left(\left(-\frac1\e,-1\right)\cup\left(3,\frac1\e\right) \right)\times(0,2)}\atop{\{|X_1-Y_1|>1\}}} \frac{(u_\star(X_1)-u_\star(Y_1))^2}{|X_1-Y_1|^{1+\alpha}}\,dX_1\,dY_1,
\end{eqnarray*}
which is equal to a constant depending on~$\alpha$ and~$u_\star$.

Plugging these pieces of information into~\eqref{dweiout59698765434yujkjhgfdsdfghytrgwsqabhtgrfedyhtgrfed},
we conclude that
\begin{eqnarray*}&&
\e^{\alpha-n}
\iint_{{\mathcal{Q}}(\Omega)} \frac{(u_\e(x)-u_\e(y))^2}{|x-y|^{n+\alpha}}\,dx\,dy\\&\simeq&
\frac1{\e^{n-1}}+\frac1{\e^{n-1}}
\iint_{(-2,2)\times(2,+\infty)}\frac{(u_\star(X_1)-1)^2}{|X_1-Y_1|^{1+\alpha}}\,dX_1\,dY_1
\\&&\qquad
+\frac1{\e^{n-1}}
\iint_{ \left\{ |X_1|\in\left(2,\frac1\e\right)\right\}\times(2,+\infty)}
\frac{(u_\star(X_1)-1)^2}{|X_1-Y_1|^{1+\alpha}}\,dX_1\,dY_1
\\&\simeq&
\frac1{\e^{n-1}}+\frac1{\e^{n-1}}
\iint_{ \left(-\frac1\e,-2\right) \times(2,+\infty)}
\frac{dX_1\,dY_1}{|X_1-Y_1|^{1+\alpha}}\\&=&\begin{dcases}
\e^{\alpha-n}&{\mbox{ if }}\alpha\in(0,1),\\
\e^{1-n}|\ln\e|&{\mbox{ if }}\alpha=1,\\
\e^{1-n}&{\mbox{ if }}\alpha\in(1,2),
\end{dcases}\end{eqnarray*}
which establishes~\eqref{jyhtgrfu7y6tg5rfet8569b7c5tnxm4utxc46}.

\section*{Acknowledgments}

The authors are members of AustMS.
It is a pleasure to thank Matteo Cozzi, David Perrella and David Pfefferl\'e for interesting conversations.
SD is supported by the Australian Research Council DECRA DE180100957 ``PDEs, free boundaries and
applications''. EV is supported by the Australian Laureate Fellowship FL190100081 ``Minimal
surfaces, free boundaries and partial differential equations''.

\vfill

\end{document}